\documentclass[]{article}

\usepackage{hyperref}
\usepackage{graphicx}	
\usepackage{subfig}	
\usepackage[shortlabels]{enumitem}	
\usepackage{booktabs}	
\usepackage{multirow}	
\usepackage{amsmath}	
\usepackage{subeqnarray}
\usepackage{cases}      
\usepackage{amsthm}     
\usepackage{amsfonts,amssymb}
\usepackage{color}
\usepackage{algorithm}   
\usepackage{algorithmic} 
\usepackage{setspace}
\usepackage{bm}
\usepackage{amssymb}
\usepackage[T1]{fontenc}

\numberwithin{equation}{section}

\newtheorem{lemma}{Lemma}[section]
\newtheorem{theorem}{Theorem}[section]
\newtheorem{definition}{Definition}[section]
\newtheorem{remark}{Remark}[section]

\newtheorem{proposition}{Proposition}[section]

\usepackage{geometry}
\title{\bf Non-convex Pose Graph Optimization in SLAM via Proximal Linearized Riemannian ADMM}
\author{Xin Chen\footnote{LMIB, School of Mathematical Sciences, Beihang University, Beijing, 100191, China, chenxin2020@buaa.edu.cn, handr@buaa.edu.cn},\;\;\;
	Chunfeng Cui\footnote{Corresponding author. LMIB, School of Mathematical Sciences, Beihang University, Beijing, 100191, China, chunfengcui@buaa.edu.cn},\;\;\;
	Deren Han\footnotemark[1],\;\;\;
	and\; 
	Liqun Qi\footnote{Department of Mathematics, School of Science, Hangzhou Dianzi University, Hangzhou 310018, China; Department of Applied Mathematics, The Hong Kong Polytechnic University, Hung Hom, Kowloon, HongKong, maqilq@polyu.edu.hk}
	}
\date{}

\begin{document}

\maketitle

\begin{abstract}
	Pose graph optimization  is a well-known technique for solving the pose-based simultaneous localization and mapping (SLAM) problem.
	In this paper, we represent the rotation and translation by a unit quaternion and a three-dimensional vector, and propose a new   model based on the von Mises-Fisher distribution. The constraints derived from the unit quaternions are spherical manifolds, and the projection onto the constraints can be calculated by normalization.
	Then a proximal linearized Riemannian alternating direction method of multipliers, denoted by PieADMM, is developed to solve the proposed model, which not only has low memory requirements, but also can update the poses in parallel. 
	Furthermore, we establish the sublinear  iteration complexity   of PieADMM for finding the stationary point of our model. 
	The efficiency of our proposed algorithm is demonstrated by numerical experiments on two synthetic and four 3D SLAM benchmark datasets.
	\par\
\end{abstract}

\section{Introduction} \label{introduction}
Simultaneous localization and mapping (SLAM) \cite{cadena2016past,smith1990estimating} is a crucial technology that allows mobile robots to navigate autonomously through
partially or fully unknown environments. It consists in the concurrent estimation of the state of a robot  with on-board sensors, and the construction of the map of the environment that the sensors are detecting. SLAM can be categorized into different types based on the categories of sensors and mapping techniques, such as visual SLAM \cite{davison2007monoslam}, laser SLAM \cite{hess2016real}, inertial SLAM \cite{weiss2011real}, etc.

The classical approaches for solving the SLAM problem can be categorized as filter-based \cite{whyte2006slam} or graph-based \cite{juric2021comparison} methods. In the first twenty years since the SLAM problem was proposed in 1986, the filter-based methods with probabilistic formulations had achieved accurate estimation. However, updating the covariance matrix is computationally expensive in large-scale problems. The graph-based methods, first introduced in 1997 by Lu and Milios \cite{lu1997globally}, {were} cheap with the growth of the graph. With the increase of the computational power, optimization algorithms for graph-based SLAM have received widespread attention, 
compared with the classical filter-based methods, such as extended Kalman filters \cite{smith1990estimating}, Rao–Blackwellized particle filters \cite{montemerlo2002fastslam}, and information filters \cite{thrun2004simultaneous}. Wilbers et al. \cite{wilbers2019comparison} have shown the graph-based localization achieves a higher accuracy than the particle filter.

Pose graph optimization (PGO) \cite{carlone2015initialization,juric2021comparison} can be modeled as a non-convex optimization problem underlying graph-based SLAM, in which it associates each pose with a vertex and each measurement with an edge of a graph, and needs to estimate a number of unknown poses from noisy relative measurements. The pose in 3D space typically consists of rotation and translation, where the rotation can be formulated using Euler angles, axis-angle ($so(3)$), special orthogonal group ($SO(3)$) or quaternion ($\mathbb{Q}$), and the translation is specified by a three-dimension vector $\bm{t}$. Besides, the overall pose can also be represented using special Euclidean group ($SE(3)$), Lie algebra ($se(3)$) or dual quaternion ($\mathbb{DQ}$). Different modeling methods will produce different constraints, such as no constraints in $se(3)$, matrix orthogonal and determinant constraints in $SE(3)$, or spherical constraints in $\mathbb{Q}$. Selecting a simple representation that is compatible with the problem structure will lead to an easier to solve and more accurate model. 

\subsection{Literature Review}
\indent In the last twenty years, a lot of models have been developed in terms of the different statistical distribution of noise and representation of poses. At the same time, many efficient optimization algorithms have also been proposed to solve these models. We list several results in Table \ref{summary_SLAM} and then give a comment.

From the perspective of the models, the statistical distribution of rotational noise is typically categorized into Gaussian or isotropic von Mises-Fisher (vMF) distribution, while translational noise is uniformly characterized as Gaussian noise. Based on maximum likelihood estimator, the Gaussian noise on $se(3)$ can directly derive an unconstrained nonlinear least square model  \cite{dellaert2006square,grisetti2009nonlinear,kaess2012isam2,kummerle2011g2o}.
Similarly, Cheng et al. \cite{cheng2016dual} established a least square model based on unit dual quaternion, and proposed a more efficient method to compute the Jacobian matrices. By eliminating two variables, their model is also unconstrained. 
Another modeling method represented the rotations by $SO(3)$, which was assumed to obey the vMF distribution and derived a model with orthogonal and determinant constraints \cite{carlone2015lagrangian,fan2019generalized,rosen2019se}. 
Since $se(3)$ needs to be transformed to describe the process of motion,
the expression of the objective function modeled by $SO(3)$ or $\mathbb{Q}$ and a three-dimension vector in \cite{carlone2015lagrangian,fan2019generalized,rosen2019se} is more concise compared with the unconstrained model in \cite{dellaert2006square,grisetti2009nonlinear,kaess2012isam2,kummerle2011g2o}; however, the incorporation of constraints introduces challenges.

From the perspective of the algorithms, several   efficient and accurate methods are proposed for solving large-scale problems in SLAM. The first-order optimization methods such as stochastic gradient descent \cite{grisetti2009nonlinear,olson2006fast} can reduce the complexity of gradient calculation and solve the unconstrained optimization problem effectively.  
The  algorithms with faster convergence rate, such as Gauss-Newton method \cite{dellaert2006square}, Levenberg–Marquardt method \cite{kummerle2011g2o}, trust-region method \cite{rosen2012incremental}, had also been introduced to solve this problem. 
Instead of computing the matrix inverse,  \cite{dellaert2006square,kaess2012isam2,kaess2008isam} used the matrix factorization techniques, such as QR or Cholesky factorization, to reduce the complexity, and   proposed an incremental version. 
Grisetti et al. \cite{grisetti2010tutorial} and Wagner et al. \cite{wagner2011rapid} proposed the manifold-based Gauss-Newton algorithms, in which the Jacobian matrices had a sparse structure and the update process avoided the expensive storage of large-scale systems of linear equations.

\begin{table*}[tbp]
	\centering
	\footnotesize
	\renewcommand{\arraystretch}{1.4}
	\caption{Related works for solving the SLAM problem. The algorithms SGD, G-N, L-M, TRM, SDR and GPM represent stochastic gradient descent, Gauss-Newton, Levenberg–Marquardt, Trust-Region methods, semi-definite
		relaxation and generalized proximal methods, respectively. The ``\textnormal{m}" in algorithms indicates manifold optimization methods. }
	\resizebox{\linewidth}{!}{
		\begin{tabular}{cccccccc}
			\toprule[1pt]
			Paper & Space & Variables & Distribution & Initialization & Algorithms & Convergence & Complexity 
			\\
			\midrule[1pt]
			Olson et al. \cite{olson2006fast}
			& 2D
			& $SE(2)$  
			& Gaussian     
			& Odometry               
			& SGD          
			& -             
			& -           
			\\
			Grisetti et al. \cite{grisetti2009nonlinear}
			& 3D
			& $se(3)$
			& Gaussian          
			& Odometry               
			& SGD           
			& -            
			& -           
			\\
			Dellaert and Kaess \cite{dellaert2006square}
			& 3D
			& $se(3)$
			& Gaussian          
			& Odometry               
			& SAM        
			& -            
			& -        
			\\
			Kaess et al. \cite{kaess2008isam,kaess2012isam2}
			& 3D
			& $se(3)$
			& Gaussian          
			& Odometry               
			& iSAM          
			& -            
			& -        
			\\
			Grisetti et al. \cite{grisetti2010tutorial}
			& 3D
			& $\mathbb{AU}$  
			& Gaussian       
			& Odometry              
			& mG-N        
			& -            
			& -        
			\\
			Rosen et al. \cite{rosen2012incremental}
		    & 3D
		    & $se(3)$
		    & Gaussian         
		    & Odometry              
		    & TRM          
		    & -            
		    & -    
			\\
			Wagner et al. \cite{wagner2011rapid}
			& 3D
			& $SE(3)$
			& Gaussian 
			& Odometry  
			& mG-N \& mL-M    
			& -         
			& -      
			\\
			K\"{u}mmerle et al. \cite{kummerle2011g2o}
			& 3D
			& $se(3)$
			& Gaussian         
			& Odometry               
			& L-M    
			& -            
			& -        
			\\
			Cheng et al. \cite{cheng2016dual}
			& 3D
			& $\mathbb{DQ}$
			& Gaussian         
			& Odometry               
			& G-N    
			& -            
			& -        
			\\
			Liu et al. \cite{liu2012sdr}
			& 2D
			& $SE(2)$
			& Gaussian         
			& Odometry               
			& SDR    
			& -            
			& -        
			\\
			Rosen et al. \cite{rosen2015convex}
			& 3D
			& $SE(3)$
			& Gaussian        
			& Convex Relaxation                
			& L-M    
			& -            
			& -        
			\\
			Carlone et al. \cite{carlone2015lagrangian}
			& 3D
			& $SE(3)$
			& Gaussian+vMF         
			& Chord               
			& Dual method   
			& -            
			& -        
			\\
			Rosen et al. \cite{rosen2019se}
			& 3D
			& $SE(3)$
			& Gaussian+vMF         
			& Chord               
			& SE-Sync   
			& \checkmark            
			& -        
			\\
			Fan and Murphey \cite{fan2019generalized}
			& 3D
			& $SE(3)$
			& Gaussian+vMF         
			& Chord               
			& GPM   
			& \checkmark            
			& -        
			\\
			This paper
			& 3D
			& $\mathbb{AU}$
			& Gaussian+vMF         
			& Chord               
			& PieADMM    
			& \checkmark            
			& \checkmark       
			\\
			\bottomrule[1pt]  
	\end{tabular}}
	\label{summary_SLAM}
\end{table*}

However,  second-order algorithms have fast convergence rate only  in the local region, and may  return a local minima for non-convex problems. Later work had focused on finding a better initial point and confirming the optimality of the solution.
Rosen et al.~\cite{rosen2012incremental} presented a robust incremental least-squares estimation based on Powell’s Dog-Leg trust-region method and improved the numerical stability. Carlone et al.~\cite{carlone2015lagrangian} derived a quadratic constrained quadratic programming and verified the optimal solution by checking the dual gap. In \cite{rosen2015convex},  a convex relaxation is proposed by expanding the feasible set to its convex closure, which can effectively overcome the difficulty of initial point selection of non-convex problems. Furthermore, Rosen et al.~\cite{rosen2019se} relaxed the model into a semidefinite program, and proved that the minimizer of its relaxation provides an exact maximum-likelihood estimate so long as the noise falls below a certain critical threshold. Fan and Murphey \cite{fan2019generalized} proposed an upper bound of PGO, and solved it by the generalized proximal method which can converge to the first-order critical points and do not rely on Riemannian gradients. 
Another approach to find a better local minima or global minima depended on initialization techniques \cite{carlone2015initialization,Hartley2012RA}. They pointed that the non-convex rotation estimation were the actual reason why SLAM is a difficult problem, and the translations had a minor influence on the rotation estimate. Therefore, computing a good rotation estimation will improve the performance of algorithms.

\subsection{Contributions}
\indent In this paper, based on augmented unit quaternion \cite{qi2023Augmented} which uses the unit quaternions and three-dimension vectors to represent the rotations and translations, respectively, we propose a new non-convex pose graph optimization model  under von Mises-Fisher distribution \cite{sra2012short}:
\begin{equation*}
	\begin{aligned}
		\underset{\{\tilde{q}_i\},\{\bm{t}_{i}\}}{\min}	&~~
		\sum_{(i,j)\in \mathcal{E}} 
		\|\tilde{t}_j - \tilde{t}_i-\tilde{q}_{i}\tilde{t}_{ij}\tilde{q}_{i}^{*}\|_{\Sigma_{1}}^2
		+\|\tilde{q}_{j}^{*}\tilde{q}_{i}\tilde{q}_{ij}-1\|_{\Sigma_{2}}^2, \\
		\operatorname{s.t.} ~~& \quad 
		\tilde{q}_i \in \mathbb{U},~ \bm{t}_{i} \in \mathbb{R}^3,~ i=1,2,\dots,n.
	\end{aligned}
\end{equation*} 
Here, the constraint $\mathbb{U}$ is a unit sphere.
Compared with the special orthogonal group constraints, the projection of unit quaternions can be calculated by simple normalization without singular value decomposition. In addition, since a unit quaternion in $\mathbb{U}$ requires storing four elements in the computation process, while a rotation matrix in $SO(3)$ requires storing nine elements, the quaternion representation has lower data storage.
The established model constitutes a nonlinear least square problem with a quartic objective function. The non-convex structures and spherical manifold constraints may make the model difficult to solve. By introducing redundant variables, it can be reformulated into a multi-linear least square problem:
\begin{equation*}
	\begin{aligned}
		\underset{\{\tilde{p}_i\},\{\tilde{q}_i\}, \{\bm{t}_{i}\}}{\min}	&
		\sum_{(i,j)\in \mathcal{E}} 
		\|\tilde{t}_j - \tilde{t}_i-\tilde{q}_{i}\tilde{t}_{ij}\tilde{p}_{i}^{*}\|_{\Sigma_{1}}^2
		+\|\tilde{p}_{j}^{*}\tilde{q}_{i}\tilde{q}_{ij}-1\|_{\Sigma_{2}}^2, \\
		\operatorname{s.t.} ~~~~& 
		\tilde{p}_i=\tilde{q}_i,~
		\tilde{p}_i \in \mathbb{U},~
		\tilde{q}_i \in \mathbb{R}^{4},~ 
		\bm{t}_{i} \in \mathbb{R}^3,~ 
		i=1,\dots,n.
	\end{aligned}
\end{equation*}
The variables $\tilde{p}_{i}$ and $\tilde{q}_{i}$ in the objective function are coupled, while each individual variable retains linearity and convexity. 
One valid method for solving this category of problems is splitting algorithms \cite{Boyd2011,han2022survey}, and they have been proven effective in fields such as compressive sensing, robust principal component analysis and tensor decomposition.

 Alternating direction method of multipliers (ADMM) is an effective splitting algorithm for solving large-scale problems, which has been applied in many fields, such as the image alignment problem \cite{peng2012rasl}, the robust principal component analysis model \cite{tao2011recovering}, phase retrieval \cite{wen2012alternating} and
background/foreground extraction problem \cite{yang2015alternating}. The classic ADMM is proposed by Glowinski and
Marrocco \cite{glowinski1975admm} and Gabay and Mercier \cite{gabay1976admm} for solving the linearly
constrained convex optimization problem with two blocks of variables. Furthermore, ADMM had also been extended to multi-block \cite{cai2017convergence,chen2016admm,lin2015global}, non-convex \cite{guo2017convergence,themelis2020douglas}, and nonseparable \cite{hong2016convergence,chen2019extended,cui2016convergence} cases. 
Hong et al.~\cite{hong2016convergence} discussed the  linearly constrained consensus and sharing problems which had the following nonseparable structure:  
\begin{equation*}
	\begin{aligned}
		\underset{x_1,x_2,\dots,x_{n}}{\min} ~&~ H(x_{1},x_{2},\dots,x_{n})+\sum_{i=1}^{n}\theta_{i}(x_{i})
		\\
		\operatorname{s.t.} ~\quad&~ \sum_{i=1}^{n}A_{i}x_{i}=b , \quad x_i \in \mathcal{X}_{i}, ~i=1,2,\dots,n.
	\end{aligned}
\end{equation*}
 Although the variables in $H(\cdot)$ were coupled, the experiment results shown that ADMM still had a good performance when the subproblems had closed-form solutions. In theory, the iterative sequence can converge to a critical point of the augmented Lagrangian function. 

 In recent years, many scholars have devoted themselves to multi-block linear equality constrained problem with Riemannian submanifold constraints. 
There is no guarantee that the manifold constraints can be strictly satisfied for the traditional methods in Euclidean space. 
In this case, it is usually necessary to project the last iteration point on the constraint sets.
However, the manifold optimization algorithm can ensure the feasibility of the iterative sequence on the manifold. Zhang et al.~\cite{zhang2020primal} considered the gradient-based Riemannian ADMM and gave an iteration complexity. Li et al.\cite{li2022riemannian} discussed the nonsmoooth problem solving by Moreau envelope smoothing technique. We refer readers to \cite{lu2018nonconvex,chen2024nonsmooth} for more details.

We propose a proximal linearized Riemannian alternating direction method of multipliers (PieADMM) for the non-convex pose graph optimization, which updates the other variables using the most recent partial information.  Our subproblems not only have closed-form solutions, but also can be computed in parallel, which results in a low time complexity per update. The above  superiority is verified in the large-scale numerical experiments. Theoretically, the convergence analysis is established to complement our findings.

Now we summarize our contributions in this paper as follows:
\begin{enumerate}
	\item[(i)] We propose a non-convex pose graph optimization model based on augmented unit quaternion and vMF distribution, in which the data storage is low-cost and the projection of unit quaternions can be calculated by normalization.
	
	\item[(ii)] We propose a PieADMM
	which has closed-form solutions in its subproblems, and update them in parallel.
	
	\item[(iii)] Based on the first-order optimality conditions on manifolds, we define an $\epsilon$-stationary solution of our model.
	Then, we establish the iteration complexity of $O(1/\epsilon^{2})$ of PieADMM for
	finding an $\epsilon$-stationary solution.
	
	\item[(iv)] We test our algorithm on two synthetic datasets with different data scales  and four 3D SLAM benchmark datasets. Numerical experiments verify the effectiveness of our method.
\end{enumerate}

The remaining parts of this paper are organized as follows.  In the next section, we review some basic properties of quaternions and Riemannian submanifolds. In Sect.~\ref{sec-AUQ model}, we present an augmented unit quaternion model which is constrained on the sphere manifolds. The PieADMM is proposed in Sect.~\ref{sec-proximal linearized ADMM for PGO Model}. In Sect.~\ref{sec-Convergence Analysis}, we establish the iteration complexity of $O(1/\epsilon^{2})$ of PieADMM for
finding an $\epsilon$-stationary solution. Finally, in Sect.~\ref{sec-numerical experiments}, we present the numerical results.

\section{Notation and preliminaries} \label{preliminaries}
In this section, we introduce some basic notations, definitions and lemmas for this paper.

The fields of real numbers, quaternion numbers and unit quaternion numbers are denoted by $\mathbb{R}$, $\mathbb{Q}$ and $\mathbb{U}$, respectively. 
Throughout this paper, scalars, vectors, matrices, and quaternions are denoted by lowercase letters
(e.g., $x$), boldface lowercase letters (e.g., $\bm{x}$), boldface capital letters (e.g., $X$), and lowercase letters with tilde (e.g., $\tilde{x}$), respectively.
The special orthogonal group $SO(3)$ is the set of three-dimensional rotations which is formally defined by
$SO(3):=\{R \in \mathbb{R}^{3 \times 3} \mid R^\top R = I_{3}, \operatorname{det}(R)=1\}$, and the special Euclidean group is the set of poses defined by $SE(3)\triangleq\{(R,\bm{t}): R\in SO(3),\bm{t}\in\mathbb{R}^{3}\}$.

The notation $\left\| \cdot\right\| $ denotes the $2$-norm of vectors or the Frobenius norm of matrices. Let $M$ be a positive definite linear operator; we use $\|\bm{x}\|_{M}:=\sqrt{\left\langle \bm{x},M\bm{x}\right\rangle }$ to denote its $M$-norm; and $\sigma_{\min}(M)$ and $\sigma_{\max}(M)$  denote the smallest and largest eigenvalue
of $M$, respectively. For symmetric matrices $M_{1}, M_{2} \in \mathbb{R}^{n \times n}$,
$M_{1}\succ M_{2}$ and $M_{1}\succeq M_{2}$ means that $M_{1}- M_{2}$ is positive definite and positive semidefinite, respectively.

The followings are some preliminaries about quaternions (Section \ref{sec-Quaternion}), manifolds (Section \ref{sec-Manifolds}) and optimization theory over manifolds (Section \ref{sec-Optimization Theory over Manifolds}), respectively.

\subsection{Quaternion and pose} \label{sec-Quaternion}
A quaternion number $\tilde{q} \in \mathbb{Q}$, proposed by Hamilton, has the form $\tilde{q}=q_0+q_1 \mathbf{i}+q_2 \mathbf{j}+q_3 \mathbf{k}$, where $q_0, q_1, q_2, q_3 \in \mathbb{R}$ and $\mathbf{i}, \mathbf{j}, \mathbf{k}$ are three imaginary units.
We may also write $\tilde{q}=[q_0,q_1,q_2,q_3]=[q_0, \bm{q}] \in \mathbb{R}^4$ as the vector representation where $\bm{q}=[q_1,q_2,q_3]\in \mathbb{R}^3$ for convenience. We note that we also regard the above representation as a column
vector and its transpose $[q_0, \bm{q}]^\top$ a row vector.
The sum of $\tilde{p}$ and $\tilde{q}$ is defined as
$
\tilde{p}+\tilde{q}=[p_0+q_0,\bm{p}+\bm{q}].
$
The product of $\tilde{p}$ and $\tilde{q}$ is defined by
$$
\tilde{p} \tilde{q}=\left[p_0 q_0-\bm{p} \cdot \bm{q}, p_0 \bm{q}+q_0 \bm{p}+\bm{p} \times \bm{q}\right],
$$
where $\bm{p} \cdot \bm{q}$ is the $\operatorname{dot}$ product, and $\bm{p} \times \bm{q}$ is the cross product of $\bm{p}$ and $\bm{q}$. Thus, in general, $\tilde{p} \tilde{q} \neq \tilde{q} \tilde{p}$, and we have $\tilde{p} \tilde{q}=\tilde{q} \tilde{p}$ if and only if $\bm{p} \times \bm{q}=\vec{0}$, i.e., either $\bm{p}=\vec{0}$ or $\bm{q}=\vec{0}$, or $\bm{p}=\alpha \bm{q}$ for some real number $\alpha$ (see \cite{qi2023Augmented}). The multiplication of quaternions is associative and distributive over vector addition, but is not commutative.

The conjugate of $\tilde{q}$ is the quaternion $\tilde{q}^*=q_0-q_1 \mathbf{i}-q_2 \mathbf{j}-q_3 \mathbf{k}$. Then, $(\tilde{p} \tilde{q})^*=\tilde{q}^* \tilde{p}^*$ for any $\tilde{p}, \tilde{q} \in \mathbb{Q}$.  The magnitude of $\tilde{q}$ is defined by 
$
\|\tilde{q}\|=\sqrt{\tilde{q} \tilde{q}^*}=\sqrt{\tilde{q}^* \tilde{q}}.
$
And $\tilde{q}$ is invertible if and only if $\|\tilde{q}\|$ is positive. In this case, we have
$
\tilde{q}^{-1}=\tilde{q}^{*}/\|\tilde{q}\|.
$

The quaternion $\tilde{q} \in \mathbb{Q}$ is called a unit quaternion if 
$$
\|\tilde{q}\|=\sqrt{q_0^2+q_1^2+q_2^2+q_3^2}=1.
$$ 
Denote the set of all unit quaternions by $\mathbb{U}$, which can be regarded as a unit sphere in $\mathbb{R}^{4}$.
Equivalently, a unit quaternion has the following form:  
$$
\tilde{q}=[\cos(\theta/2),\sin(\theta/2)\bm{n}],
$$
where $\bm{n}=(n_{x},n_{y},n_{z})$ is a unit vector and $\theta$ is an angle. Then, we will show that $\theta\bm{n}\in so(3)$, $\tilde{q}\in\mathbb{U}$ and $R\in SO(3)$ can all represent rotation and realize motion . 
Let a vector $\bm{t}_{1}\in \mathbb{R}^{3}$ rotates $\theta$ radians around axis $\bm{n}$ to reach $\bm{t}_{2}\in \mathbb{R}^{3}$. This process can be represented by quaternion as $$
[0,\bm{t}_{2}]=\tilde{q}[0,\bm{t}_{1}]\tilde{q}^{*}.
$$ 
Using rotation matrix in $SO(3)$, we also have 
$\bm{t}_{2}=R\bm{t}_{1}$, where 
$$
R =\cos(\theta)I_{3}+(1-\cos(\theta))\bm{n}\bm{n}^{\top}+\sin(\theta)\bm{n}^{\land}, \quad
\text{and} \quad
\bm{n}^\land=\left(\begin{array}{ccc}
	0 & -n_z & n_y  \\
	n_z & 0 & -n_x\\
	-n_y & n_x & 0 
\end{array}\right).
$$
The relationship between rotation matrix and unit quaternion is given in the next lemma.
\begin{lemma}
	\cite{kuipers1999quaternions} Given a unit quaternion $\tilde{q}=[q_0,\bm{q}]=[q_0,q_1,q_2,q_3] \in \mathbb{U}$ and
	a vector $\bm{t} \in \mathbb{R}^{3}$. Then $[0,R\bm{t}]=\tilde{q}[0,\bm{t}]\tilde{q}^{*}$, where the rotation matrix $R \in SO(3)$ satisfies
	\begin{equation*}
		R=\bm{q}\bm{q}^\top+q_0^{2}I+2q_0\bm{q}^\land+(\bm{q}^\land)^2.
	\end{equation*} 
	\label{lem-SO3-quaternion}
\end{lemma}
\vspace{-1em}
If the rotation matrix $R$ is compound motion of two rotations, i.e., $R=R_{2}R_{1}$, then the corresponding quaternion $\tilde{q}$ can be formulated as $\tilde{q}=\tilde{q}_{2}\tilde{q}_{1}$. 
Next, we show a lemma which can simplify the product of two quaternions by multiplication between matrix and vector.
Given any $\tilde{a}=[a_0,a_1,a_2,a_3]$, we define 
$$
M(\tilde{a})=\left(\begin{array}{rrrr}
	a_0 & -a_1 & -a_2 & -a_3 \\
	a_1 & a_0 & -a_3 & a_2 \\
	a_2 & a_3 & a_0 & -a_1 \\
	a_3 & -a_2 & a_1 & a_0
\end{array}\right), \qquad 
W(\tilde{a})=\left(\begin{array}{rrrr}
	a_0 & -a_1 & -a_2 & -a_3 \\
	a_1 & a_0 & a_3 & -a_2 \\
	a_2 & -a_3 & a_0 & a_1 \\
	a_3 & a_2 & -a_1 & a_0
\end{array}\right) .
$$
\begin{lemma}
	\cite{chen2024regularization} For any $\tilde{a} = [a_0,a_1,a_2,a_3]\in \mathbb{Q}$ and $\tilde{b} = [b_0,b_1,b_2,b_3]\in \mathbb{Q}$, the following statements hold
	\begin{enumerate}
		\item[(a).] $M(\tilde{a}^*)=M(\tilde{a})^\top$, $W(\tilde{a}^*)=W(\tilde{a})^\top$.
		\item[(b).] $\tilde{a}\tilde{b}=M(\tilde{a})\tilde{b}=W(\tilde{b})\tilde{a}$.
		\item[(c).] $M(\tilde{a})^\top M(\tilde{a})=W(\tilde{a})^\top W(\tilde{a})=\|\tilde{a}\|^2 I_{4}$, where $I_{4}$ is the identity matrix of size $4 \times 4$.
	\end{enumerate}
	\label{lem-MW}
\end{lemma}
A quaternion $\tilde{q}=[0,q_1,q_2,q_3] \in \mathbb{Q}$ is called a vector quaternion which has the following properties.
\begin{lemma}
	\cite{qi2023Augmented} Given a quaternion $\tilde{q} \in \mathbb{Q}$, the following statements hold
	\begin{enumerate}
		\item[(a).] $\tilde{q}$ is a vector quaternion if and only if $\tilde{q}=-\tilde{q}^{*}$.
		\item[(b).] if $\tilde{q}$ is a vector quaternion, then $\tilde{p}\tilde{q}\tilde{p}^{*}$ is still a vector quaternion for any quaternion $\tilde{p} \in \mathbb{Q}$.
	\end{enumerate}
	\label{lem-vectorquaternion}
\end{lemma}

Following \cite{qi2023Augmented}, we denote $\mathbb{AU} \triangleq\{(\tilde{q},\bm{t}):\tilde{q}\in \mathbb{U},\bm{t}\in \mathbb{R}^{3}\}$ as augmented unit quaternions.
An augmented unit quaternion vector $(\tilde{\bm{p}},\bm{t})=(\tilde{p}_{1},\bm{t_1},\tilde{p}_{2},\bm{t_2},\dots,\tilde{p}_{n},\bm{t_n}) \in \mathbb{AU}^{n}$ is an $n$-component vector, such that each component $(\tilde{p}_{i},\bm{t_i})$ is an augmented unit quaternion.
In fact, $\mathbb{AU}^{n}$ is consist of the Cartesian product of $n$ spheres and $3n$-dimensional Euclidean space $\mathbb{R}^{3n}$, which is embedded in Euclidean space $\mathbb{R}^{7n}$. As will be mentioned in the Section \ref{sec-Manifolds}, the set $\mathbb{AU}^{n}$ is a Riemannian submanifold. An augmented unit quaternion optimization problem can be formulated as
\begin{equation*}
	\min f(\tilde{\bm{p}},\bm{t}) \quad \text{s.t.}~(\tilde{\bm{p}},\bm{t}) \in \mathbb{AU}^{n},
\end{equation*}
which is also a $7n$-dimensional equality constrained optimization problem, with $7n$ real variables and $n$ spherical equality constraints.

\subsection{Manifolds}
	\label{sec-Manifolds}
	Suppose $\mathcal{M}$ is a differentiable manifold, then for any $x \in \mathcal{M}$, there exists a chart $(U, \varphi)$ in which $U$ is an open set with $x\in U \subseteq \mathcal{M} $ and $\phi$ is  a homeomorphism between $U$ and an open set $\varphi(U)$ in Euclidean space. This coordinate transform enables us to locally treat a manifold as a Euclidean space. Next, we show the definition of tangent space which can help us get a linearized approximation around a point.
	
	\begin{definition}(Tangent Space)  
		Let $\mathcal{M}$ be a subset of a linear space, For all $x \in \mathcal{M}$, the tangent space is defined by
		\begin{equation*}
			T_{x}\mathcal{M}=
			\left\lbrace \gamma'(0) \mid \gamma:I \rightarrow \mathcal{M} ~\text{is smooth and} ~\gamma(0)=x
			\right\rbrace ,
		\end{equation*}
		where $I$ is any open interval containing zeros. That is, $v \in T_{x}\mathcal{M}$ if and only if there exists a smooth curve on $\mathcal{M}$ passing through $x$ with velocity $v$.
		\label{def-Tangent Space}
	\end{definition}
	
	Define the set of all functions differentiable at point $x$ to be $\mathcal{F}_{x}$. An alternative but more general way of defining tangent space is by viewing a tangent vector $v \in T_{x}\mathcal{M}$ as an operator mapping $f \in \mathcal{F}_{x}$ to $v[f]\in \mathbb{R}$, which satisfies $v[f]=\frac{d(f(\gamma(t)))}{dt}\big|_{t=0}$. In other words, $v[f]$ computes the directional derivative of $f$ at $x$ along direction $v$. 
	The tangent bundle of a manifold $\mathcal{M}$ is the disjoint union of
	the tangent spaces of $\mathcal{M}$, i.e., 
	$$T\mathcal{M} = \{(x,v):x\in \mathcal{M} ~\text{and}~ v \in T_{x}\mathcal{M}\}.$$ 
	
	\begin{definition}(Differential) 
		The differential of $F: \mathcal{M}_{1} \rightarrow \mathcal{M}_{2}$ at the point $x \in \mathcal{M}_{1} $ is the linear map $DF(x): T_{x}\mathcal{M}_{1} \rightarrow T_{F(x)}\mathcal{M}_{2}$ defined by:
		$$
		DF(x)[v]=\frac{d}{dt}F(\gamma(t))\bigg|_{t=0} = (F \circ \gamma)'(0),
		$$
		where $\gamma$ is a smooth curve on $\mathcal{M}_{2}$ passing through $x$ at $t = 0$ with velocity $v$, and it satisfies
		$$
		(DF(x)[v])[f]=v[f \circ F], ~\text{for all} ~v \in T_{x}\mathcal{M}_{1}, ~\text{and}~ \forall f \in \mathcal{F}_{F(x)}.
		$$
		\label{def-Differential}
	\end{definition}
	\vspace{-1.5em}
	\begin{definition}(Riemannian Manifold) 
		A Riemannian manifold is a manifold which equip each tangent space of itself with a Riemannian metric $\left\langle \cdot,\cdot \right\rangle_{x} $ on $\mathcal{M}$ , and the metric varies smoothly with $x$.
		\label{def-Riemannian Manifold}
	\end{definition}
	
	When $\mathcal{M}$ is an embedded submanifolds of a linear space, $T_{x}\mathcal{M}$ coincides with the linear subspace. Furthermore, let $\left\langle \cdot,\cdot \right\rangle $ be the Euclidean metric on the linear space, the metric on $\mathcal{M}$ defined at each $x$ by restriction, $\left\langle u,v \right\rangle_{x} =\left\langle u,v \right\rangle $ for $u,v \in T_{x}\mathcal{M}$,  is a Riemannian metric. At this point, we call $\mathcal{M}$ a Riemannian submanifold. 
	Let $\mathcal{M}_{1}$, $\mathcal{M}_{2}$ be two Riemannian submanifolds, $\mathcal{M}_{1} \times \mathcal{M}_{2}$ is also a Riemannian submanifold with  tangent spaces given by 
	$$
	T_{(x_{1},x_{2})}(\mathcal{M}_{1} \times \mathcal{M}_{2})
	=T_{x_{1}}\mathcal{M}_{1} \times T_{x_{2}}\mathcal{M}_{2}.
	$$
	
	\begin{definition}(Riemannian Gradient) 
		Let $f \in \mathcal{F}_{x}$. The
		Riemannian gradient of $f$ is the vector field $\operatorname{grad}f$ on $\mathcal{M}$ uniquely defined by the following identities:
		$$
		\forall (x,v) \in T\mathcal{M}, \quad Df(x)[v]
		=\left\langle v,\operatorname{grad} f(x) \right\rangle_{x}. 
		$$
		\label{def-Riemannian Gradient}
	\end{definition}
	\vspace{-2em}
	For an $m$-dimensional Riemannian submanifold $\mathcal{M}$, by defining $v = \gamma'(0)$ and 
	$v[f]=\left\langle \gamma'(0) , \nabla f(x) \right\rangle  $, we have 
	$$
	\operatorname{grad} f(x) = \operatorname{Proj}_{T_{x}\mathcal{M}}(\nabla f(x)),
	$$
	where $\operatorname{Proj}_{T_{x}\mathcal{M}}$ is the Euclidean projection operator onto the subspace $T_{x}\mathcal{M}$. Let $(U, \varphi)$ be a chart at $x\in \mathcal{M}$, then there exists a set of basis vectors $\{\bm{a}_{i}\}_{i=1}^{m}$ of $T_{x}\mathcal{M}$, such that for $v=\sum_{i=1}^{m} v_{i}\bm{a}_{i}$, we have 
	$$
	\hat{v} = D\varphi(x)[v]=(v_{1},\dots,v_{m}), \quad \hat{f}=f\circ \varphi, \quad \hat{x}=f(x),
	$$
	where $\hat{o}$ denote the Euclidean counterpart of an object $o$ in $\mathcal{M}$. If we define the Gram matrix $G_{x}(i,j)=\left\langle \bm{a}_{i},\bm{a}_{j}\right\rangle_{x} $, then 
	$\left\langle u,v\right\rangle_{x}=\hat{u}^{\top}G_{x}\hat{v} $.
	
	When we regard the unit quaternion $\mathbb{U}$ as a sphere $\mathcal{S}^{3}$ embedded in $\mathbb{R}^{4}$, it is a Riemannian submanifold with the inherited metric $\left\langle u,v\right\rangle_{x}=\left\langle u,v\right\rangle=u^{\top}v$ and tangent space 
	$$
	T_{x}\mathbb{U}=\{v\in\mathbb{R}^{4}:x^{\top} v=0\}.
	$$ 
	The Riemannian gradient of $f \in \mathcal{F}_{x}$ is 
	$$
	\operatorname{grad}f(x)=(I_{4}-xx^{\top})\nabla f(x).
	$$

\subsection{Optimization theory over manifolds}
\label{sec-Optimization Theory over Manifolds}
Let $C$ be a nonempty closed subset in $\mathbb{R}^{n}$ and $x \in C$, the tangent cone is defined by $T_C(x):=\left\{y: \exists t_k \downarrow 0, x_k \rightarrow x, y_k \rightarrow y,~ \text{s.t.}~x_k=x+t_k y_k \in C\right\}$. The normal cone is the dual cone of $-T_C(x)$ or the polar of $T_C(x)$, which is defined as $N_C(x):=\left\{z: z^{\top} y \leq 0, \forall y \in T_C(x)\right\}$. Suppose $S$ is a closed subset on the Riemannian manifold $\mathcal{M},(U, \varphi)$ is a chart at point $x \in S$, then by using coordinate transform , the Riemannian tangent cone can be defined as
$$
T_S(x):=[D \varphi(x)]^{-1}\left[T_{\varphi(S \cap U)}(\varphi(x))\right] .
$$
Consequently, the Riemannian normal cone (see \cite{yang2014optimality}) can be computed by
$$
N_S(x)=[D \varphi(x)]^{-1}\left[G_x^{-1} N_{\varphi(U \cap S)}(\varphi(x))\right].
$$	
\begin{lemma} 
	\cite{zhang2020primal}
	Consider the following optimization problem over a manifold:
	\begin{equation}
		\min ~f(x) \quad \text{s.t.}~x\in S\subset\mathcal{M},
		\label{manifold problem1}
	\end{equation}
	$x^{*}$ is a stationary point if it satisfies
	$$
	0 \in \partial f(x^{*}) +N_{S}(x^{*}).
	$$
	where $\partial f$ is the Riemannian Clarke subdifferential defined in \cite{yang2014optimality}. If f is differentiable, then the first-order optimality condition reduces to
	$$
	-\operatorname{grad} f(x^{*})\in N_{S}(x^{*}),
	$$
	or, equivalently 
	$$
	0=\operatorname{dist}(-\operatorname{grad} f(x^{*}),N_{S}(x^{*})).
	$$
	\label{lem-first-order optimality condition1}
\end{lemma}
\vspace{-2em}
If there exists  equality constraints in problem \eqref{manifold problem1}, i.e. the feasible set can be described as $S=\Omega\cap X$, where $\Omega=\{x\in \mathcal{M}\mid c_{i}(x)=0,i=1,\dots m\}$ and $X$ is an nonempty bounded set. Then we have the following lemma:
\begin{lemma}
	\cite{zhang2020primal}
	Suppose that $x^{*} \in \mathcal{M}$, $c_{i}(x^{*})=0,i=1,\dots m$, and
	\begin{equation}
		0 \in \partial f(x^{*})+ N_{\Omega}(x^{*}) + N_{X}(x^{*}),
	\end{equation}
	where
	$$
	N_{\Omega}\left(x^*\right)=\left\{\sum_{i=1}^m \lambda_i \operatorname{grad} c_i\left(x^*\right) \mid \lambda \in \mathbb{R}^m\right\},
	$$
	then $x^{*}$ is a stationary solution of problem \eqref{manifold problem1}.
	\label{lem-first-order optimality condition2}
\end{lemma}

\section{The Augmented Unit Quaternion Model} \label{sec-AUQ model}
In this section, we propose a new PGO model based on augmented unit quaternions. The traditional PGO models represent the unknown poses and the measurements in $SE(3)$ or $se(3)$. Instead of orthogonal rotation matrix $SO(3)$ or aixs-angle representation $so(3)$, we use unit quaternions in $\mathbb{U}$ to represent rotations, which is a simple spherical manifold and can be solved efficiently by proximal operators.

\begin{figure}[t]
	\centering
	\subfloat[$\kappa=1$]{\label{VMF_k_1}	\includegraphics[width=0.3\linewidth]{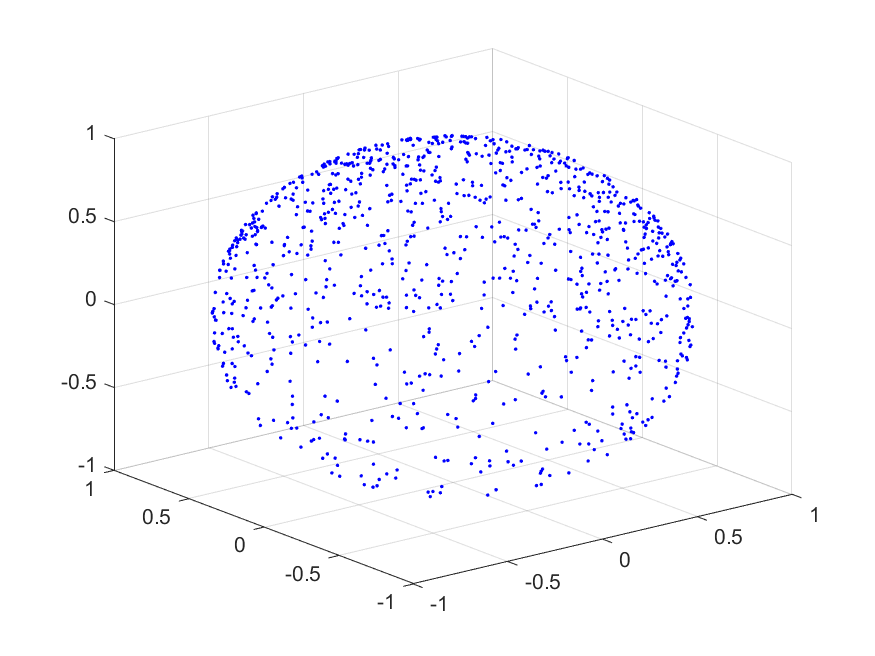}}
	\subfloat[$\kappa=5$]{\label{VMF_k_5}	\includegraphics[width=0.3\linewidth]{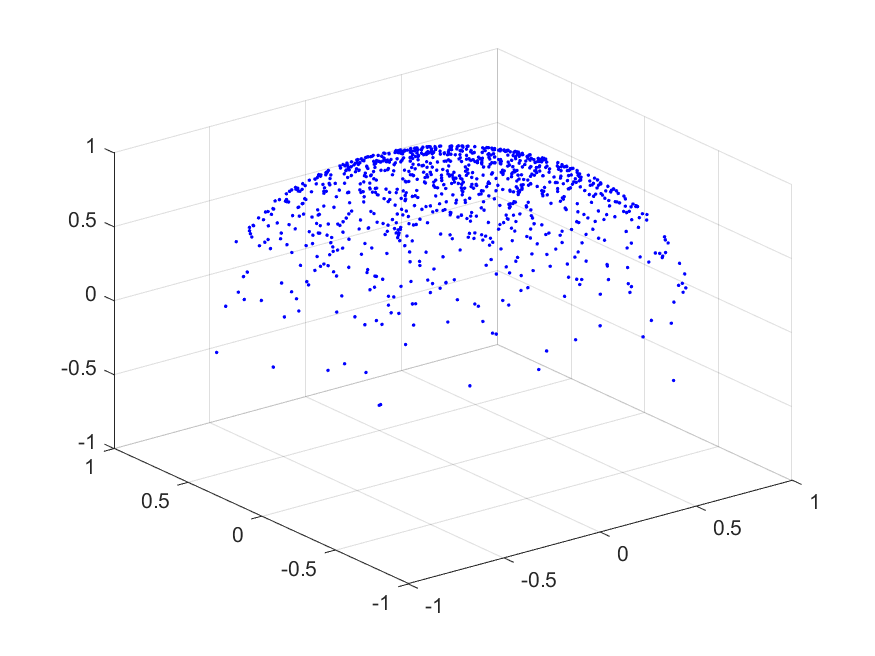}}
	\subfloat[$\kappa=20$]{\label{VMF_k_20}	\includegraphics[width=0.3\linewidth]{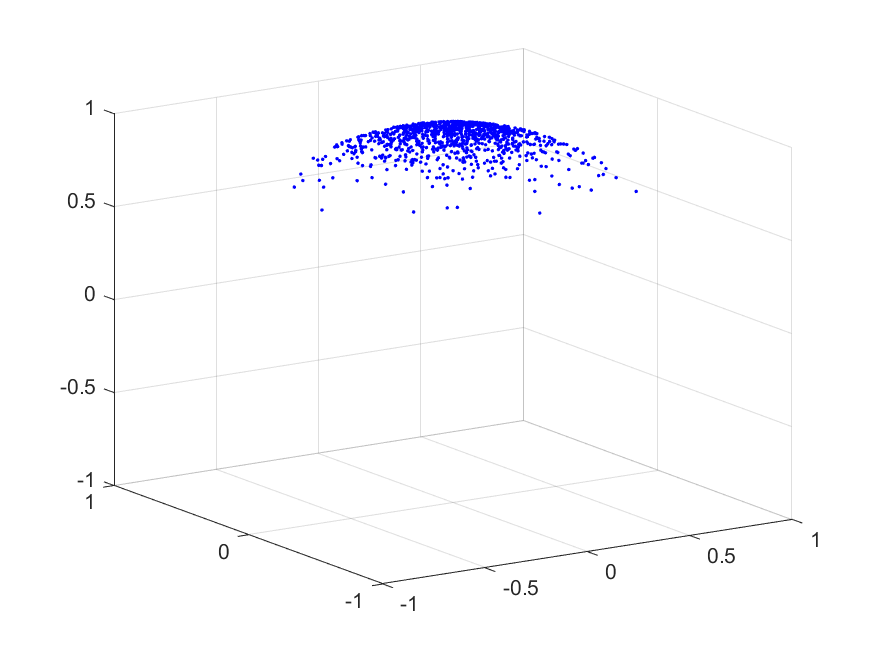}}
	\caption{The illustration of the  von Mises-Fisher distribution with $\bm{\mu}=[0,0,1]^{\top}$ and $\kappa=1,~5,~20$, respectively.}
	\label{fig-VMF}
\end{figure}

PGO can be visualized as a directed graph $\mathcal{G}=(\mathcal{V},\mathcal{E})$, see \cite{carlone2015lagrangian,rosen2015convex,grisetti2009nonlinear},  in which each vertex $i \in \mathcal{V} $ corresponds to a robot pose $(\tilde{q}_{i},\bm{t}_{i}) \in \mathbb{AU}$, and each directed edge $(i,j) \in \mathcal{E}$ corresponds to a relative measurement $(\tilde{q}_{ij},\bm{t}_{ij}) \in \mathbb{AU}$. We define $n=|\mathcal{V}|$ and $m=|\mathcal{E}|$ which indicate the number of vertices and edges, respectively. The goal is to estimate the unknown poses from the noisy measurements. Now, we assume the following generative model for the relative pose measurements: 
\begin{equation*}
	\begin{aligned}
		\bm{t}_{ij}&=R_{i}^{\top}(\bm{t}_{j}-\bm{t}_{i})+\bm{t}_{\epsilon},    \quad \text{where}~ \bm{t}_{\epsilon} \sim  ~ \mathcal{N}(0,\Omega_{1})\\
		\tilde{q}_{ij}&=\tilde{q}_{i}^{*}\tilde{q}_{j}\tilde{q}_{\epsilon}, \quad\quad\quad\quad\quad~ \text{where}~\tilde{q}_{\epsilon}\sim  ~ \text{vMF}([1,0,0,0],\kappa)
	\end{aligned}
\end{equation*}
where $R_{i}$ and $\tilde{q}_{i}$ are the rotation representation of vertex $i$ in $SO(3)$ and unit quaternion, respectively.
``$\mathcal{N}(\bm{\mu},\Omega)$'' denotes a Gaussian distribution with mean $\bm{\mu}$ and covariance matrix $\Omega$. 
``$\text{vMF}(\bm{\mu},\kappa)$'' denotes a $d$-dimensional von Mises-Fisher distribution where $\bm{\mu} \in \mathbb{S}^{d-1}$ and $\kappa \geq 0$ are mean direction and concentration parameters, respectively. It is one of the most commonly used distributions to model data distributed on the surface
of the unit hypersphere \cite{fisher1953dispersion,sra2012short,hornik2014maximum} and can be considered a circular analogue of the normal distribution.  Its probability density function is given by 
\begin{equation*}
	p(\bm{x};\bm{\mu}, \kappa)=c_{d}(\kappa)e^{\kappa \bm{\mu}^\top \bm{x}}
\end{equation*}
where $\bm{x}\in \mathbb{S}^{d-1}$ is a $d$-dimensional unit vector, and the normalizing constant $c_{d}(\kappa)$ has the form 
\begin{equation*}
	c_{d}(\kappa)=\frac{\kappa^{d/2-1}}{(2\pi)^{d/2}I_{d/2-1}(\kappa)}
\end{equation*}
where $I_{r}(\cdot)$ is the modified Bessel function of the first kind of order $r$. In fact, as the concentration parameter $\kappa$ increases, the vMF distribution becomes increasingly concentrated at the mean direction $\bm{\mu}$. When $\kappa=0$, it corresponds to the uniform distribution on $\mathbb{S}^{d-1}$. When $\kappa\rightarrow+\infty$, the distribution approximates to a Gaussian distribution with mean $\bm{\mu}$ and covariance  $1/\kappa$. We show the von Mises-Fisher distribution in Fig.~\ref{fig-VMF}.

The parameter $d$ is the spatial dimension. In the case $d=2$, the von Mises-Fisher  distribution degenerates to the von Mises distribution. And when $d=3$, it is called the Fisher distribution.
We set $d=4$ which is corresponding to the unit quaternion. Then, we can build our model by the maximum likelihood estimator (MLE) or minimizing the negative log-likelihood:
\begin{equation*}
	\underset{\{\tilde{q}_i\} \in \mathbb{U},\{\bm{t}_{i}\} \in \mathbb{R}^3}
	{\min}
	\sum_{(i,j)\in \mathcal{E}} -\log l (\tilde{q}_{ij}) -\log l (\bm{t}_{ij}) ,
\end{equation*}
where $l(\tilde{q}_{ij})$ and $l (\bm{t}_{ij})$ are likelihood functions defined by:
\begin{align*}
	l(\tilde{q}_{ij}) & = c_{4}(\kappa)
	\exp(\kappa \left\langle [1,0,0,0], \tilde{q}_{j}^{*}\tilde{q}_{i}\tilde{q}_{ij}\right\rangle),
	\\
	l (\bm{t}_{ij}) & = \frac{1}{(2\pi)^{3/2}|\Omega_{1}|^{1/2}} 
	\exp\left(-\frac{1}{2}
	\|\bm{t}_{ij}-R_{i}^{\top}(\bm{t}_{j}-\bm{t}_{i})\|_{\Omega_{1}^{-1}}^{2}
	\right),
\end{align*}
respectively. 
The negative log-likelihood for $\bm{t}_{ij}$ with Gaussian distributions is 
\begin{equation*}
	-\log l (\bm{t}_{ij})=\frac{1}{2}\|\bm{t}_{ij}-R_{i}^{\top}(\bm{t}_{j}-\bm{t}_{i})\|_{\Omega_{1}^{-1}}^{2}+const,
\end{equation*} 
where $\|\bm{x}\|_{\Omega}^{2}=\bm{x}^{\top}\Omega\bm{x}$. The negative log-likelihood for $\tilde{q}_{ij}$ with vMF distributions is 
\begin{equation*}
	\begin{aligned}
		-\log l (\tilde{q}_{ij})
		&=-\kappa \left\langle [1,0,0,0], \tilde{q}_{\epsilon}\right\rangle +const\\
		&=-\kappa \left\langle [1,0,0,0], \tilde{q}_{j}^{*}\tilde{q}_{i}\tilde{q}_{ij}\right\rangle +const\\
		&=\frac{\kappa}{2}\|\tilde{q}_{j}^{*}\tilde{q}_{i}\tilde{q}_{ij}-1\|^{2}+const,
	\end{aligned}
\end{equation*} 
where the last equality follows from the fact that $\tilde{q}_{j}^{*}\tilde{q}_{i}\tilde{q}_{ij}$ is still a unit quaternion. 

Based on the assumptions of noise and Lemma \ref{lem-SO3-quaternion}, we can get the following optimization problem:
\begin{equation*}
	\begin{aligned}
		\underset{\{\tilde{q}_i\},\{\bm{t}_{i}\}}{\min}	&~
		\sum_{(i,j)\in \mathcal{E}} 
		\|\tilde{t}_j - \tilde{t}_i-\tilde{q}_{i}\tilde{t}_{ij}\tilde{q}_{i}^{*}\|_
		{\Omega_2}^2
		+\|\tilde{q}_{j}^{*}\tilde{q}_{i}\tilde{q}_{ij}-1\|_{
			\kappa I_{4}}^2, \\
		\operatorname{s.t.} ~~& ~ 
		\tilde{q}_i \in \mathbb{U},~ \bm{t}_{i} \in \mathbb{R}^3,~ i=1,2,\dots,n,
	\end{aligned}
\end{equation*}
where $\Omega_2 = \left( \begin{array}{cc}
	c & 0\\
	0 & \Omega_{1}^{-1}
\end{array}\right)$, and $c>0$ is an arbitrary scalar. By generalizing the covariance matrix $\Omega_{1}$ and concentration parameters $\kappa$,
we propose our PGO model in the following:
\begin{eqnarray}
	\nonumber	&&\underset{\{\tilde{q}_i\},\{\bm{t}_{i}\}}{\min}	~
		\sum_{(i,j)\in \mathcal{E}} 
		\|\tilde{t}_j - \tilde{t}_i-\tilde{q}_{i}\tilde{t}_{ij}\tilde{q}_{i}^{*}\|_{\Sigma_{1}}^2
		+\|\tilde{q}_{j}^{*}\tilde{q}_{i}\tilde{q}_{ij}-1\|_{\Sigma_{2}}^2, \\
		&&\operatorname{s.t.} ~ ~~
		\tilde{q}_i \in \mathbb{U},~ \bm{t}_{i} \in \mathbb{R}^3,~ i=1,2,\dots,n, 
	\label{PGO-model1}
\end{eqnarray}
where $\Sigma_{1}$, $\Sigma_{2} \in \mathcal{S}_{+}^{n}$ are positive semi-definite matrices. $\tilde{t}_{i}=[0,\bm{t}_{i}]$ is a pure quaternion. The model \eqref{PGO-model1} is an unconstrained quartic polynomial optimization problem on manifold. With the development of manifold optimization, many algorithms with convergence guarantees can deal with the model \eqref{PGO-model1}, such as proximal Riemannian gradient method \cite{gabay1982minimizing}, Riemannian conjugate gradient method \cite{smith1994optimization}, and Riemannian Newton method \cite{hu2018adaptive}. We refer the readers to \cite{hu2020brief} for more details.

By exploiting the structure of model \eqref{PGO-model1}, we focus on the splitting method, whose subproblems are usually easier to solve. 
With respect to the variable $\bm{t}_{i}$, model \eqref{PGO-model1} constitutes a quadratic optimization problem for which an explicit solution can be directly computed. With respect to variable $\tilde{q}_{i}$, the optimization problem is quartic, which may make the model difficult to solve. 
By introducing auxiliary variables $\tilde p_i$, $i=1,\dots,n$, we can get an equivalent model of \eqref{PGO-model1}:
\begin{equation}
	\begin{aligned}
		\underset{\{\tilde{p}_i\},\{\tilde{q}_i\}, \{\bm{t}_{i}\}}{\min}	&~
		\sum_{(i,j)\in \mathcal{E}} 
		\|\tilde{t}_j - \tilde{t}_i-\tilde{q}_{i}\tilde{t}_{ij}\tilde{p}_{i}^{*}\|_{\Sigma_{1}}^2
		+\|\tilde{p}_{j}^{*}\tilde{q}_{i}\tilde{q}_{ij}-1\|_{\Sigma_{2}}^2, \\
		\operatorname{s.t.} ~\quad&  
		\tilde{p}_i=\tilde{q}_i,~
		\tilde{p}_i \in \mathbb{U},~
		\tilde{q}_i \in \mathbb{R}^{4},~ 
		\bm{t}_{i} \in \mathbb{R}^3,~ 
		i=1,\dots,n.
	\end{aligned}
	\label{PGO-model2}
\end{equation}
When considering the variables separately, the model \eqref{PGO-model2} is reformulated to a multi-linear least square problem, whose subproblems are easier to solve than that of the  nonlinear least square  problem \eqref{PGO-model1}. The variables $\tilde{p}_{i}$ and $\tilde{q}_{i}$ in the objective function are coupled, while each individual variable retains linearity and convexity.  Using the alternating update strategy, it can be seen in the next section that the subproblems  have closed-form solutions, and can also be solved in parallel corresponding to the structure of the directed graph $\mathcal{E}$, which will greatly improve the algorithm efficiency.

\section{Proximal Linearized Riemannian ADMM} 
\label{sec-proximal linearized ADMM for PGO Model}
In this section, we propose a proximal linearized Riemannian ADMM algorithm to solve PGO model \eqref{PGO-model2}. 

Let $\tilde{\bm{p}}=(\tilde{p}_{1},\tilde{p}_{2},\dots,\tilde{p}_{n}) \in \mathbb{U}^{n}$, $\tilde{\bm{q}}=(\tilde{q}_{1},\tilde{q}_{2},\dots,\tilde{q}_{n}) \in \mathbb{R}^{4n}$, and $\bm{t}=(\bm{t}_{1},\bm{t}_{2},\dots,\bm{t}_{n}) \in \mathbb{R}^{3n}$, where $\mathbb{U}^{n}$ represents the Cartesian product of $n$ unit quaternion sets. We define
\begin{align}
	f(\tilde{\bm{p}},\tilde{\bm{q}},\bm{t})
	&=\sum_{(i,j)\in \mathcal{E}} 
	\|\tilde{t}_j - \tilde{t}_i-\tilde{q}_{i}\tilde{t}_{ij}\tilde{p}_{i}^{*}\|_{\Sigma_{1}}^2,
	\label{function_f1}\\
	g(\tilde{\bm{p}},\tilde{\bm{q}})
	&=\sum_{(i,j)\in \mathcal{E}} 
	\|\tilde{p}_{j}^{*}\tilde{q}_{i}\tilde{q}_{ij}-1\|_{\Sigma_{2}}^2.
	\label{function_g1}
\end{align}
Then the augmented Lagrangian function of PGO model \eqref{PGO-model2} is 
\begin{equation}
		\mathcal{L}_{\beta}(\tilde{\bm{p}},\tilde{\bm{q}},\bm{t},\bm{\lambda})=
		f(\tilde{\bm{p}},\tilde{\bm{q}},\bm{t})+g(\tilde{\bm{p}},\tilde{\bm{q}})
		+\sum_{i=1}^{n} \left\lbrace 
		\delta_{\mathbb{U}}(\tilde{p}_{i})
		-\left\langle \bm{\lambda}_{i}, \tilde{p}_{i}-\tilde{q}_{i} \right\rangle
		+\frac{\beta}{2}\|\tilde{p}_{i}-\tilde{q}_{i}\| ^2\right\rbrace ,
\end{equation}
where $\bm{\lambda} \in \mathbb{R}^{4n}$ is the Lagrange multiplier and $\beta>0$ is a penalty parameter. The function $\delta_{\mathbb{U}}(\cdot) $ is the indicator function of $\mathbb{U}$ which is defined as
\begin{equation*}
	\delta_{\mathbb{U}}(\tilde{p}):=
	\left\lbrace 
	\begin{aligned}
		&0,       \qquad \text{if}~ \tilde{p} \in \mathbb{U},\\
		&+\infty, ~ \text{otherwise}.
	\end{aligned}
	\right. 
\end{equation*}
The iterative scheme of classical ADMM is given by 
\begin{subequations}
	\begin{numcases}{}
		\tilde{\bm{p}}^{k+1}=
		\mathop{\arg\min}\limits_{\tilde{\bm{p}}\in\mathbb U^n}
		\mathcal{L}_{\beta}(\tilde{\bm{p}},\tilde{\bm{q}}^{k},\bm{t}^{k},\bm{\lambda}^{k})
		,
		\label{subproblem_p_1_c}\\
		\tilde{\bm{q}}^{k+1}=
		\mathop{\arg\min}\limits_{\tilde{\bm{q}}}
		\mathcal{L}_{\beta}(\tilde{\bm{p}}^{k+1},\tilde{\bm{q}},\bm{t}^{k},\bm{\lambda}^{k})
		,
		\label{subproblem_q_1_c}\\
		\bm{t}^{k+1}=
		\mathop{\arg\min}\limits_{\bm{t}}
		\mathcal{L}_{\beta}(\tilde{\bm{p}}^{k+1},\tilde{\bm{q}}^{k+1},\bm{t},\bm{\lambda}^{k})
		,
		\label{subproblem_t_1_c}\\
		\bm{\lambda}^{k+1}=\bm{\lambda}^{k}-\beta(\tilde{\bm{p}}^{k+1}-\tilde{\bm{q}}^{k+1}).
		\label{subproblem_lambda_1_c}
	\end{numcases}
\end{subequations}
 In the above algorithm, the optimization subproblem \eqref{subproblem_q_1_c} and \eqref{subproblem_t_1_c} can be solved  efficiently.
However, because of   the manifold constraints, there is no closed-form solution for \eqref{subproblem_p_1_c}. 

 The linearized technique can help us overcome this difficulty and get a closed-form solution.
In \cite{ouyang2015accelerated,gao2018information}, the authors linearized the whole augmented Lagrangian function and obtained easier subproblems. However, this may result in the slow convergence. Instead, we only use the linearization of $f$ and $g$, and keep the quadratic term $\|\tilde{p}_{i}-\tilde{q}_{i}\|^{2}$. Consequently, the linearized augmented Lagrangian function about $\tilde{\bm{p}}$ is defined as 
\begin{equation}
		\begin{aligned}
			\mathcal{L}_{\beta}^{k}(\tilde{\bm{p}},\tilde{\bm{q}}^{k},\bm{t}^{k},\bm{\lambda}^{k})
			&=f(\tilde{\bm{p}}^{k},\tilde{\bm{q}}^{k},\bm{t}^{k})
			+\left\langle \nabla_{\tilde{\bm{p}}} f(\tilde{\bm{p}}^{k},\tilde{\bm{q}}^{k},\bm{t}^{k}),\tilde{\bm{p}}-\tilde{\bm{p}}^{k}
			\right\rangle \\
			&\quad +g(\tilde{\bm{p}}^{k},\tilde{\bm{q}}^{k})
			+\left\langle \nabla_{\tilde{\bm{p}}}
			g(\tilde{\bm{p}}^{k},\tilde{\bm{q}}^{k}),\tilde{\bm{p}}-\tilde{\bm{p}}^{k}
			\right\rangle \\
			&\quad + \sum_{i=1}^{n} \left\lbrace 
			\delta_{\mathbb{U}}(\tilde{p}_{i})
			-\left\langle \bm{\lambda}_{i}^{k}, \tilde{p}_{i}-\tilde{q}_{i}^{k} \right\rangle
			+\frac{\beta}{2}\|\tilde{p}_{i}-\tilde{q}_{i}^{k}\| ^2\right\rbrace ,
		\end{aligned}
	\label{eq-Lk}	
 \end{equation}
		
 Furthermore, an extra proximal term \cite{he2002new} not only can guarantee the uniqueness of the solution, but also provide a quantifiable descending of augmented Lagrangian function, which will help us analyze the convergence. 

Finally, the iterative scheme of proximal linearized ADMM is given by
\begin{subequations}
	\begin{numcases}{}
		\tilde{\bm{p}}^{k+1}=
		\mathop{\arg\min}\limits_{\tilde{\bm{p}}\in\mathbb U^n}
		\mathcal{L}_{\beta}^{k}(\tilde{\bm{p}},\tilde{\bm{q}}^{k},\bm{t}^{k},\bm{\lambda}^{k})
		\!+\!\frac{1}{2}\|\tilde{\bm{p}}-\tilde{\bm{p}}^{k}\|_{H_1}^2,
		\label{subproblem_p_1}\\
		\tilde{\bm{q}}^{k+1}=
		\mathop{\arg\min}\limits_{\tilde{\bm{q}}}
		\mathcal{L}_{\beta}(\tilde{\bm{p}}^{k+1},\tilde{\bm{q}},\bm{t}^{k},\bm{\lambda}^{k})
		\!+\!\frac{1}{2}\|\tilde{\bm{q}}-\tilde{\bm{q}}^{k}\|_{H_2}^2,
		\label{subproblem_q_1}\\
		\bm{t}^{k+1}=
		\mathop{\arg\min}\limits_{\bm{t}}
		\mathcal{L}_{\beta}(\tilde{\bm{p}}^{k+1},\tilde{\bm{q}}^{k+1},\bm{t},\bm{\lambda}^{k})
		\!+\!\frac{1}{2}\|\bm{t}-\bm{t}^{k}\|_{H_3}^2,
		\label{subproblem_t_1}\\
		\bm{\lambda}^{k+1}=\bm{\lambda}^{k}-\beta(\tilde{\bm{p}}^{k+1}-\tilde{\bm{q}}^{k+1}),
		\label{subproblem_lambda_1}
	\end{numcases}
\end{subequations}
where 
$H_1$, $H_2$, $H_3 \succ 0$ are positive definite matrices with block diagonal structures. Here, $\frac{1}{2}\|\tilde{\bm{p}}-\tilde{\bm{p}}^{k}\|_{H_1}^2$ can be split as 
$\frac{1}{2}\sum_{i=1}^{n}\|\tilde{p}_{i}-\tilde{p}_{i}^{k}\|_{H_{1,i}}^{2}$. 
In other words, each block of $H_{1}$ is still a diagonal matrix with the form  $H_{1,i} = \tau_{1,i}I_{4}$. 
We assume $H_2$ has similar structures.

\subsection{Subproblems}
In this section, we give the closed-form solutions to the subproblems \eqref{subproblem_p_1}-\eqref{subproblem_t_1} and an algorithm framework in parallel.

First of all, we partition the given directed graph $\mathcal{G}=(\mathcal{V},\mathcal{E})$ according to the vertices.
We define $\mathcal{E}_{i}^{in}=\{(l,i)\in\mathcal{E} \}$ for all $l \in \mathcal{V}$, and $\mathcal{E}_{i}^{out}=\{(i,j)\in\mathcal{E} \}$ for all $j \in \mathcal{V}$. In other words, $\mathcal{E}_{i}^{in}$ represents all directed edges that pointing to vertex $i$, while $\mathcal{E}_{i}^{out}$ is the opposite. Then we have the properties that 
$$
\mathcal{E}=\bigcup_{i\in\mathcal{V}}\left( \mathcal{E}_{i}^{in}\cup \mathcal{E}_{i}^{out} \right),\quad
\mathcal{E}_{i}^{in} \cap \mathcal{E}_{i}^{out} = \varnothing \text{ for all } i\in\mathcal{V},\quad
\text{ and } \mathcal{E}_{i}^{in} \cap \mathcal{E}_{j}^{in} = \varnothing ,\text{ for all } i\neq j.
$$

For the $\tilde{\bm{p}}$-subproblem, according to Lemma \ref{lem-MW}, we can transform the multiplication between two quaternions into the multiplication between a matrix and a vector, and rewrite the functions $f,g$ in \eqref{function_f1} and \eqref{function_g1} as
\begin{align}
	f(\tilde{\bm{p}},\tilde{\bm{q}},\bm{t})
	&=\sum_{(i,j)\in \mathcal{E}} 
	\|M(\tilde{q}_{i})M(\tilde{t}_{ij})D\tilde{p}_{i}-(\tilde{t}_j - \tilde{t}_i) \|_{\Sigma_{1}}^2,
	\label{function_f2}\\
	g(\tilde{\bm{p}},\tilde{\bm{q}})
	&=\sum_{(i,j)\in \mathcal{E}} 
	\|W(\tilde{q}_{ij})W(\tilde{q}_{i})\tilde{p}_{j}-1\|_{\Sigma_{2}}^2.
	\label{function_g2}
\end{align}
where the matrix $D=\text{diag}(1,-1,-1,-1)$ is a diagonal matrix of size $4 \times 4$. The gradient of functions $f,g$ can be calculated as
\begin{align}
		\nabla_{\tilde{p}_{i}}f(\tilde{\bm{p}},\tilde{\bm{q}},\bm{t})
		&\!=\!\!\!\!\sum_{(i,j)\in \mathcal{E}_{i}^{out}} 
		\!\!\!\!2D M(\tilde{t}_{ij})^{\!\top\!} M(\!\tilde{q}_{i}\!)^{\!\top\!} \Sigma_{1} 
		\left( M(\tilde{q}_{i})M(\tilde{t}_{ij})D\tilde{p}_{i}\!-\!(\tilde{t}_j \!-\! \tilde{t}_i) \right), 
		\nonumber\\
		\nabla_{\tilde{p}_{i}}g(\tilde{\bm{p}},\tilde{\bm{q}})
		&\!=\!\sum_{(l,i)\in \mathcal{E}_{i}^{in}} 
		2W(\tilde{q}_{i})^{\top}W(\tilde{q}_{ij})^{\top} \Sigma_{2}
		\left( W(\tilde{q}_{ij})W(\tilde{q}_{i})\tilde{p}_{j}-1 \right).
		\nonumber
\end{align}

The subproblem \eqref{subproblem_p_1} of $\tilde{\bm{p}}$ can be written as 	
\begin{align}
	\tilde{\bm{p}}^{k+1} &\!=\! 
	\mathop{\arg\min}\limits_{\tilde{\bm{p}} \in \mathbb{U}^{n}}	
	\sum_{i=1}^{n} \left\lbrace 
	\left\langle \nabla_{\tilde{p}_{i}}f(\tilde{\bm{p}}^{k},\tilde{\bm{q}}^{k},\bm{t}^{k})
	\!+\!\nabla_{\tilde{p}_{i}}g(\tilde{\bm{p}}^{k},\tilde{\bm{q}}^{k}), \tilde{p}_{i}\right\rangle 
	\right. 	
	\nonumber\\
	& \quad +\left. 
	\frac{\beta}{2}\|\tilde{p}_{i}-(\tilde{q}_{i}^{k}+\frac{1}{\beta}\bm{\lambda}_{i}^{k})\|^{2}
	+\frac{1}{2}\|\tilde{p}_{i}-\tilde{p}_{i}^{k}\|_{H_{1,i}}^2 \right\rbrace \nonumber.
\end{align}
Since $\tilde{p}_{i}$ are fully separable, we can update them in parallel. For $i=1,2,\dots,n$, we have
\begin{align}
		\tilde{p}_{i}^{k+1} &= \operatorname{Proj}_{\mathbb{U}}
		\left( (\beta I_4{}+H_{1,i})^{-1}\left( 
		\beta \tilde{q}_{i}^{k} + \bm{\lambda}_{i}^{k} + H_{1,i}\tilde{p}_{i}^{k}
		-\nabla_{\tilde{p}_{i}}f(\tilde{\bm{p}}^{k},\tilde{\bm{q}}^{k},\bm{t}^{k})
		-\nabla_{\tilde{p}_{i}}g(\tilde{\bm{p}}^{k},\tilde{\bm{q}}^{k})
		\right) \right) ,
		\label{subproblem_p_2}
\end{align}
where the operator $\operatorname{Proj}_{\mathbb{U}}(\bm{x})=\frac{\bm{x}}{\|\bm{x}\|}$ is the projection on $\mathbb{U}$ when $x$ is non-zero.

The subproblem \eqref{subproblem_q_1} can be written as
\begin{align}
		\tilde{\bm{q}}^{k+1} &= 
		\mathop{\arg\min}\limits_{\tilde{\bm{q}} \in \mathbb{R}^{4n}}
		\sum_{(i,j)\in \mathcal{E}} 
		\|W(\tilde{p}_{i}^{k+1})^{\top}W(\tilde{t}_{ij})\tilde{q}_{i}-(\tilde{t}_j^{k} - \tilde{t}_i^{k}) \|_{\Sigma_{1}}^2
		+\|W(\tilde{q}_{ij})M(\tilde{p}_{j}^{k+1})^{\top}\tilde{q}_{i}-1\|_{\Sigma_{2}}^2.
		\nonumber\\	
		&\quad +\sum_{i=1}^{n} \left\lbrace 
		\frac{\beta}{2}\|\tilde{q}_{i}-(\tilde{p}_{i}^{k+1}-\frac{1}{\beta}\bm{\lambda}_{i}^{k})\|^{2}
		+\frac{1}{2}\|\tilde{q}_{i}-\tilde{q}_{i}^{k}\|_{H_{2,i}}^2 \right\rbrace .
		\nonumber
	\end{align}
	Similarly, we can update $\tilde{q}_{i}$ in parallel. 
	For $i=1,2,\dots,n$, we have
	\begin{align}
		\tilde{q}_{i}^{k+1} &= 
		\left( \sum_{(i,j)\in \mathcal{E}_{i}^{out}}
		2\left\lbrace
		G_{1,ij}^{\top}\Sigma_{1}G_{1,ij}
		+G_{2,ij}^{\top}\Sigma_{2}G_{2,ij}
		\right\rbrace +\beta I_{4}+H_{2,i}
		\right) ^{-1}
		\nonumber\\
		&\left( \sum_{(i,j)\in \mathcal{E}_{i}^{out}} 
		2\left\lbrace 
		G_{1,ij}^{\top}\Sigma_{1}(\tilde{t}_j^{k} - \tilde{t}_i^{k})
		+G_{2,ij}^{\top}\Sigma_{2}1		
		\right\rbrace +\beta\tilde{p}_{i}^{k+1}-\bm{\lambda}_{i}^{k}+H_{2,i}\tilde{q}_{i}^{k}
		\right),
		\label{subproblem_q_2}
	\end{align}
	where $G_{1,ij}=W(\tilde{p}_{i}^{k+1})^{\top}W(\tilde{t}_{ij})$, and $G_{2,ij}=W(\tilde{q}_{ij})M(\tilde{p}_{j}^{k+1})^{\top}$. 

Then,  denote 
$$
\Sigma_{1}=\left( \begin{array}{*{20}{c}}
	\sigma_{11}& \bm{\sigma}_{12}^{\top}\\
	\bm{\sigma}_{21}& \hat{\Sigma}_{1}
\end{array}
\right). 
$$
For the $\bm{t}$-subproblem, there is
\begin{align}
	\bm{t}^{k+1}&=
	\mathop{\arg\min}\limits_{\bm{t} \in \mathbb{R}^{3n}}
	f(\tilde{\bm{p}}^{k+1},\tilde{\bm{q}}^{k+1},\bm{t})
	+\frac{1}{2}\|\bm{t}-\bm{t}^{k}\|_{H_3}^2
	\nonumber\\
	&=\mathop{\arg\min}\limits_{\bm{t} \in \mathbb{R}^{3n}}
	\!\!\!\sum_{(i,j)\in\mathcal{E}}\!\!\!\|\tilde{t}_j \!-\! \tilde{t}_i \!-\!\tilde{q}_{i}^{k+1}\tilde{t}_{ij}(\tilde{p}_{i}^{k+1})^{*}\!\|_{\Sigma_{1}}^{2}
	\!+\!\frac{1}{2}\!\|\bm{t}\!-\!\bm{t}^{k}\|_{H_3}^2
	\nonumber\\
	&=\mathop{\arg\min}\limits_{\bm{t} \in \mathbb{R}^{3n}}
	\sum_{(i,j)\in\mathcal{E}}
	\|Q_{ij} \bm{t}-\bm{s}_{ij}^{k+1}
	\|_{\hat{\Sigma}_{1}}^{2}
	+\frac{1}{2}\|\bm{t}-\bm{t}^{k}\|_{H_3}^2
	\nonumber\\
	&=\mathop{\arg\min}\limits_{\bm{t} \in \mathbb{R}^{3n}}
	\|Q\bm{t}-\bm{s}^{k+1}\|_{I_{m}\otimes {\hat{\Sigma}_{1}}}^{2}
	+\frac{1}{2}\|\bm{t}-\bm{t}^{k}\|_{H_3}^2,
	\label{subproblem_t_2-1}
\end{align}
where $\tilde{s}_{ij}^{k+1}=\tilde{q}_{i}^{k+1}\tilde{t}_{ij}(\tilde{p}_{i}^{k+1})^{*}$ and $\bm{s}_{ij}^{k+1}$ is the imaginary part of $\tilde{s}_{ij}^{k+1}$. $Q_{ij} \in \mathbb{R}^{3 \times 3n}$ is a block matrix consisting of three-order square matrices where the $i$-th block is $-I_{3}$, the $j$-th block is $I_{3}$ and the others are all zero.
The last equality is a compact form corresponding to the edges of $\mathcal{E}$ in which $Q \in \mathbb{R}^{3m \times 3n}$ is a matrix and $\bm{s}^{k+1} \in \mathbb{R}^{3m}$ is a vector. $\otimes$ denotes the Kronecker product. The minimizer of \eqref{subproblem_t_1} is given explicitly by
\begin{equation}
		\bm{t}^{k+1} \!= \!
		\left( 2Q^{\top}(I_{m}\otimes {\hat{\Sigma}_{1}})Q+H_{3} \right)^{-1}\!
		\left( 2Q^{\top}(I_{m}\otimes {\hat{\Sigma}_{1}})\bm{s}^{k+1}\! +\! H_{3}\bm{t}^{k} \right) .
		\label{subproblem_t_2}
\end{equation}

Now, we are ready to formally present our algorithm. The proximal linearized Riemannian ADMM for solving the PGO model \eqref{PGO-model2} can be described in Algorithm \ref{alg-PieADMM1}.

\begin{algorithm}[t]
	\setstretch{1.2}
	\caption{\textbf{P}roximal L\textbf{i}nearized Ri\textbf{e}mannian \textbf{ADMM} (PieADMM) for solving PGO Model \eqref{PGO-model2}.}
	\label{alg-PieADMM1}
	\hspace*{0.02in} 
	{ \textbf{Input:} $\mathcal{G}=(\mathcal{V},\mathcal{E})$, 
		$\tilde{q}_{ij} \in \mathbb{U}$, $\bm{t}_{ij} \in \mathbb{R}^{3}$ for all $(i,j) \in \mathcal{E}$, $\beta>0$, $H_{1},~H_{2},~H_{3}\succ0$ and initial points $\tilde{\bm{q}}^{1} \in \mathbb{U}^{n}$, $\bm{t}^{1} \in \mathbb{R}^{3n}$, ${\bm{\lambda}}^{1} \in \mathbb{R}^{4n}$}. 
	Set $k=1$.
	\begin{algorithmic}[1]
		\WHILE{stopping criteria is not satisfied}
		\STATE Compute $\tilde{p}_{i}^{k+1}$ for $i=1,\dots,n$ by \eqref{subproblem_p_2};
		\STATE Compute $\tilde{q}_{i}^{k+1}$ for $i=1,\dots,n$ by \eqref{subproblem_q_2};
		\STATE Compute $\bm{t}^{k+1}$ by \eqref{subproblem_t_2} ;
		\STATE $\bm{\lambda}^{k+1}=\bm{\lambda}^{k}-\beta(\tilde{\bm{p}}^{k+1}-\tilde{\bm{q}}^{k+1})$;
		\STATE $k=k+1$;
		\ENDWHILE
	\end{algorithmic}
	\hspace*{0.02in} {\bf Output: $\tilde{\bm{p}}^{k+1} \in \mathbb{U}^{n}$ and $\bm{t}^{k+1} \in \mathbb{R}^{3n}$.}
\end{algorithm}

\section{Convergence Analysis} 
\label{sec-Convergence Analysis}
The PGO model is a nonconvex nonseparable optimization problem with linear equality and manifold constraints, which is difficult to find a global minimum from arbitrary initial points.  By introducing the first-order optimality condition (Section \ref{sec-First-Order Optimality Conditions}), we can prove that our PieADMM algorithm converges to an $\epsilon$-stationary solution with limited number of iterations (Section \ref{sec-Iteration Complexity}), which can also guide us in choosing parameters.

\subsection{First-Order Optimality Conditions} 
\label{sec-First-Order Optimality Conditions}
In this section, we give the first-order optimality condition and $\epsilon$-stationary solution of PGO model \eqref{PGO-model2}.
\begin{theorem}(Optimality Conditions)
	If there exists a Lagrange multiplier $\bm{\lambda}^{*} $ such that 
		\begin{equation}
			\left\lbrace 
			\begin{aligned}
				&0\in\operatorname{Proj}_{T_{\tilde{\bm{p}}^{*}}\mathbb{U}^{n}}
				\left(\nabla_{\tilde{\bm{p}}}f(\tilde{\bm{p}}^{*},\tilde{\bm{q}}^{*},\bm{t}^{*})
				\!+\!\nabla_{\tilde{\bm{p}}}g(\tilde{\bm{p}}^{*},\tilde{\bm{q}}^{*})
				\!-\!\bm{\lambda}^{*}\right) \!+\! N_{\mathbb{U}^{n}}(\tilde{\bm{p}}^{*}),\\
				&0=\nabla_{\tilde{\bm{q}}}f(\tilde{\bm{p}}^{*},\tilde{\bm{q}}^{*},\bm{t}^{*})
				+\nabla_{\tilde{\bm{q}}}g(\tilde{\bm{p}}^{*},\tilde{\bm{q}}^{*})
				+\bm{\lambda}^{*},\\
				&0=\nabla_{\bm{t}}f(\tilde{\bm{p}}^{*},\tilde{\bm{q}}^{*},\bm{t}^{*}),\\
				&0=\tilde{\bm{p}}^{*}-\tilde{\bm{q}}^{*},
			\end{aligned}
			\right. 
			\label{eq-optimality of PGO}
		\end{equation}
	then $(\tilde{\bm{p}}^{*},\tilde{\bm{q}}^{*},\bm{t}^{*})$ is a stationary point of the PGO model \eqref{PGO-model2}.
	\label{thm-optimality of PGO}
\end{theorem}
\begin{proof}
	Let $\Omega=\left\lbrace (\tilde{\bm{p}},\tilde{\bm{q}},\bm{t})\in \mathbb{U}^{n} \times \mathbb{R}^{4n} \times \mathbb{R}^{3n}
	\mid \tilde{\bm{p}}=\tilde{\bm{q}}\right\rbrace$, and $X=\mathbb{U}^{n} \times \mathbb{R}^{4n} \times \mathbb{R}^{3n}$ in Lemma \eqref{lem-first-order optimality condition2}, we can get \eqref{eq-optimality of PGO}.
\end{proof}

Hence, an $\epsilon$-stationary solution of PGO model \eqref{PGO-model2} can be naturally defined as follows.
\begin{definition}($\epsilon$-stationary solution)
	Solution $(\tilde{\bm{p}}^{*},\tilde{\bm{q}}^{*},\bm{t}^{*})$ is said to be an $\epsilon$-stationary solution of PGO model \eqref{PGO-model2} if there exists a Lagrange multiplier $\bm{\lambda}^{*} $ such that 
	\begin{equation}
		\left\lbrace 
		\begin{aligned}
			&\operatorname{dist}\left( 
				\operatorname{Proj}_{T_{\tilde{\bm{p}}^{*}}\mathbb{U}^{n}}
				\left(-\nabla_{\tilde{\bm{p}}}f(\tilde{\bm{p}}^{*},\tilde{\bm{q}}^{*},\bm{t}^{*})
				-\nabla_{\tilde{\bm{p}}}g(\tilde{\bm{p}}^{*},\tilde{\bm{q}}^{*})
				+\bm{\lambda}^{*}\right), N_{\mathbb{U}^{n}}(\tilde{\bm{p}}^{*})\right) \leq \epsilon,
			\\
			&\|\nabla_{\tilde{\bm{q}}}f(\tilde{\bm{p}}^{*},\tilde{\bm{q}}^{*},\bm{t}^{*})
			+\nabla_{\tilde{\bm{q}}}g(\tilde{\bm{p}}^{*},\tilde{\bm{q}}^{*})
			+\bm{\lambda}^{*}\|\leq \epsilon,\\
			&\|\nabla_{\bm{t}}f(\tilde{\bm{p}}^{*},\tilde{\bm{q}}^{*},\bm{t}^{*})
			\|\leq \epsilon,\\
			&\|\tilde{\bm{p}}^{*}-\tilde{\bm{q}}^{*}\|\leq \epsilon.
		\end{aligned}
		\right. 
		\label{eq-e stationary solution of PGO}
	\end{equation}
	\label{def-e stationary solution of PGO}
\end{definition}

\subsection{Iteration Complexity}
\label{sec-Iteration Complexity}
In this section, we establish the global convergence of our proximal linearized Riemannian ADMM algorithm. In addition, we also show the iteration complexity of $O(1/\epsilon^{2})$ to  reach an $\epsilon$-stationary solution.

First of all, we summarize the properties of our PGO model in Proposition \ref{thm_PGO_property}.
\begin{proposition}
		The functions $f(\tilde{\bm{p}},\tilde{\bm{q}},\bm{t})$ and $g(\tilde{\bm{p}},\tilde{\bm{q}})$ in PGO model \eqref{PGO-model2} satisfy the following properties:
		\begin{enumerate}
			\item[(a).] $f(\tilde{\bm{p}},\tilde{\bm{q}},\bm{t})$ and $g(\tilde{\bm{p}},\tilde{\bm{q}})$ are all bounded from below in the feasible region. We denote the lower bounds by 
			$$
			f^{*}=\inf_{\tilde{\bm{p}} \in \mathbb{U}^{n}, \tilde{\bm{q}}\in \mathbb{R}^{4n}, \bm{t}\in \mathbb{R}^{3n}} f(\tilde{\bm{p}},\tilde{\bm{q}},\bm{t}),
			$$
			and
			$$
			g^{*}=\inf_{\tilde{\bm{p}} \in \mathbb{U}^{n}, \tilde{\bm{q}}\in \mathbb{R}^{4n}} g(\tilde{\bm{p}},\tilde{\bm{q}}).
			$$
			\item[(b).] For any fixed $\tilde{\bm{q}}, \bm{t}$, the partial gradient 
			$\nabla_{\tilde{\bm{p}}}f(\tilde{\bm{p}},\tilde{\bm{q}},\bm{t})$ is globally Lipschitz with constant $L_{f,\tilde{\bm{p}}}(\tilde{\bm{q}})>0$, that is 
			\begin{equation}
				\|\nabla_{\tilde{\bm{p}}}f(\tilde{\bm{p}}_{1},\tilde{\bm{q}},\bm{t})
				-\nabla_{\tilde{\bm{p}}}f(\tilde{\bm{p}}_{2},\tilde{\bm{q}},\bm{t})\|
				\leq L_{f,\tilde{\bm{p}}}(\tilde{\bm{q}})
				\|\tilde{\bm{p}}_{1}-\tilde{\bm{p}}_{2}\|, 
				~\forall \tilde{\bm{p}_{1}},\tilde{\bm{p}_{2}} \in \mathbb{R}^{4n},
				\label{eq_PGO_property2.1}
			\end{equation}
			or equivalently,
			\begin{equation}
				f(\tilde{\bm{p}_{1}},\tilde{\bm{q}},\bm{t})
				\leq f(\tilde{\bm{p}_{2}},\tilde{\bm{q}},\bm{t})
				+\left\langle \nabla_{\tilde{\bm{p}}}f(\tilde{\bm{p}_{2}},\tilde{\bm{q}},\bm{t})
				, \tilde{\bm{p}_{1}}-\tilde{\bm{p}_{2}}\right\rangle 
				+\frac{L_{f,\tilde{\bm{p}}}(\tilde{\bm{q}})}{2}
				\|\tilde{\bm{p}}_{1}-\tilde{\bm{p}}_{2}\|^{2},
				~\forall \tilde{\bm{p}_{1}},\tilde{\bm{p}_{2}} \in \mathbb{R}^{4n}.
				\label{eq_PGO_property2.2}
			\end{equation}
			In addition, $\nabla_{\tilde{\bm{q}}}f(\tilde{\bm{p}},\tilde{\bm{q}},\bm{t}), 
			\nabla_{\bm{t}}f(\tilde{\bm{p}},\tilde{\bm{q}},\bm{t}), 
			\nabla_{\tilde{\bm{p}}}g(\tilde{\bm{p}},\tilde{\bm{q}})$, 
			and $\nabla_{\tilde{\bm{q}}}g(\tilde{\bm{p}},\tilde{\bm{q}})$ are also globally Lipschitz with constant 
			$L_{f,\tilde{\bm{q}}}(\tilde{\bm{p}}),
			L_{f,\bm{t}},
			L_{g,\tilde{\bm{p}}}(\tilde{\bm{q}})$, and
			$L_{g,\tilde{\bm{q}}}(\tilde{\bm{p}})$, respectively.
			
			\item[(c).] If $\tilde{\bm{q}}$ lies in a bounded subset, then the Lipschitz constants of the partial gradient 
			$\nabla_{\tilde{\bm{p}}}f(\tilde{\bm{p}},\tilde{\bm{q}},\bm{t})$ and $\nabla_{\tilde{\bm{p}}}g(\tilde{\bm{p}},\tilde{\bm{q}})$ have uniform upper bounds
			, respectively, i.e.,
			\begin{align}
				\sup_{\tilde{\bm{q}}} L_{f,\tilde{\bm{p}}}(\tilde{\bm{q}})\leq L_{f,\tilde{\bm{p}}}
				\quad&\text{and}\quad
				\sup_{\tilde{\bm{q}}} L_{g,\tilde{\bm{p}}}(\tilde{\bm{q}})\leq L_{g,\tilde{\bm{p}}},
				\label{eq_PGO_property3.1}
			\end{align}
			and if $\tilde{\bm{p}}$ lies in a bounded subset, we also have
			\begin{align}
				\sup_{\tilde{\bm{p}}} L_{f,\tilde{\bm{q}}}(\tilde{\bm{p}})\leq L_{f,\tilde{\bm{q}}}
				\quad&\text{and}\quad
				\sup_{\tilde{\bm{p}}} L_{g,\tilde{\bm{q}}}(\tilde{\bm{p}})\leq L_{g,\tilde{\bm{q}}}.
				\label{eq_PGO_property3.2}
			\end{align}
			
			\item[(d).] The gradient of $f(\tilde{\bm{p}},\tilde{\bm{q}},\bm{t})$ is Lipschitz continuous on bounded subset of $\mathbb{R}^{4n}\times\mathbb{R}^{4n}\times \mathbb{R}^{3n}$ with Lipschitz constant $L_{f}>0$, i.e., for any $(\tilde{\bm{p}}_{1},\tilde{\bm{q}}_{1},\bm{t}_{1})$ and $(\tilde{\bm{p}}_{2},\tilde{\bm{q}}_{2},\bm{t}_{2}) \in 
			\mathbb{R}^{4n}\times\mathbb{R}^{4n}\times \mathbb{R}^{3n}$, it holds that 
			\begin{equation}
				\|\nabla f(\tilde{\bm{p}}_{1},\tilde{\bm{q}}_{1},\bm{t}_{1})
				-\nabla f(\tilde{\bm{p}}_{2},\tilde{\bm{q}}_{2},\bm{t}_{2})\| 
				\leq L_{f}\|(\tilde{\bm{p}}_{1}-\tilde{\bm{p}}_{2},
				\tilde{\bm{q}}_{1}-\tilde{\bm{q}}_{2},
				\bm{t}_{1}-\bm{t}_{2})\|.
				\label{eq_PGO_property4}
			\end{equation}
			Similarly, the gradient of $g(\tilde{\bm{p}},\tilde{\bm{q}})$ is Lipschitz continuous with Lipschitz constant $L_{g}>0$.
		
		\end{enumerate}
	\label{thm_PGO_property}
\end{proposition}
\begin{proof}
		(a). Due to the non-negativity of $\ell_2$-norm, we have $f^{*}\geq0,~g^{*}\geq0$.
		
		(b). We can rewrite the function $f$ as
		\begin{align}
			f(\tilde{\bm{p}},\tilde{\bm{q}},\bm{t})
			&=\sum_{i=1}^{n} \sum_{(i,j)\in \mathcal{E}_{i}^{out}} 
			\|M(\tilde{q}_{i})M(\tilde{t}_{ij})D\tilde{p}_{i}-(\tilde{t}_j - \tilde{t}_i) \|_{\Sigma_{1}}^2.\nonumber
		\end{align}
		Then, we have
		\begin{align}
			\nabla_{\tilde{p}_{i}}f(\tilde{\bm{p}},\tilde{\bm{q}},\bm{t})
			=\sum_{(i,j)\in \mathcal{E}_{i}^{out}} 
			(M(\tilde{q}_{i})M(\tilde{t}_{ij})D)^{\top}\Sigma_{1}
			\left[  M(\tilde{q}_{i})M(\tilde{t}_{ij})D\tilde{p}_{i}-(\tilde{t}_j - \tilde{t}_i)\right],  \nonumber
		\end{align}
	and
		\begin{align}
			\nabla_{\tilde{p}_{i}}^{2}f(\tilde{\bm{p}},\tilde{\bm{q}},\bm{t})
			=\sum_{(i,j)\in \mathcal{E}_{i}^{out}}  (M(\tilde{q}_{i})M(\tilde{t}_{ij})D)^{\top}\Sigma_{1}M(\tilde{q}_{i})M(\tilde{t}_{ij})D.\nonumber
		\end{align}
	It follows from Lemma \ref{lem-MW}, and we have
		\begin{align}
			\sigma_{\max}(\nabla_{\tilde{p}_{i}}^{2}f(\tilde{\bm{p}},\tilde{\bm{q}},\bm{t}))
			\leq \sigma_{\max}(\Sigma_{1})
			\sum_{(i,j)\in \mathcal{E}_{i}^{out}} \|\tilde{q}_{i}\|^{2}\|\tilde{t}_{ij}\|^{2}
			:=L_{f,\tilde{p}_{i}}(\tilde{\bm{q}}).
			\label{eq-thm5.2-1}
		\end{align}
	Let $L_{f,\tilde{\bm{p}}}(\tilde{\bm{q}}):=\max_{i}\{L_{f,\tilde{p}_{i}}(\tilde{\bm{q}})\}$, it is trivial that 
		\begin{align}
			\sigma_{\max}(\nabla_{\tilde{\bm{p}}}^{2}f(\tilde{\bm{p}},\tilde{\bm{q}},\bm{t}))
			\leq L_{f,\tilde{\bm{p}}}(\tilde{\bm{q}}).
			\nonumber
		\end{align} 
	which implies that our result holds.
	From \eqref{subproblem_t_2-1}, we obtain
	\begin{align}
		f(\tilde{\bm{p}},\tilde{\bm{q}},\bm{t})
		=\|Q\bm{t}-\bm{s}\|_{I_{m}\otimes {\hat{\Sigma}_{1}}}^{2}+c(\tilde{\bm{p}},\tilde{\bm{q}}),
		\nonumber
	\end{align}
	and $Q$ does not depend on $\tilde{\bm{p}}$, $\tilde{\bm{q}}$. The globally Lipschitz constant of $\nabla_{\bm{t}}f(\tilde{\bm{p}},\tilde{\bm{q}},\bm{t})$ is not hard to verify.
	The other globally Lipschitz constants are proven analogously.
	
	(c). From \eqref{eq-thm5.2-1}, the result can be founded obviously.
	
	(d). The result comes from that $\nabla f(\tilde{\bm{p}},\tilde{\bm{q}},\bm{t})$ is continuous and well defined on the closure of any bounded subset. The proof is completed.
	\end{proof}
	
	\medskip
	Before presenting the main results, we show the first-order optimality conditions of each subproblem, which is fundamental to the following analysis:
	\begin{subequations}
		\begin{numcases}{}
			0=\operatorname{dist}\left( \operatorname{Proj}_{T_{\tilde{\bm{p}}^{k+1}}\mathbb{U}^{n}}
			\left( 
			\nabla_{\tilde{\bm{p}}}f(\tilde{\bm{p}}^{k},\tilde{\bm{q}}^{k},\bm{t}^{k})
			+\nabla_{\tilde{\bm{p}}}g(\tilde{\bm{p}}^{k},\tilde{\bm{q}}^{k})
			-\bm{\lambda}^{k}+\beta(\tilde{\bm{p}}^{k+1}-\tilde{\bm{q}}^{k})
			\right.\right.  \nonumber\\
			 \qquad \left. \left. 
			+H_{1}(\tilde{\bm{p}}^{k+1}-\tilde{\bm{p}}^{k})\right),
			-N_{\mathbb{U}^{n}}(\tilde{\bm{p}}^{k+1})
			\right), \label{eq-subsubgradient-a}\\
			0=\nabla_{\tilde{\bm{q}}}f(\tilde{\bm{p}}^{k+1},\tilde{\bm{q}}^{k+1},\bm{t}^{k})
			+\nabla_{\tilde{\bm{q}}}g(\tilde{\bm{p}}^{k+1},\tilde{\bm{q}}^{k+1})
			+\bm{\lambda}^{k}-\beta(\tilde{\bm{p}}^{k+1}-\tilde{\bm{q}}^{k+1})
			+H_{2}(\tilde{\bm{q}}^{k+1}-\tilde{\bm{q}}^{k}),  
			\label{eq-subsubgradient-b}\\
			0=\nabla_{\bm{t}}f(\tilde{\bm{p}}^{k+1},\tilde{\bm{q}}^{k+1},\bm{t}^{k+1})
			+H_{3}(\bm{t}^{k+1}-\bm{t}^{k}).
			\label{eq-subsubgradient-c}
		\end{numcases} 
	\end{subequations}
	
	Then, we estimate the upper bound of iterative residuals of dual variable in the following Lemma. 
	\begin{lemma}
		Let $\{(\tilde{\bm{p}}^{k},\tilde{\bm{q}}^{k},\bm{t}^{k}, \bm{\lambda}^{k})\}$ be the sequence generated by Algorithm \ref{alg-PieADMM1} which is assumed to be bounded, then
		\begin{align}
			\|\bm{\lambda}^{k+1}-\bm{\lambda}^{k}\|
			&\leq 4(L_{f}^{2}+L_{g}^{2})\|\tilde{\bm{p}}^{k+1}-\tilde{\bm{p}}^{k}\|^{2}
			+4L_{f}^{2}\|\bm{t}^{k}-\bm{t}^{k-1}\|^{2}\nonumber\\
			& \quad +4\|\tilde{\bm{q}}^{k+1}-\tilde{\bm{q}}^{k}\|_{(L_{f}^{2}+L_{g}^{2})I +H_{2}^{\top}H_{2}}^{2}
			+4\|\tilde{\bm{q}}^{k}-\tilde{\bm{q}}^{k-1}\|_{H_{2}^{\top}H_{2}}^{2}.	
			\label{eq-bounded lambda}
		\end{align}
		\label{lem-bounded lambda}
	\end{lemma}
	\begin{proof}
		Combining \eqref{eq-subsubgradient-b} and \eqref{subproblem_lambda_1}, we have 
		\begin{equation}
			\bm{\lambda}^{k+1}=
			-\nabla_{\tilde{\bm{q}}}f(\tilde{\bm{p}}^{k+1},\tilde{\bm{q}}^{k+1},\bm{t}^{k})
			-\nabla_{\tilde{\bm{q}}}g(\tilde{\bm{p}}^{k+1},\tilde{\bm{q}}^{k+1})
			-H_{2}(\tilde{\bm{q}}^{k+1}-\tilde{\bm{q}}^{k}).
			\label{eq-lem5.1-1}
		\end{equation}
	Hence,
		\begin{align}
			&\|\bm{\lambda}^{k+1}-\bm{\lambda}^{k}\|^{2}\nonumber\\
			&\leq 4\|\nabla_{\tilde{\bm{q}}}f(\tilde{\bm{p}}^{k+1},\tilde{\bm{q}}^{k+1},\bm{t}^{k})
			-\nabla_{\tilde{\bm{q}}}f(\tilde{\bm{p}}^{k},\tilde{\bm{q}}^{k},\bm{t}^{k-1})\|^{2}
			+4\|\nabla_{\tilde{\bm{q}}}g(\tilde{\bm{p}}^{k+1},\tilde{\bm{q}}^{k+1})
			-\nabla_{\tilde{\bm{q}}}g(\tilde{\bm{p}}^{k},\tilde{\bm{q}}^{k})\|^{2}\nonumber\\
			& \quad 
			+4\|H_{2}(\tilde{\bm{q}}^{k+1}-\tilde{\bm{q}}^{k})\|^{2}
			+4\|H_{2}(\tilde{\bm{q}}^{k}-\tilde{\bm{q}}^{k-1})\|^{2}\nonumber\\
			&\leq 4(L_{f}^{2}+L_{g}^{2})\|\tilde{\bm{p}}^{k+1}-\tilde{\bm{p}}^{k}\|^{2}
			+4L_{f}^{2}\|\bm{t}^{k}-\bm{t}^{k-1}\|^{2}\nonumber\\
			& \quad +4\|\tilde{\bm{q}}^{k+1}-\tilde{\bm{q}}^{k}\|_{(L_{f}^{2}+L_{g}^{2})I +H_{2}^{\top}H_{2}}^{2}
			+4\|\tilde{\bm{q}}^{k}-\tilde{\bm{q}}^{k-1}\|_{H_{2}^{\top}H_{2}}^{2},\nonumber	
		\end{align}
	where the last inequality follows from \eqref{eq_PGO_property4}.
	\end{proof}

	Now we define the following merit function, which will play a crucial role in our analysis:
	\begin{equation}
		\Phi^{k} = \mathcal{L}_{\beta}
		(\tilde{\bm{p}}^{k},\tilde{\bm{q}}^{k},\bm{t}^{k},\bm{\lambda}^{k})
		+\|\tilde{\bm{q}}^{k}-\tilde{\bm{q}}^{k-1}\|_{\frac{4}{\beta}H_{2}^{\top}H_{2}}^2
		+\|\bm{t}^{k}-\bm{t}^{k-1}\|_{\frac{4}{\beta}L_{f}^2I}^2.
		\label{eq-def-Phi}
	\end{equation}
	For the ease of analysis, we also define
	\begin{align}
		P_{1}:&=\frac{1}{2}(H_1-L_{f,\tilde{\bm{p}}}I-L_{g,\tilde{\bm{p}}}I)   -\frac{4}{\beta}(L_{f}^{2}+L_{g}^2)I\nonumber\\
		P_{2}:&=\frac{1}{2}H_2-\frac{4}{\beta}(L_{f}^{2}+L_{g}^2)I-\frac{8}{\beta}H_{2}^{\top}H_{2}\nonumber\\
		P_{3}:&=\frac{1}{2}H_3-\frac{4}{\beta}L_{f}^2I\nonumber
	\end{align}
	and
	\begin{align}
		P:=\left( {\begin{array}{*{20}{c}}
				{{P_1}}&0&0\\
				0&{{P_2}}&0\\
				0&0&{{P_3}}
		\end{array}} \right),
	\qquad 
	\bm{w}:=\left( {\begin{array}{*{20}{c}}
			\tilde{\bm{p}}\\\tilde{\bm{q}}\\\bm{t}
	\end{array}} \right)
	\nonumber
	\end{align}
	Next, we prove that $\{\Phi^{k}\}$ is bounded from below in Lemma \ref{lem-lower bound} and monotonically nonincreasing in Lemma \ref{lem-nonincreasing}.
	\begin{lemma}
		Let $\{(\tilde{\bm{p}}^{k},\tilde{\bm{q}}^{k},\bm{t}^{k}, \bm{\lambda}^{k})\}$ be the sequence generated by Algorithm \ref{alg-PieADMM1} which is assumed to be bounded. If $\beta>\frac{8}{7}(L_{f,\tilde{\bm{q}}}+L_{g,\tilde{\bm{q}}})$, then $\Phi^{k+1}$ is bounded from below, i.e.,  
		\begin{equation}
			\Phi^{k+1}\geq f^{*}+g^{*},
			\quad \forall~ k\geq 0.
			\label{eq-lower bound}
		\end{equation}
		\label{lem-lower bound}
	\end{lemma}
	\begin{proof}
		It follows from \eqref{eq-lem5.1-1}, we have
		\begin{align}
			\Phi^{k+1}&=
			f(\tilde{\bm{p}}^{k+1},\tilde{\bm{q}}^{k+1},\bm{t}^{k+1})
			+g(\tilde{\bm{p}}^{k+1},\tilde{\bm{q}}^{k+1})
			-\left\langle\bm{\lambda}^{k+1},\tilde{\bm{p}}^{k+1}-\tilde{\bm{q}}^{k+1}\right\rangle
			\nonumber\\
			&\quad +\frac{\beta}{2}\|\tilde{\bm{p}}^{k+1}-\tilde{\bm{q}}^{k+1}\|^{2}
			+\|\tilde{\bm{q}}^{k+1}-\tilde{\bm{q}}^{k}\|_{\frac{4}{\beta}H_{2}^{\top}H_{2}}^2
			+\|\bm{t}^{k+1}-\bm{t}^{k}\|_{\frac{4}{\beta}L_{f}^2I}^2\nonumber\\
			&=f(\tilde{\bm{p}}^{k+1},\tilde{\bm{q}}^{k+1},\bm{t}^{k+1})
			+g(\tilde{\bm{p}}^{k+1},\tilde{\bm{q}}^{k+1})
			+\frac{\beta}{2}\|\tilde{\bm{p}}^{k+1}-\tilde{\bm{q}}^{k+1}\|^{2}
			\nonumber\\
			&\quad +\left\langle			
			\nabla_{\tilde{\bm{q}}}f(\tilde{\bm{p}}^{k+1},\tilde{\bm{q}}^{k+1},\bm{t}^{k+1})
			+\nabla_{\tilde{\bm{q}}}g(\tilde{\bm{p}}^{k+1},\tilde{\bm{q}}^{k+1}),
			\tilde{\bm{p}}^{k+1}-\tilde{\bm{q}}^{k+1}\right\rangle
			\nonumber\\
			&\quad +\left\langle
			\nabla_{\tilde{\bm{q}}}f(\tilde{\bm{p}}^{k+1},\tilde{\bm{q}}^{k+1},\bm{t}^{k})
			-\nabla_{\tilde{\bm{q}}}f(\tilde{\bm{p}}^{k+1},\tilde{\bm{q}}^{k+1},\bm{t}^{k+1}),
			\tilde{\bm{p}}^{k+1}-\tilde{\bm{q}}^{k+1}
			\right\rangle\nonumber\\
			&\quad
			+\left\langle H_{2}(\tilde{\bm{q}}^{k+1}-\tilde{\bm{q}}^{k}),\tilde{\bm{p}}^{k+1}-\tilde{\bm{q}}^{k+1} \right\rangle 
			+\|\tilde{\bm{q}}^{k+1}-\tilde{\bm{q}}^{k}\|_{\frac{4}{\beta}H_{2}^{\top}H_{2}}^2
			+\|\bm{t}^{k+1}-\bm{t}^{k}\|_{\frac{4}{\beta}L_{f}^2I}^2.
		\end{align}
		Since 
		\begin{align}
			&\left\langle 
			\nabla_{\tilde{\bm{q}}}f(\tilde{\bm{p}}^{k+1},\tilde{\bm{q}}^{k+1},\bm{t}^{k})
			-\nabla_{\tilde{\bm{q}}}f(\tilde{\bm{p}}^{k+1},\tilde{\bm{q}}^{k+1},\bm{t}^{k+1}),
			\tilde{\bm{p}}^{k+1}-\tilde{\bm{q}}^{k+1}
			\right\rangle \nonumber\\
			&\geq -\frac{4}{\beta}\|\nabla_{\tilde{\bm{q}}}f(\tilde{\bm{p}}^{k+1},\tilde{\bm{q}}^{k+1},\bm{t}^{k})
			-\nabla_{\tilde{\bm{q}}}f(\tilde{\bm{p}}^{k+1},\tilde{\bm{q}}^{k+1},\bm{t}^{k+1})\|^{2}
			-\frac{\beta}{16}\|\tilde{\bm{p}}^{k+1}-\tilde{\bm{q}}^{k+1}\|^{2}
			\nonumber\\
			&\geq -\frac{4L_{f}^2}{\beta}\|\bm{t}^{k+1}-\bm{t}^{k}\|^{2}
			-\frac{\beta}{16}\|\tilde{\bm{p}}^{k+1}-\tilde{\bm{q}}^{k+1}\|^{2},
		\end{align}
		and
		\begin{align}
			&\left\langle H_{2}(\tilde{\bm{q}}^{k+1}-\tilde{\bm{q}}^{k}),\tilde{\bm{p}}^{k+1}-\tilde{\bm{q}}^{k+1} \right\rangle \nonumber\\
			&\geq -\frac{4}{\beta}\|\tilde{\bm{q}}^{k+1}-\tilde{\bm{q}}^{k}\|_{H_{2}^{\top}H_{2}}^2
			-\frac{\beta}{16}\|\tilde{\bm{p}}^{k+1}-\tilde{\bm{q}}^{k+1}\|^{2},
		\end{align}
		and the partial gradients $\nabla_{\tilde{\bm{q}}}f(\tilde{\bm{p}},\tilde{\bm{q}},\bm{t})
		, \nabla_{\tilde{\bm{q}}}g(\tilde{\bm{p}},\tilde{\bm{q}})$ are globally Lipschitz which satisfy \eqref{eq_PGO_property2.2}, we can get
		\begin{align}
			\Phi^{k+1}&\geq
			f(\tilde{\bm{p}}^{k+1},\tilde{\bm{p}}^{k+1},\bm{t}^{k+1})
			+g(\tilde{\bm{p}}^{k+1},\tilde{\bm{p}}^{k+1})
			+\left( \frac{3\beta}{8}-\frac{L_{f,\tilde{\bm{q}}}}{2} -\frac{L_{g,\tilde{\bm{q}}}}{2}\right) \|\tilde{\bm{p}}^{k+1}-\tilde{\bm{q}}^{k+1}\|^{2}\nonumber\\
			&\geq f^{*}+g^{*},
		\end{align}
		when $\beta>\frac{4}{3}(L_{f,\tilde{\bm{q}}}+L_{g,\tilde{\bm{q}}})$. The proof is completed.
	\end{proof}
	\begin{lemma}
		Let $\{(\tilde{\bm{p}}^{k},\tilde{\bm{q}}^{k},\bm{t}^{k}, \bm{\lambda}^{k})\}$ be the sequence generated by Algorithm \ref{alg-PieADMM1} which is assumed to be bounded and $H_{i}, i=1,2,3$, satisfy 
		\begin{equation*}
			H_{1} \succ L_{f,\tilde{\bm{p}}}I+L_{g,\tilde{\bm{p}}}I,~  H_{2}\succ0,~\text{and}~ H_{3}\succ0.		
		\end{equation*}
		If
		\begin{equation}
			\beta >\max\left\lbrace \frac{8(L_{f}^{2}+L_{g}^2)}{\sigma_{\min}(H_{1})-L_{f,\tilde{\bm{p}}}-L_{g,\tilde{\bm{p}}}},
			\frac{8(L_{f}^{2}+L_{g}^2)+16\sigma_{\max}^{2}(H_{2})}{\sigma_{\min}(H_{2})},
			\frac{8L_{f}^{2}}{\sigma_{\min}(H_{3})}
			\right\rbrace ,
			\label{eq-range of beta2}
		\end{equation}
		we have
		\begin{align}
			\Phi^{k}-\Phi^{k+1}&\geq 
			\|\bm{w}^{k+1}-\bm{w}^{k}\|_{P}^{2},
			\label{eq-nonincreasing}
		\end{align}
		and the right-hand side is non-negative. 
	\label{lem-nonincreasing}
	\end{lemma}
	\begin{proof}
		From the subproblem \eqref{subproblem_p_1}, we have
		\begin{align}
			\mathcal{L}_{\beta}
			(\tilde{\bm{p}}^{k},\tilde{\bm{q}}^{k},\bm{t}^{k},\bm{\lambda}^{k})
			&= \mathcal{L}_{\beta}^{k}
			(\tilde{\bm{p}}^{k},\tilde{\bm{q}}^{k},\bm{t}^{k},\bm{\lambda}^{k}) \nonumber\\
			&\geq \mathcal{L}_{\beta}^{k}
			(\tilde{\bm{p}}^{k+1},\tilde{\bm{q}}^{k},\bm{t}^{k},\bm{\lambda}^{k})
			+\frac{1}{2}\|\tilde{\bm{p}}^{k+1}-\tilde{\bm{p}}^{k}\|_{H_1}^2 \nonumber\\
			&\geq f(\tilde{\bm{p}}^{k+1},\tilde{\bm{q}}^{k},\bm{t}^{k})
			-\frac{L_{f,\tilde{\bm{p}}}}{2}\|\tilde{\bm{p}}^{k+1}-\tilde{\bm{p}}^{k}\|^2
			+g(\tilde{\bm{p}}^{k+1},\tilde{\bm{q}}^{k})
			-\frac{L_{g,\tilde{\bm{p}}}}{2}\|\tilde{\bm{p}}^{k+1}-\tilde{\bm{p}}^{k}\|^2 \nonumber\\
			&\quad-\left\langle
			\bm{\lambda}^{k},\tilde{\bm{p}}^{k+1}-\tilde{\bm{q}}^{k}\right\rangle
			+\frac{\beta}{2}\|\tilde{\bm{p}}^{k+1}-\tilde{\bm{q}}^{k}\|^2
			+\frac{1}{2}\|\tilde{\bm{p}}^{k+1}-\tilde{\bm{p}}^{k}\|_{H_1}^2 \nonumber\\
			&=\mathcal{L}_{\beta}
			(\tilde{\bm{p}}^{k+1},\tilde{\bm{q}}^{k},\bm{t}^{k},\bm{\lambda}^{k})
			+\frac{1}{2}\|\tilde{\bm{p}}^{k+1}-\tilde{\bm{p}}^{k}\|_
			{H_1-L_{f,\tilde{\bm{p}}}I-L_{g,\tilde{\bm{p}}}I}^2,
			\label{eq-lem5.3-1}
		\end{align}
	where the second inequality follows from \eqref{eq-Lk} and \eqref{eq_PGO_property2.2}. From the subproblem \eqref{subproblem_q_1} and \eqref{subproblem_t_1}, we obtain
	\begin{align}
		\mathcal{L}_{\beta}
		(\tilde{\bm{p}}^{k+1},\tilde{\bm{q}}^{k},\bm{t}^{k},\bm{\lambda}^{k})
		-\mathcal{L}_{\beta}
		(\tilde{\bm{p}}^{k+1},\tilde{\bm{q}}^{k+1},\bm{t}^{k+1},\bm{\lambda}^{k})
		\geq \frac{1}{2}\|\tilde{\bm{q}}^{k+1}-\tilde{\bm{q}}^{k}\|_{H_2}^2
		+\frac{1}{2}\|\bm{t}^{k+1}-\bm{t}^{k}\|_{H_3}^2.
		\label{eq-lem5.3-2}
	\end{align}
		Moreover, according to \eqref{subproblem_lambda_1}, 
		\begin{equation}
			\mathcal{L}_{\beta}
			(\tilde{\bm{p}}^{k+1},\tilde{\bm{q}}^{k+1},\bm{t}^{k+1},\bm{\lambda}^{k})
			-\mathcal{L}_{\beta}
			(\tilde{\bm{p}}^{k+1},\tilde{\bm{q}}^{k+1},\bm{t}^{k+1},\bm{\lambda}^{k+1})
			=-\frac{1}{\beta}\|\bm{\lambda}^{k+1}-\bm{\lambda}^{k}\|^2.
			\label{eq-lem5.3-3}
		\end{equation}
		Combining \eqref{eq-lem5.3-1}-\eqref{eq-lem5.3-3} and Lemma \ref{lem-bounded lambda} yields that 
		\begin{align}
			&~\mathcal{L}_{\beta}
			(\tilde{\bm{p}}^{k},\tilde{\bm{q}}^{k},\bm{t}^{k},\bm{\lambda}^{k})
			-\mathcal{L}_{\beta}
			(\tilde{\bm{p}}^{k+1},\tilde{\bm{q}}^{k+1},\bm{t}^{k+1},\bm{\lambda}^{k+1})\nonumber\\
			&\geq \frac{1}{2}\|\tilde{\bm{p}}^{k+1}-\tilde{\bm{p}}^{k}\|_{H_1-L_{f,\tilde{\bm{p}}}I-L_{g,\tilde{\bm{p}}}I}^2
			+\frac{1}{2}\|\tilde{\bm{q}}^{k+1}-\tilde{\bm{q}}^{k}\|_{H_2}^2
			+\frac{1}{2}\|\bm{t}^{k+1}-\bm{t}^{k}\|_{H_3}^2
			-\frac{1}{\beta}\|\bm{\lambda}^{k+1}-\bm{\lambda}^{k}\|^2\nonumber\\
			&\geq \|\tilde{\bm{p}}^{k+1}-\tilde{\bm{p}}^{k}\|_{\frac{1}{2}(H_1-L_{f,\tilde{\bm{p}}}I-L_{g,\tilde{\bm{p}}}I)   -\frac{4}{\beta}(L_{f}^{2}+L_{g}^2)I}^2
			+\|\bm{t}^{k+1}-\bm{t}^{k}\|_{\frac{1}{2}H_3}^2
			-\|\bm{t}^{k}-\bm{t}^{k-1}\|_{\frac{4}{\beta}L_{f}^2I}^2
			\nonumber\\
			& \quad +\|\tilde{\bm{q}}^{k+1}-\tilde{\bm{q}}^{k}\|
			_{\frac{1}{2}H_2-\frac{4}{\beta}(L_{f}^{2}+L_{g}^2)I-\frac{4}{\beta}H_{2}^{\top}H_{2}}^2
			-\|\tilde{\bm{q}}^{k}-\tilde{\bm{q}}^{k-1}\|_{\frac{4}{\beta}H_{2}^{\top}H_{2}}^{2}
			, \nonumber
		\end{align}
		which implies that
		\begin{align}
			\Phi^{k}-\Phi^{k+1}&\geq \|\tilde{\bm{p}}^{k+1}-\tilde{\bm{p}}^{k}\|_{P_{1}}^2
			+\|\tilde{\bm{q}}^{k+1}-\tilde{\bm{q}}^{k}\|_{P_{2}}^2
			+\|\bm{t}^{k+1}-\bm{t}^{k}\|_{P_{3}}^2. \nonumber
		\end{align}
		It is not hard to verify that when $\beta$ satisfies \eqref{eq-range of beta2}, $\{\Phi^{k}\}$ is monotonically nonincreasing. The proof is completed.
	\end{proof}
	Now we are ready to establish the iteration complexity of Algorithm \ref{alg-PieADMM1} for finding an $\epsilon$-stationary solution of PGO model \eqref{PGO-model2}. For the ease of analysis, we define
	\begin{align}
		\kappa_{1}&=\frac{4}{\beta^{2}}(L_{f}^{2}+L_{g}^{2}+\sigma_{\max}^{2}(H_{2})),
		\quad \kappa_{2}=\sigma_{\max}^{2}(H_{3}),\nonumber\\
		\kappa_{3}&=\max\{2L_{f}^{2},2\sigma_{\max}^{2}(H_{2})\},
		\quad \kappa_{4}=(L_{f}+L_{g}+\beta+\sigma_{\max}(H_{1}))^{2},\nonumber
	\end{align}
	and
	\begin{align}
		\nu&= \sigma_{\min}(P)
		=\min\left\lbrace \frac{1}{2}(\sigma_{\min}(H_1)-L_{f,\tilde{\bm{p}}}-L_{g,\tilde{\bm{p}}})-\frac{4}{\beta}(L_{f}^{2}+L_{g}^2),\right. 
		\nonumber\\
		&\qquad\qquad \left.  \frac{1}{2}\sigma_{\min}(H_2)-\frac{4}{\beta}(L_{f}^{2}+L_{g}^2)
		-\frac{8}{\beta}\sigma_{\max}^{2}(H_2),
		~\frac{1}{2}\sigma_{\min}(H_3)-\frac{4}{\beta}L_{f}^2
		\right\rbrace. \nonumber
	\end{align}	 
	\begin{theorem}
		Suppose that the sequence $\{(\tilde{\bm{p}}^{k},\tilde{\bm{q}}^{k},\bm{t}^{k},\bm{\lambda}^{k})\}$ is generated by Algorithm \ref{alg-PieADMM1} and $H_{i}, i=1,2,3$, satisfy 
		\begin{equation*}
			H_{1} \succ L_{f,\tilde{\bm{p}}}I+L_{g,\tilde{\bm{p}}}I,~  H_{2}\succ0,~\text{and}~ H_{3}\succ0.		
		\end{equation*}
		and $\beta$ satisfies
		\begin{equation}
			\beta >\max\left\lbrace
			\frac{4}{3}(L_{f,\tilde{\bm{q}}}+L_{g,\tilde{\bm{q}}}), \frac{8(L_{f}^{2}+L_{g}^2)}{\sigma_{\min}(H_{1})-L_{f,\tilde{\bm{p}}}-L_{g,\tilde{\bm{p}}}},
			\frac{8(L_{f}^{2}+L_{g}^2)+16\sigma_{\max}^{2}(H_{2})}{\sigma_{\min}(H_{2})},
			\frac{8L_{f}^{2}}{\sigma_{\min}(H_{3})}
			\right\rbrace .
			\label{eq-range of beta3}
		\end{equation}
		Then, we can find an $\epsilon$-stationary solution $(\tilde{\bm{p}}^{\hat{k}},\tilde{\bm{q}}^{\hat{k}},\bm{t}^{\hat{k}},\bm{\lambda}^{\hat{k}})$ of PGO model \eqref{PGO-model2}, where 
		$\bm{w}^{\hat{k}}=(\tilde{\bm{p}}^{\hat{k}},\tilde{\bm{q}}^{\hat{k}},\bm{t}^{\hat{k}})$ be the first iteration that satisfies
		\begin{align}
			\theta_{\hat{k}}:=\|\bm{w}^{\hat{k}}-\bm{w}^{\hat{k}-1}\|^{2}+\|\bm{w}^{\hat{k}-1}-\bm{w}^{\hat{k}-2}\|^{2}
			\leq \epsilon^{2}/\max\{\kappa_{1},\kappa_{2},\kappa_{3},\kappa_{4}\},\nonumber
		\end{align}
		Moreover, $\hat{k}$ is no more than 
		\begin{equation*}
			T:=\bigg\lceil\frac{2\max\{\kappa_{1},\kappa_{2},\kappa_{3},\kappa_{4}\}}{\nu\epsilon^{2}}(\Phi^{1}-f^{*}-g^{*})\bigg\rceil.
		\end{equation*}
		\label{thm-iteration complexity}
	\end{theorem}
	\begin{proof}
		By summing \eqref{eq-nonincreasing} over $k=1,2,\dots,T$, we have
		\begin{align}
			\Phi^{1}-\Phi^{T+1}
			\geq 
			\sum_{k=1}^{T}\|\bm{w}^{k+1}-\bm{w}^{k}\|_{P}^{2}
			\geq
			\nu\sum_{k=1}^{T}\|\bm{w}^{k+1}-\bm{w}^{k}\|^{2},\nonumber
		\end{align}
	which implies
	\begin{align}
		\min_{2\leq k \leq T+1}\theta_{k}
		&\leq\frac{1}{T}\sum_{k=2}^{T+1}\theta_{k}\nonumber\\
		&=\frac{1}{T}\sum_{k=2}^{T+1}\|\bm{w}^{k+1}-\bm{w}^{k}\|^{2}
		+\frac{1}{T}\sum_{k=1}^{T}\|\bm{w}^{k+1}-\bm{w}^{k}\|^{2}\nonumber\\
		&\leq\frac{1}{\nu T}(\Phi^{1}+\Phi^{2}-\Phi^{T+1}-\Phi^{T+2})
		\leq\frac{2}{\nu T}(\Phi^{1}-f^{*}-g^{*}).
		\label{eq-thm5.3-1}
	\end{align}
	By Lemma \ref{lem-bounded lambda},
	\begin{align}
		\|\tilde{\bm{p}}^{k+1}-\tilde{\bm{q}}^{k+1}\|^{2}
		&=\frac{1}{\beta^{2}}\|\bm{\lambda}^{k+1}-\bm{\lambda}^{k}\|^{2}
		\nonumber\\
		&\leq \frac{4}{\beta^{2}}
		\left\lbrace (L_{f}^{2}+L_{g}^{2})\|\tilde{\bm{p}}^{k+1}-\tilde{\bm{p}}^{k}\|^{2}
		+L_{f}^{2}\|\bm{t}^{k}-\bm{t}^{k-1}\|^{2}\right. \nonumber\\
		& \quad \left. +\|\tilde{\bm{q}}^{k+1}-\tilde{\bm{q}}^{k}\|_{(L_{f}^{2}+L_{g}^{2})I +H_{2}^{\top}H_{2}}^{2}
		+\|\tilde{\bm{q}}^{k}-\tilde{\bm{q}}^{k-1}\|_{H_{2}^{\top}H_{2}}^{2}\right\rbrace  \nonumber\\
		&\leq \kappa_{1}\theta_{k}.
		\label{eq-thm5.3-2}
	\end{align} 
	According to \eqref{eq-subsubgradient-c},
	\begin{align}
		\|\nabla_{\bm{t}}f(\tilde{\bm{p}}^{k+1},\tilde{\bm{q}}^{k+1},\bm{t}^{k+1})\|^{2}
		=\|H_{3}(\bm{t}^{k+1}-\bm{t}^{k})\|^{2}
		\leq \kappa_{2}\theta_{k}.
		\label{eq-thm5.3-3}
	\end{align}
	From \eqref{eq-subsubgradient-b}, we have
	\begin{align}
		&\|\nabla_{\tilde{\bm{q}}}f(\tilde{\bm{p}}^{k+1},\tilde{\bm{q}}^{k+1},\bm{t}^{k+1})
		+\nabla_{\tilde{\bm{q}}}g(\tilde{\bm{p}}^{k+1},\tilde{\bm{q}}^{k+1})+\bm{\lambda}^{k+1}
		\|^{2}\nonumber\\
		&\leq
		\|\nabla_{\tilde{\bm{q}}}f(\tilde{\bm{p}}^{k+1},\tilde{\bm{q}}^{k+1},\bm{t}^{k+1})
		-\nabla_{\tilde{\bm{q}}}f(\tilde{\bm{p}}^{k+1},\tilde{\bm{q}}^{k+1},\bm{t}^{k})
		-H_{2}(\tilde{\bm{q}}^{k+1}-\tilde{\bm{q}}^{k})\|^{2}\nonumber\\
		&\leq 
		2L_{f}^{2}\|\bm{t}^{k+1}-\bm{t}^{k}\|^{2}
		+2\|\tilde{\bm{q}}^{k+1}-\tilde{\bm{q}}^{k}\|_{H_{2}^{\top}H_{2}}^{2}\nonumber\\
		&\leq \kappa_{3}\theta_{k}.
		\label{eq-thm5.3-4}
	\end{align}
	From the first-order optimality conditions of $\tilde{\bm{p}}$-subproblem \eqref{eq-subsubgradient-a}, there exists some $\bm{s}^{k+1}\in N_{\mathbb{U}^{n}}(\tilde{\bm{p}}^{k+1})$ such that 
	$$
	\operatorname{Proj}_{T_{\tilde{\bm{p}}^{k+1}}\mathbb{U}^{n}}
	\left( 
	\nabla_{\tilde{\bm{p}}}f(\tilde{\bm{p}}^{k},\tilde{\bm{q}}^{k},\bm{t}^{k})
	+\nabla_{\tilde{\bm{p}}}g(\tilde{\bm{p}}^{k},\tilde{\bm{q}}^{k})
	-\bm{\lambda}^{k}+\beta(\tilde{\bm{p}}^{k+1}-\tilde{\bm{q}}^{k})
	+H_{1}(\tilde{\bm{p}}^{k+1}-\tilde{\bm{p}}^{k})\right)+\bm{s}^{k+1}=0,
	$$
	which implies that $\bm{s}^{k+1} \in T_{\tilde{\bm{p}}^{k+1}}\mathbb{U}^{n}$. Therefore, we have
	\begin{align}
		&\operatorname{dist}\left( 
		\operatorname{Proj}_{T_{\tilde{\bm{p}}^{k+1}}  \mathbb{U}^{n} }\left\lbrace 
		-\nabla_{\tilde{\bm{p}}}f(\tilde{\bm{p}}^{k+1},\tilde{\bm{q}}^{k+1},\bm{t}^{k+1})
		-\nabla_{\tilde{\bm{p}}}g(\tilde{\bm{p}}^{k+1},\tilde{\bm{q}}^{k+1})+\bm{\lambda}^{k+1}
		\right\rbrace , N_{\mathbb{U}^{n}}(\tilde{\bm{p}}^{k+1})
		\right) \nonumber\\
		&\leq \left\| \operatorname{Proj}_{T_{\tilde{\bm{p}}^{k+1}}  \mathbb{U}^{n} }\left\lbrace 
		-\nabla_{\tilde{\bm{p}}}f(\tilde{\bm{p}}^{k+1},\tilde{\bm{q}}^{k+1},\bm{t}^{k+1})
		-\nabla_{\tilde{\bm{p}}}g(\tilde{\bm{p}}^{k+1},\tilde{\bm{q}}^{k+1})+\bm{\lambda}^{k+1}
		\right\rbrace
		-\bm{s}^{k+1}\right\| \nonumber\\
		&\leq\left\|
		-\nabla_{\tilde{\bm{p}}}f(\tilde{\bm{p}}^{k+1},\tilde{\bm{q}}^{k+1},\bm{t}^{k+1})
		-\nabla_{\tilde{\bm{p}}}g(\tilde{\bm{p}}^{k+1},\tilde{\bm{q}}^{k+1})
		+\nabla_{\tilde{\bm{p}}}f(\tilde{\bm{p}}^{k},\tilde{\bm{q}}^{k},\bm{t}^{k})
		\right.  \nonumber\\
		&\quad \left. 
		+\nabla_{\tilde{\bm{p}}}g(\tilde{\bm{p}}^{k},\tilde{\bm{q}}^{k})
		+\beta(\tilde{\bm{q}}^{k+1}-\tilde{\bm{q}}^{k})
		+H_{1}(\tilde{\bm{p}}^{k+1}-\tilde{\bm{p}}^{k})\right\| \nonumber\\
		&\leq \left\| 
		\nabla_{\tilde{\bm{p}}}f(\tilde{\bm{p}}^{k+1},\tilde{\bm{q}}^{k+1},\bm{t}^{k+1})
		-\nabla_{\tilde{\bm{p}}}f(\tilde{\bm{p}}^{k},\tilde{\bm{q}}^{k},\bm{t}^{k})\right\| 
		+\left\| \nabla_{\tilde{\bm{p}}}g(\tilde{\bm{p}}^{k+1},\tilde{\bm{q}}^{k+1})
		-\nabla_{\tilde{\bm{p}}}g(\tilde{\bm{p}}^{k},\tilde{\bm{q}}^{k})\right\| \nonumber\\
		&\quad +\left\| \beta(\tilde{\bm{q}}^{k+1}-\tilde{\bm{q}}^{k}) \right\| 
		+\left\| H_{1}(\tilde{\bm{p}}^{k+1}-\tilde{\bm{p}}^{k}) \right\| \nonumber\\
		&\leq (L_{f}+L_{g}+\beta+\sigma_{\max}(H_{1}))\|\bm{w}^{k+1}-\bm{w}^{k}\|\nonumber\\
		&\leq \sqrt{\kappa_{4}\theta_{k}},
		\label{eq-thm5.3-5}
	\end{align}
	where the second inequality follows the nonexpansive of projection operator and \eqref{subproblem_lambda_1}.
	Combining \eqref{eq-thm5.3-1}-\eqref{eq-thm5.3-5}, we can get the upper bound of iteration complexity of Algorithm \ref{alg-PieADMM1}. The proof is completed.
\end{proof}
\begin{remark}
		Since $\mathbb{U}$ is a unit sphere, $T_{\tilde{\bm{p}}_{i}^{k+1}}\mathbb{U}$ is a linear space and $N_{\mathbb{U}}(\tilde{\bm{p}}_{i}^{k+1})$ is an orthogonal complement space of $T_{\tilde{\bm{p}}_{i}^{k+1}}\mathbb{U}$. Therefore, we have $\bm{s}^{k+1}=0$ which can help us simplify the proof in Theorem \ref{thm-iteration complexity}. Without loss of generality, the iteration complexity in Theorem \ref{thm-iteration complexity} still holds for arbitrary manifold constraints.
\end{remark}

\section{Numerical experiments} \label{sec-numerical experiments}
In this section, we evaluate the effective of PieADMM for augmented unit quaternion model \eqref{PGO-model2} on different 3D pose graph datasets. As a basis for comparison, we also evaluate the performance of the manifold-based Gauss-Newton (mG-N) method and manifold-based Levenberg-Marquardt (mL-M) method \cite{wagner2011rapid}, in which the equation is solved taking advantage of sparsity. All experiments were performed on an Intel i7-10700F CPU desktop computer with 16GB of RAM and MATLAB R2022b.

\subsection{Synthetic datasets}
We test the algorithms on two synthetic datasets: (a) Circular ring, which is a single sloop with a radius of 2 and odometric edges. The constraints of the closed loop are formed by the first point coincident with the last point. The observations are scarce, which would be challenging to restore the true poses. (b)  Cube dataset, in which the robot travels on a  $2 \times 2\times 2$ grid world and random loop closures are added between nearby nodes with probability $p_{cube}$. The total number of vertices is $n=\hat{n}^{3}$ where $\hat{n}$ is the number of nodes on each side of the cube, and
the expectation of the number of edges is $\mathbb{E}(m)=2(2\hat{n}^{3}-3\hat{n}^2+1)p_{cube} + \hat{n}^3 - 1$. As the observation probability $p_{cube}$ increases,  the recovery becomes more accurate. However, the increase in the number of edges also leads to the growth in the amount of computation when updating each vertex and $t$-subproblem.

The noisy relative pose measurements are generated by

\begin{equation*}
	\begin{aligned}
		\bm{t}_{ij}&=R_{i}^{\top}(\bm{t}_{j}-\bm{t}_{i})+\bm{t}_{\epsilon},  ~ \text{where}~ \bm{t}_{\epsilon} \sim  ~ \mathcal{N}(0,\sigma_{t}^2I_{3})\\
		\tilde{q}_{ij}&=\tilde{q}_{i}^{*}\tilde{q}_{j}\tilde{q}_{\epsilon}, \quad\quad\quad\quad~~ \text{where}~\tilde{q}_{\epsilon}\sim  ~ \text{vMF}([1,0,0,0],\frac12\sigma_{r}^2),
	\end{aligned}
\end{equation*}
where $(\tilde{q}_{i}, \bm{t}_{i})$, $i=1,2,\dots,n$, are true poses. $\sigma_{t}$ and $\sigma_{r}$ represent the noise level of translation and rotation, respectively. We measure the quality of restoration by the relative error (Rel.Err.) and Normalized Root Mean Square Error (NRMSE), which are respectively defined as
\begin{align*}
	\operatorname{Rel.Err.}  & = \frac{\|\tilde{\bm{q}}-\tilde{\bm{q}}_{0}\|+\|\bm{t}-\bm{t}_{0}\|}{\|\tilde{\bm{q}}_{0}\|+\|\bm{t}_{0}\|},\\
	\operatorname{NRMSE} & = \frac{\|\tilde{\bm{q}}-\tilde{\bm{q}}_{0}\|+\|\bm{t}-\bm{t}_{0}\|}{(\max(\bm{t})-\min(\bm{t}))\sqrt{n}},
\end{align*}
where $(\tilde{\bm{q}},\bm{t})$ is the restored pose and $(\tilde{\bm{q}}_{0},\bm{t}_{0})$ is the true pose. In order to measure the accuracy of the optimal solution obtained by PieADMM, we adopt the residual defined by 
\begin{align*}
	R^{k+1}&\!=\!\frac{1}{\beta}\|\bm{\lambda}^{k+1}\!-\!\bm{\lambda}^{k}\|^{2}
	\!+\!\beta\left( \|\tilde{\bm{q}}^{k+1}\!-\!\tilde{\bm{q}}^{k}\|^{2}
	\!+\!\|\bm{t}^{k+1}\!-\!\bm{t}^{k}\|^{2}\right) .
\end{align*}
When $R^{k+1}$ converges to zero,
\begin{align*}
	\|\tilde{\bm{p}}^{k+1}-\tilde{\bm{p}}^{k}\| 
	= &  \|\tilde{\bm{p}}^{k+1}-\tilde{\bm{q}}^{k+1}+\tilde{\bm{q}}^{k+1}-\tilde{\bm{q}}^{k}+\tilde{\bm{q}}^{k}-\tilde{\bm{p}}^{k}\| \\
	 \leq &
	\|\tilde{\bm{p}}^{k+1}-\tilde{\bm{q}}^{k+1}\|
	+\|\tilde{\bm{q}}^{k+1}-\tilde{\bm{q}}^{k}\|
	+\|\tilde{\bm{q}}^{k}-\tilde{\bm{p}}^{k}\| \\
	 = &
	\frac{1}{\beta}\left( 
	\|\bm{\lambda}^{k+1}-\bm{\lambda}^{k}\|^{2} 
	+\|\bm{\lambda}^{k}-\bm{\lambda}^{k-1}\|^{2}\right) 
	+\|\tilde{\bm{q}}^{k+1}-\tilde{\bm{q}}^{k}\| 
\end{align*}
also converges to zero.
Comparing \eqref{eq-subsubgradient-a} - \eqref{eq-subsubgradient-c}  with
Theorem \ref{thm-optimality of PGO}, we have $\{(\tilde{\bm{p}}^{k+1},\tilde{\bm{q}}^{k+1},\bm{t}^{k+1}, \bm{\lambda}^{k+1})\}$ converges to an $\epsilon$-stationary solution.

Accordingly, the stopping criterion of mG-N and mL-M is defined by relative decrease of objective function value as 
$$
R^{k+1}=\frac{F(x^{k})-F(x^{k+1})}{F(x^{k+1})}.
$$
We terminate the solvers when iteration residual $R^{k+1} < tol$ or the maximum number of iterations $MaxIter$ is reached.

\begin{figure}[t]
	\centering
	\subfloat[]{\label{circle2_re}	\includegraphics[width=0.4\linewidth]{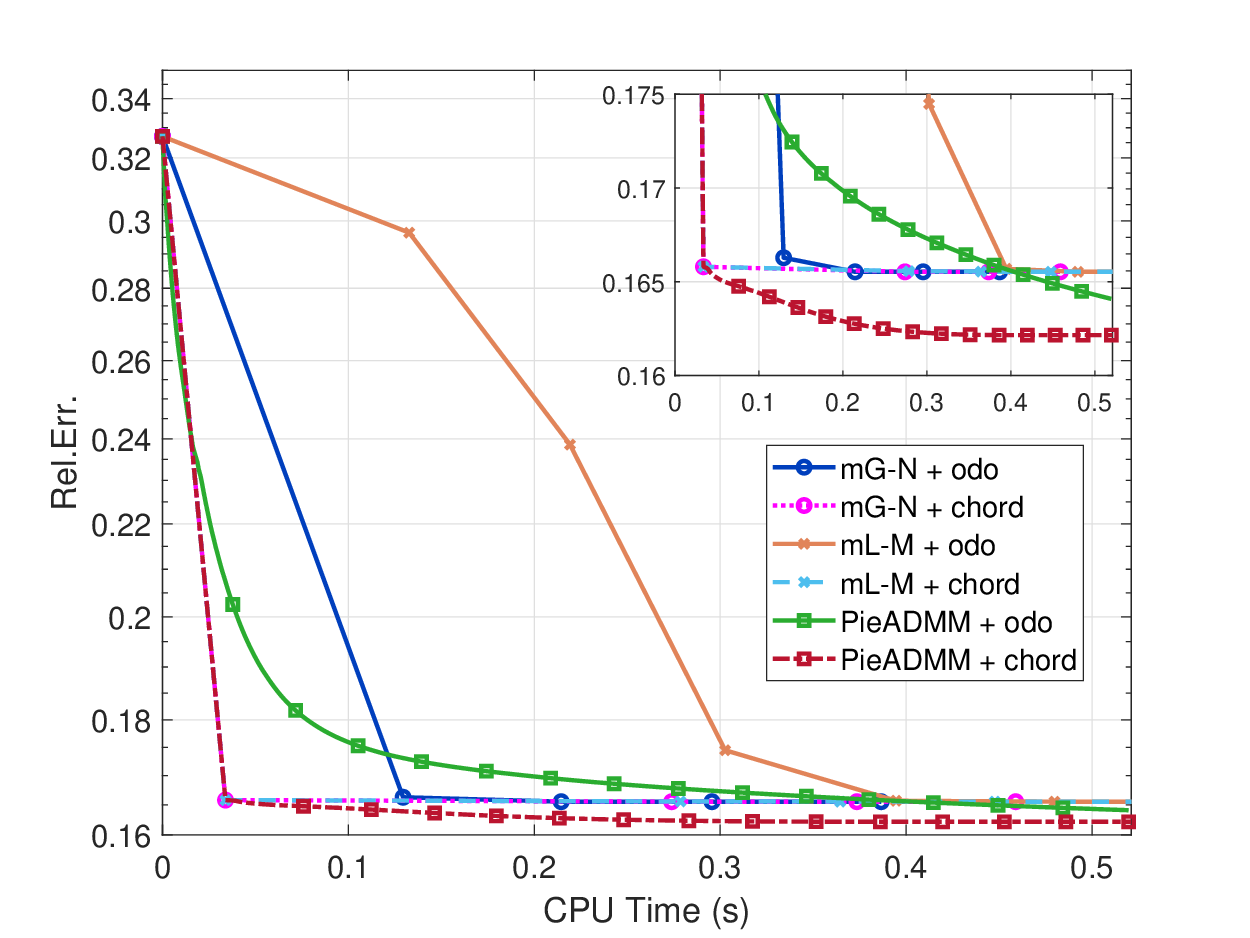}}
	\subfloat[]{\label{circle2_nrmse}	\includegraphics[width=0.4\linewidth]{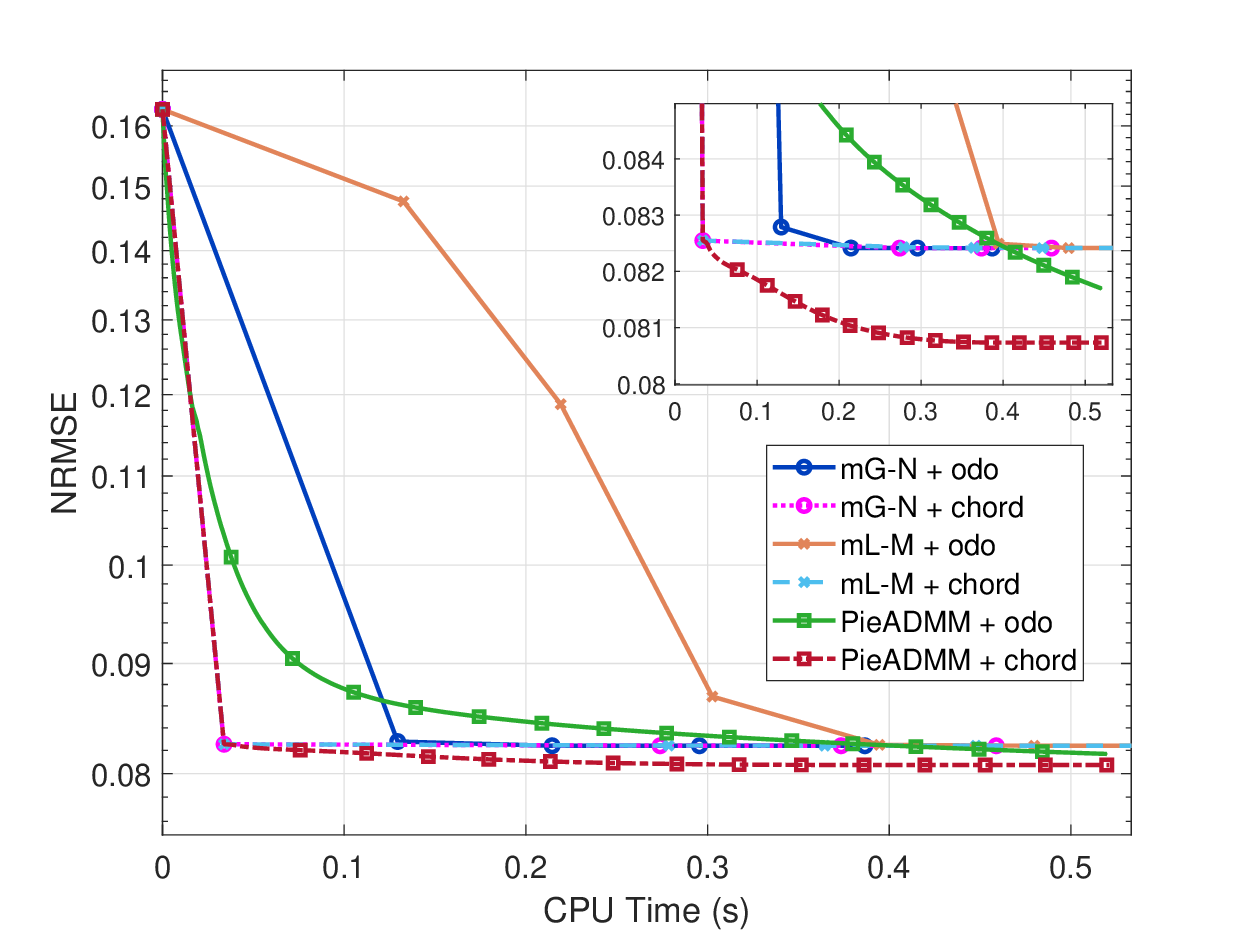}}
	
	\subfloat[]{\label{circle17_re}	\includegraphics[width=0.4\linewidth]{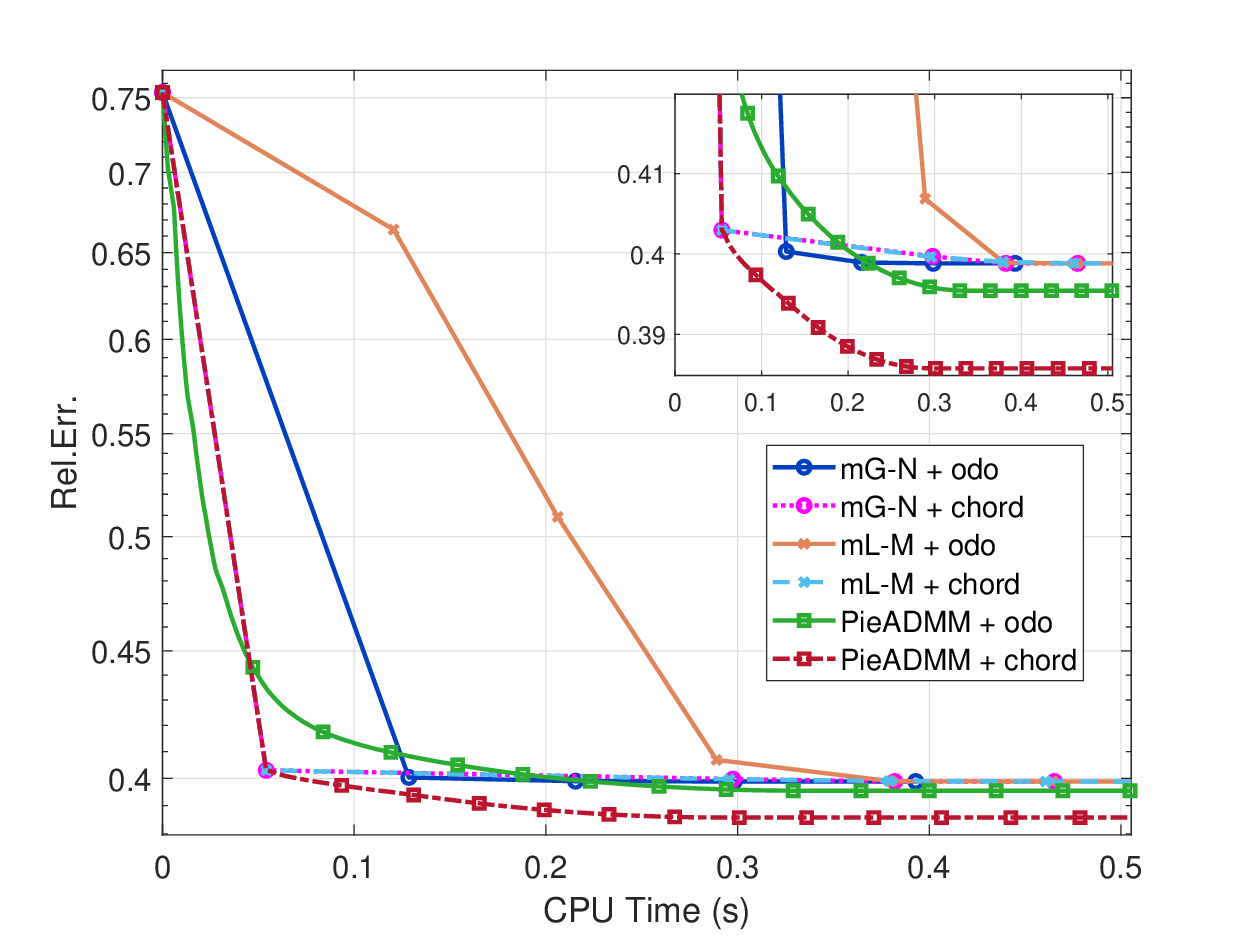}}
	\subfloat[]{\label{circle17_nrmse}	\includegraphics[width=0.4\linewidth]{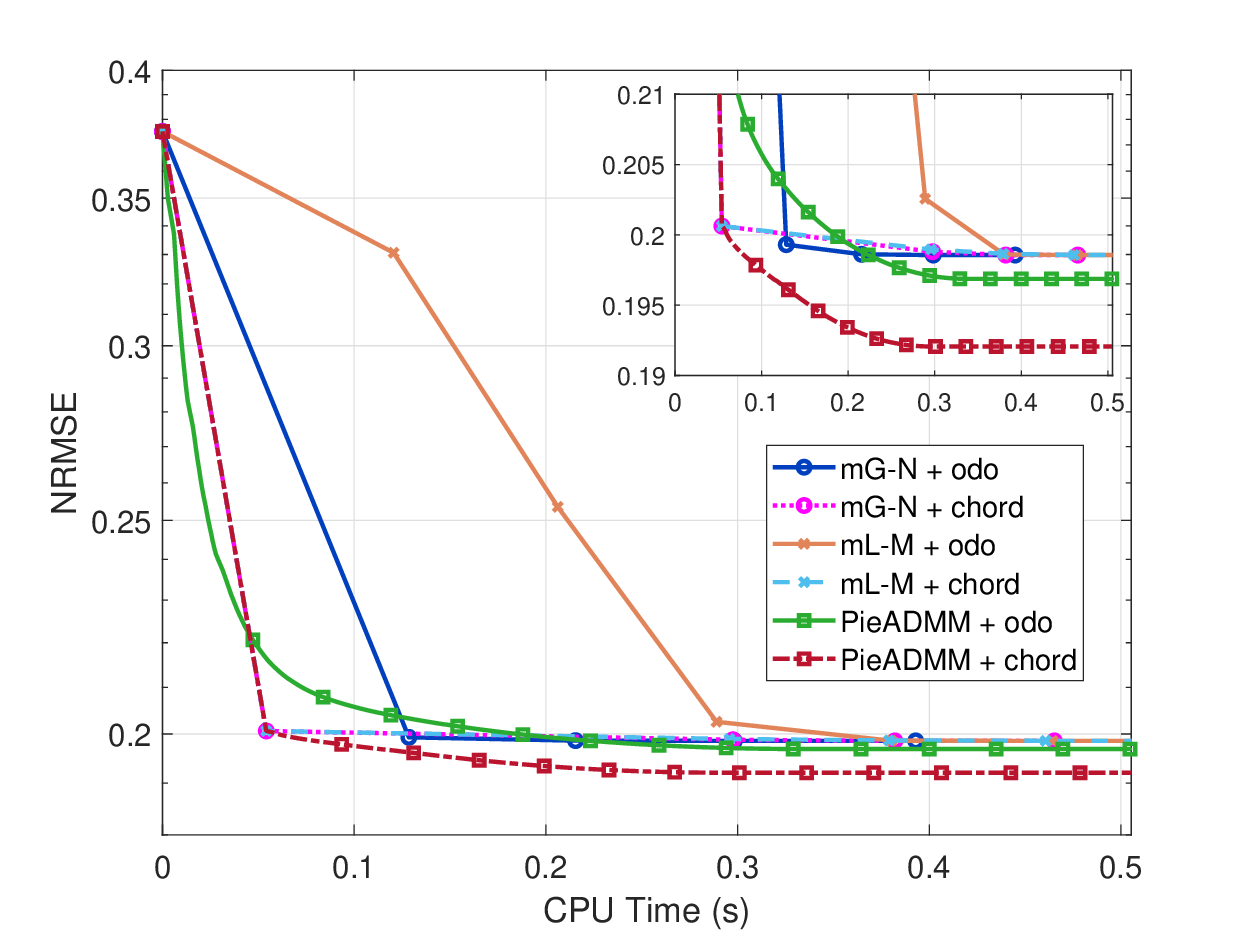}}
	\caption{Performance of the three methods versus CPU time for circular ring datasets with $n=100$, $m=100$ under different initialization techniques . The noise level of the first row is $\sigma_{r}=0.01$ and $\sigma_{t}=0.05$, and the second row is $\sigma_{r}=0.03$ and $\sigma_{t}=0.1$.}
	\label{circle_RE_NRMSE}
\end{figure}

In the experiments of synthetic datasets, we set the noise level of translation part $\sigma_{t} \in (0.01,0.3)$, and the magnitude of noise of rotation part with $\sigma_{r} \in (0.01,0.2)$. The parameters $H_{1},~H_{2},~H_{3}\succ0$ and step size $\beta>0$ satisfy Theorem \ref{thm-iteration complexity}. We also set $tol=10^{-4}$ and $MaxIter=300$ for PieADMM, and $tol=10^{-5}$ and $MaxIter=50$ for mG-N or mL-M, respectively. 
In addition, we also test odometric guess and chordal \cite{carlone2015initialization} initialization methods (henceforth referred to as `odo' and `chord') in our experiments. Without additional instructions, the default initialization is the chordal initialization. Results are averaged over $5$ runs.
\begin{figure}[t]
	\centering
	\includegraphics[width=0.5\linewidth]{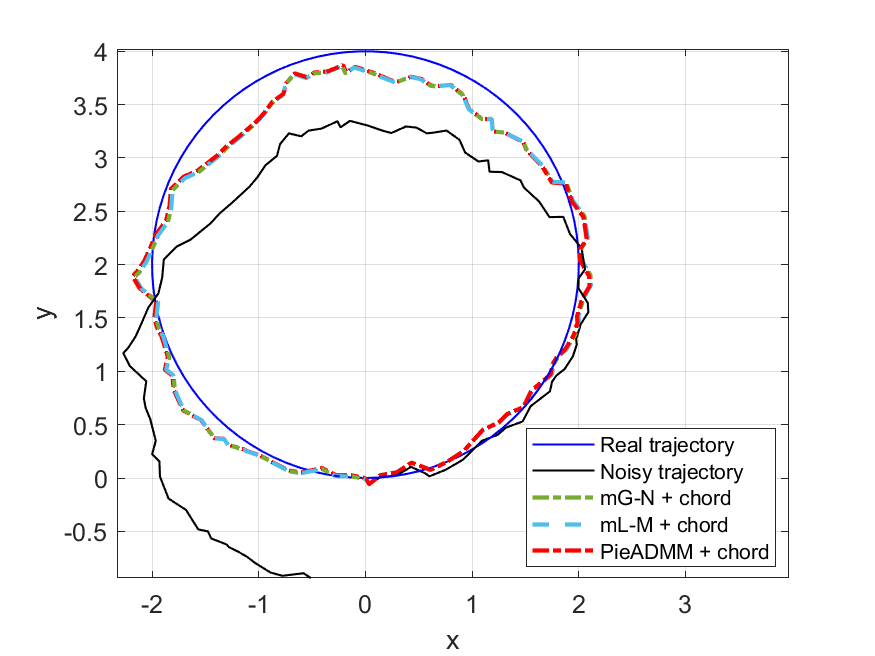}
	\caption{The trajectory of circular ring datasets with $n=100$, $m=100$, $\sigma_{r}=0.01$ and $\sigma_{t}=0.05$. The blue one is the real trajectory and the black one is the noisy trajectory. The other three dotted lines are the recovered trajectory by different algorithms.}
	\label{fig:circle2trajectory}
\end{figure}
First, we test the circular ring datasets with $n=100$, $m=100$ under different algorithms. Fig.~\ref{fig:circle2trajectory} shows the overhead view trajectory when $\sigma_{r}=0.01$, $\sigma_{t}=0.05$ and chordal initialization is adopted, and the three methods converge to the same solution in visual. We also have tested odometric guess initialization techniques. Since the recovered trajectories almost overlap and it is difficult to observe the difference, we omit them. Instead, we report the optimization process in Fig.~\ref{circle_RE_NRMSE} which records the downward trend of Rel.Err and NRMSE along with CPU time under different methods and initialization techniques. Since our PieADMM is able to update in parallel for each vertex, it can converge more quickly than others. 
Moreover, the chordal initialization can give an estimation of translation after updating the rotation, which provide a more accurate initial point than others. Under this initialization, our PieADMM can converge to a solution with lower relative error. 
The PieADMM with odometric guess initialization often not as accurate as the first several steps of mG-N, but as the iteration continues, it can achieve slightly better performance.
Thus, we take the chordal initialization as a standard initialization technique in the next experiments.

Then, we compare these algorithms under additional noise levels and list the numerical results about Rel.Err, NRMSE and CPU time in Table \ref{circle_table}. We see that PieADMM costs less time and achieves better result.

\begin{table}[H]
	\centering
	\footnotesize
	\renewcommand{\arraystretch}{1.2}
	\caption{Numerical results of different noise level of circular ring datasets  with $m=n=100$, $m=100$.}
        \resizebox{\linewidth}{!}{
	\begin{tabular}{ccccccccccc}
		\toprule[1pt]
		&    & \multicolumn{3}{c}{mG-N}   & \multicolumn{3}{c}{mL-M}   & \multicolumn{3}{c}{PieADMM} \\ 
		\cmidrule(r){3-5} \cmidrule(r){6-8} \cmidrule(r){9-11}
		$\sigma_{r}$       & $\sigma_{t}$    & Rel.Err. & NRMSE & Time\,(s) & Rel.Err. & NRMSE & Time\,(s) & Rel.Err.  & NRMSE & Time\,(s) \\
		\midrule[1pt]
		\multirow{5}{*}{0.01} & 0.01 & 0.089    & 0.0443                    & 0.31    & 0.089    & 0.0443                    & 0.52    & 0.083    & 0.0413                    & \underline{0.16}    \\
		& 0.05 & 0.166    & 0.0824                    & 0.32    & 0.166    & 0.0824                    & 0.50    & 0.164    & 0.0814                    & \underline{0.18}    \\
		& 0.1  & 0.263    & 0.1310                    & 0.32    & 0.263    & 0.1310                    & 0.50    & 0.263    & 0.1309                    & \underline{0.15}    \\
		& 0.15 & 0.361    & 0.1795                    & 0.32    & 0.361    & 0.1795                    & 0.52    & 0.357    & 0.1778                    & \underline{0.37}    \\
		& 0.2  & 0.457    & 0.2277                    & 0.40    & 0.457    & 0.2277                    & 0.50    & 0.454    & 0.2263                    & \underline{0.35}    \\
		\midrule[1pt]
		\multirow{5}{*}{0.03} & 0.01 & 0.236    & 0.1177                    & 0.32    & 0.236    & 0.1177                    & 0.50    & 0.233    & 0.1160                    & \underline{0.21}    \\
		& 0.05 & 0.306    & 0.1522                    & 0.32    & 0.306    & 0.1523                    & 0.50    & 0.304    & 0.1512                    & \underline{0.21}    \\
		& 0.1  & 0.399    & 0.1986                    & 0.32    & 0.399    & 0.1986                    & 0.51    & 0.398    & 0.1980                    & \underline{0.22}    \\
		& 0.15 & 0.492    & 0.2449                    & 0.40    & 0.492    & 0.2449                    & 0.51    & 0.489    & 0.2436                    & \underline{0.24}    \\
		& 0.2  & 0.584    & 0.2906                    & 0.40    & 0.584    & 0.2906                    & 0.59    & 0.583    & 0.2903                    & \underline{0.28}    \\
		\midrule[1pt]
		\multirow{5}{*}{0.05} & 0.01 & 0.389    & 0.1938                    & 0.32    & 0.389    & 0.1938                    & 0.50    & 0.381    & 0.1895                    & \underline{0.28}    \\
		& 0.05 & 0.453    & 0.2255                    & 0.32    & 0.453    & 0.2256                    & 0.50    & 0.446    & 0.2220                    & \underline{0.28}    \\
		& 0.1  & 0.540    & 0.2690                    & 0.40    & 0.540    & 0.2691                    & 0.50    & 0.537    & 0.2672                    & \underline{0.28}    \\
		& 0.15 & 0.628    & 0.3127                    & 0.40    & 0.628    & 0.3127                    & 0.51    & 0.626    & 0.3115                    & \underline{0.28}    \\
		& 0.2  & 0.714    & 0.3555                    & 0.40    & 0.714    & 0.3555                    & 0.58    & 0.712    & 0.3547                    & \underline{0.28}
		\\  
		\bottomrule[1pt]
	\end{tabular}}
	\label{circle_table}
\end{table}

\begin{figure}[!tbp]
	\centering
	\subfloat[]{\label{circle_n_re}	\includegraphics[width=0.3\linewidth]{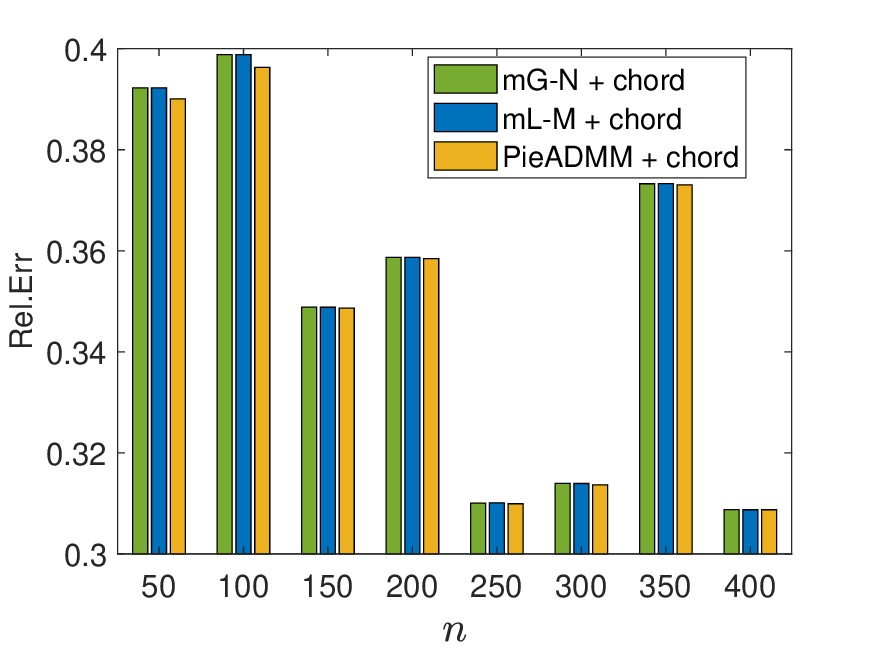}}
	\subfloat[]{\label{circle_n_nrmse}	\includegraphics[width=0.3\linewidth]{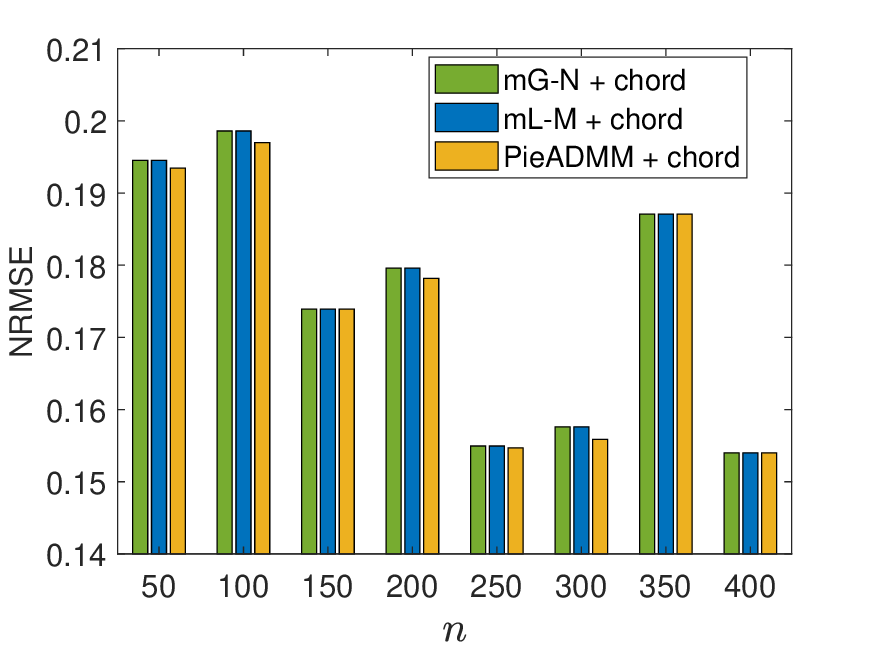}}
	\subfloat[]{\label{circle_n_time}	\includegraphics[width=0.3\linewidth]{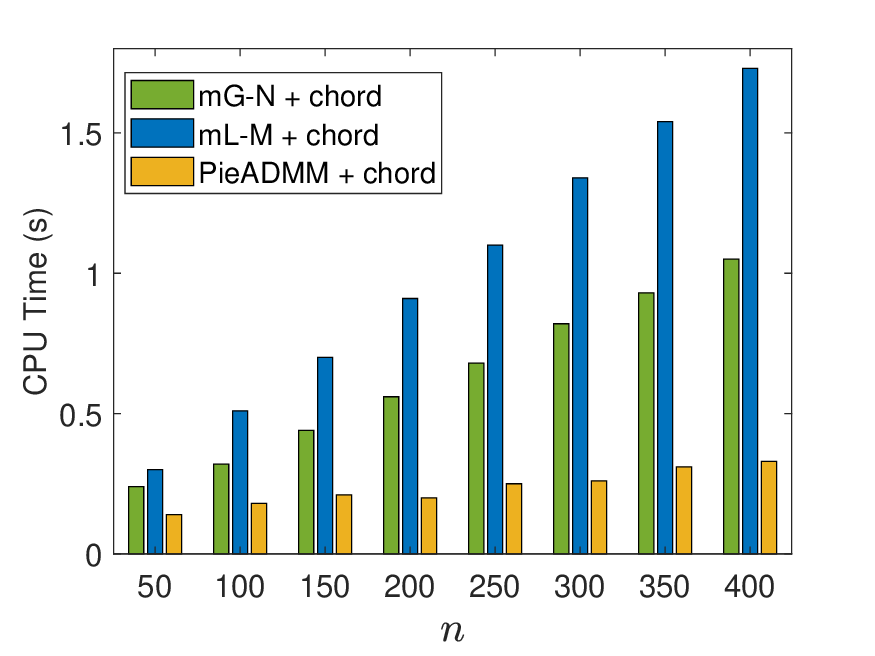}}
	\caption{Performance of the three methods for circular ring datasets under different number of poses $n$. The relative noise level is $\sigma_{r}^{rel}=0.03$ and $\sigma_{t}^{rel}=0.1$.}
	\label{circle_n}
\end{figure}

We also test the impact of the number of poses $n$.  In fact, since we limit the range of robot’s trajectory, the same level noise will cause a bigger impact when the number of vertices $n$ increase. Therefore, when comparing the influence of data size with different $n$, we use relative noise level as an unified standard, which means $\sigma_{r}=100\times\sigma_{r}^{rel}/n$ and $\sigma_{t}=100\times\sigma_{t}^{rel}/n$. The result is shown in Fig.~\ref{circle_n}. Fig.~\ref{circle_n_re} and \ref{circle_n_nrmse} show that the performance of PieADMM are flat, and sometimes slightly better, than the other two methods. However, the increasing of running time of PieADMM is much slower than them, see Fig.~\ref{circle_n_time}. It is because that the scale of $n$ almost does not affect the cost of the rotation subproblems, which can be computed in parallel. Moreover, the translation subproblem concerns only matrix multiplication, and does not depend on the inverse of the matrix.

\begin{figure}[!tbp]
	\centering
	\subfloat[Real trajectory]{\label{cube4_true_trajectory}	\includegraphics[width=0.3\linewidth]{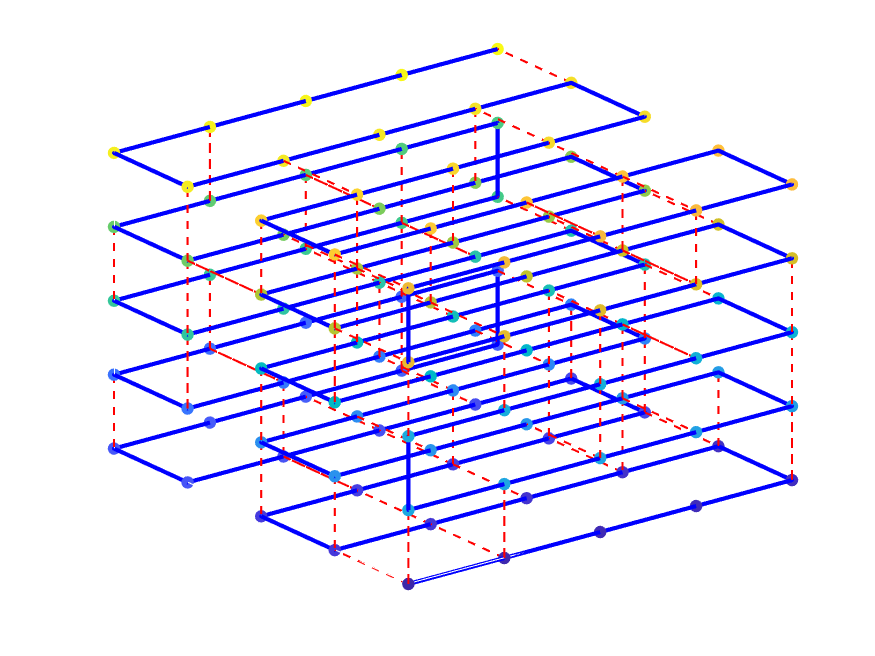}}
	\subfloat[Noisy trajectory]{\label{cube4_noisy_trajectory}	\includegraphics[width=0.3\linewidth]{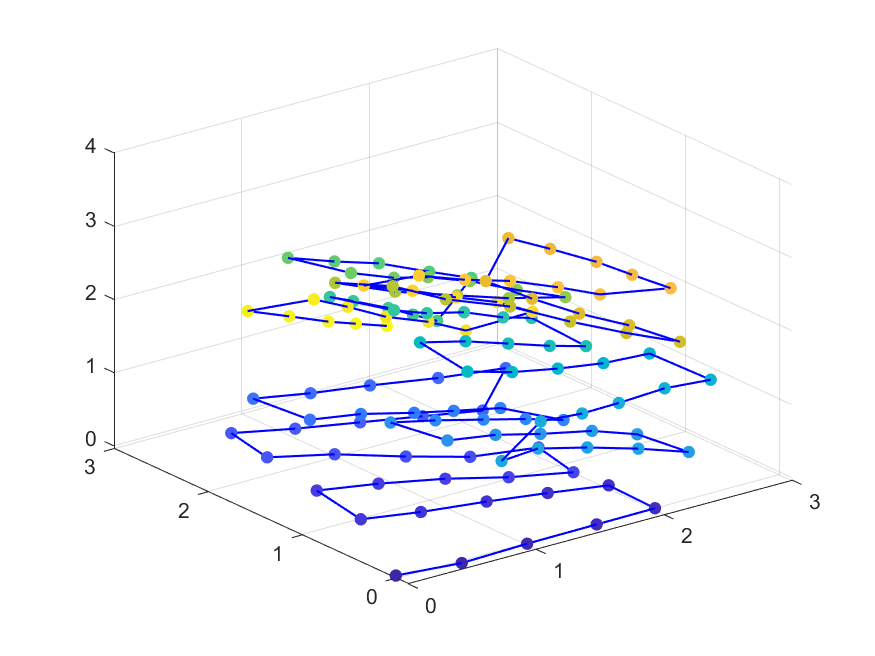}}
	\subfloat[Recovered trajectory]{\label{cube4_recovered_trajectory}	\includegraphics[width=0.3\linewidth]{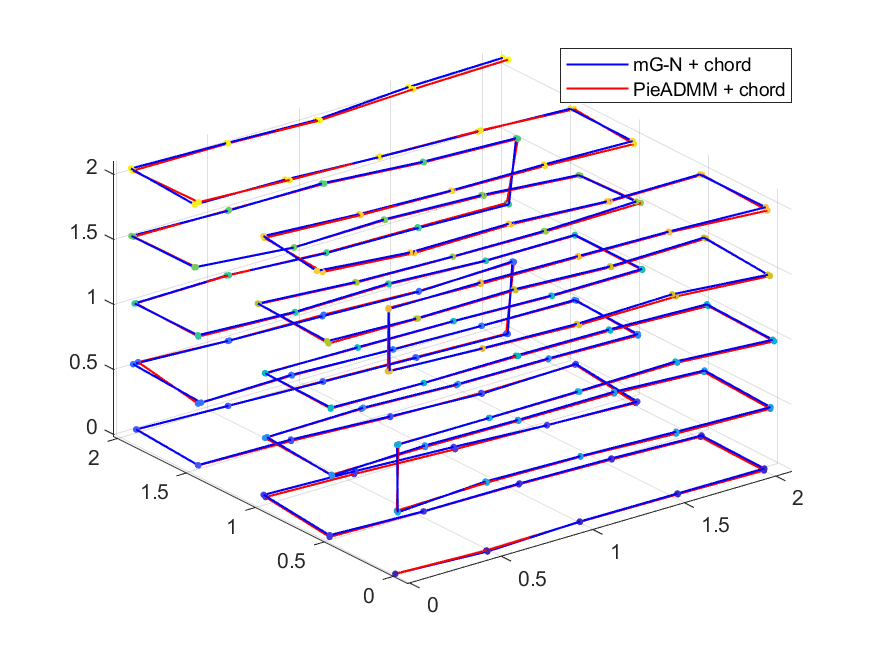}}
	
	\subfloat[Real trajectory]{\label{cube7_true_trajectory}	\includegraphics[width=0.3\linewidth]{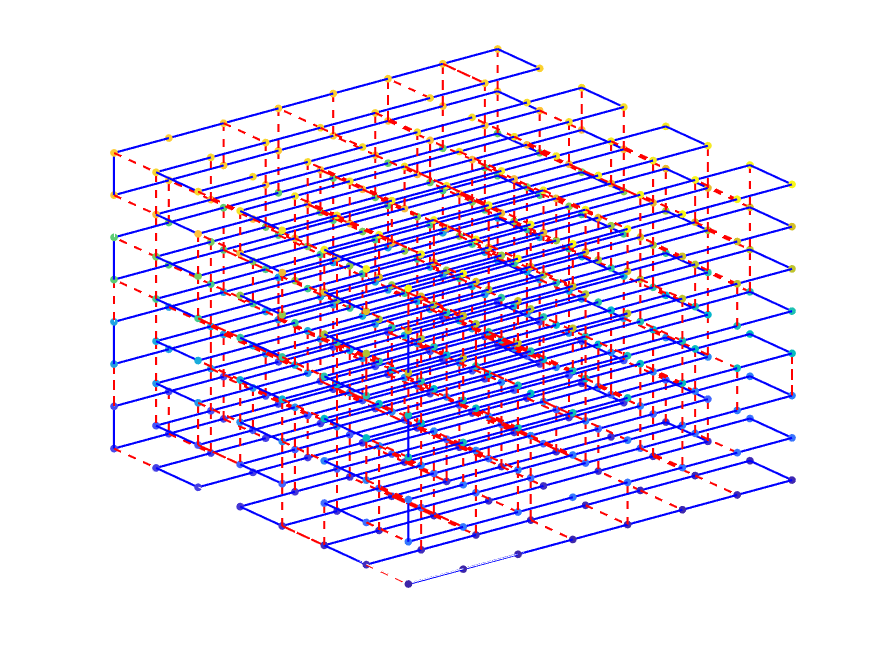}}
	\subfloat[Noisy trajectory]{\label{cube7_noisy_trajectory}	\includegraphics[width=0.3\linewidth]{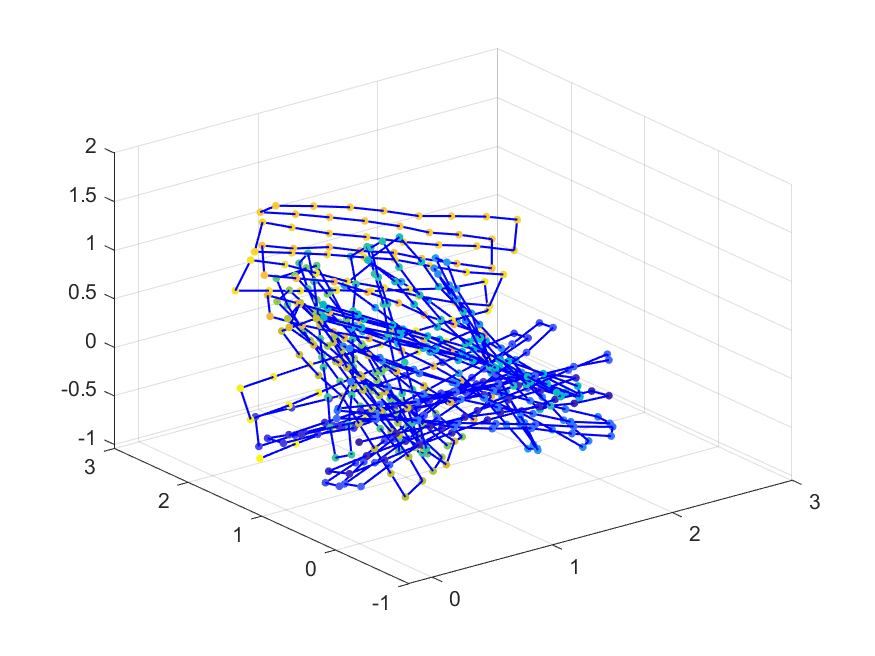}}
	\subfloat[Recovered trajectory]{\label{cube7_recovered_trajectory}	\includegraphics[width=0.3\linewidth]{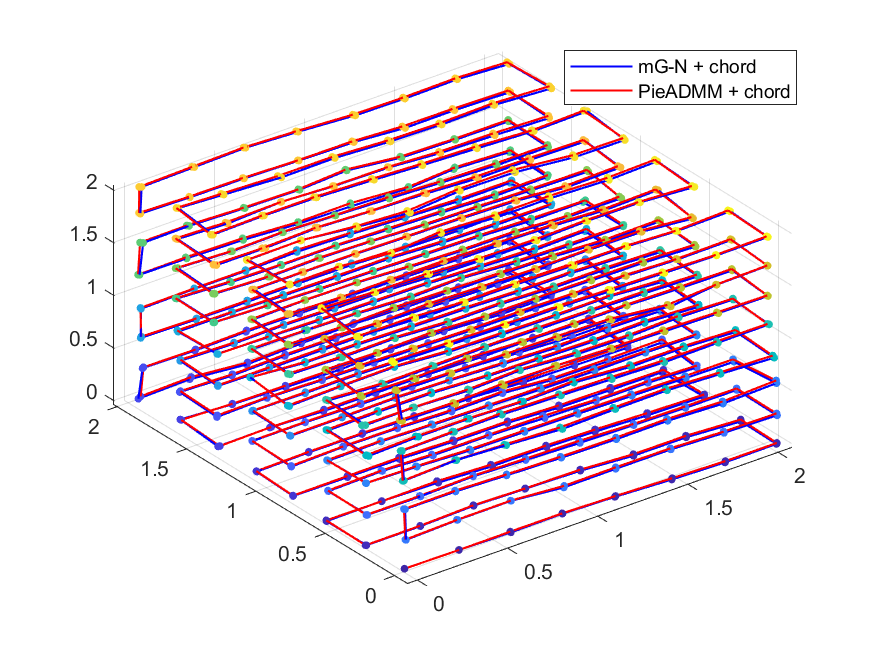}}
	\caption{The comparison of cube trajectory, where the first row is $\hat{n}=5$ and  the second row is $\hat{n}=8$, respectively. The color change of the vertices indicates the direction of the trajectory.}
	\label{cube_trajectory}
\end{figure}

\begin{figure}[!tbp]
	\centering
	\subfloat[$\hat{n}=5$]{\label{cube4_re}	\includegraphics[width=0.4\linewidth]{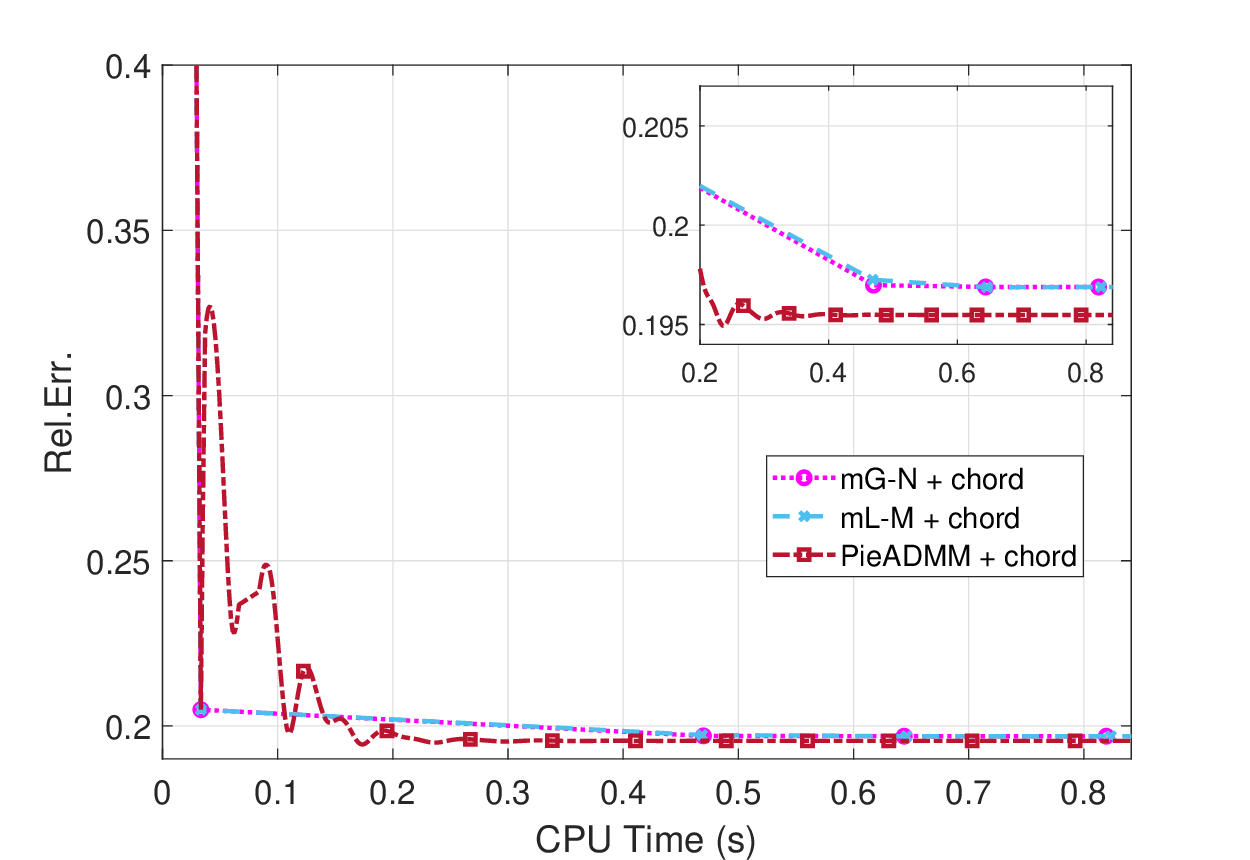}}
	\subfloat[$\hat{n}=8$]{\label{cube7_re}	\includegraphics[width=0.4\linewidth]{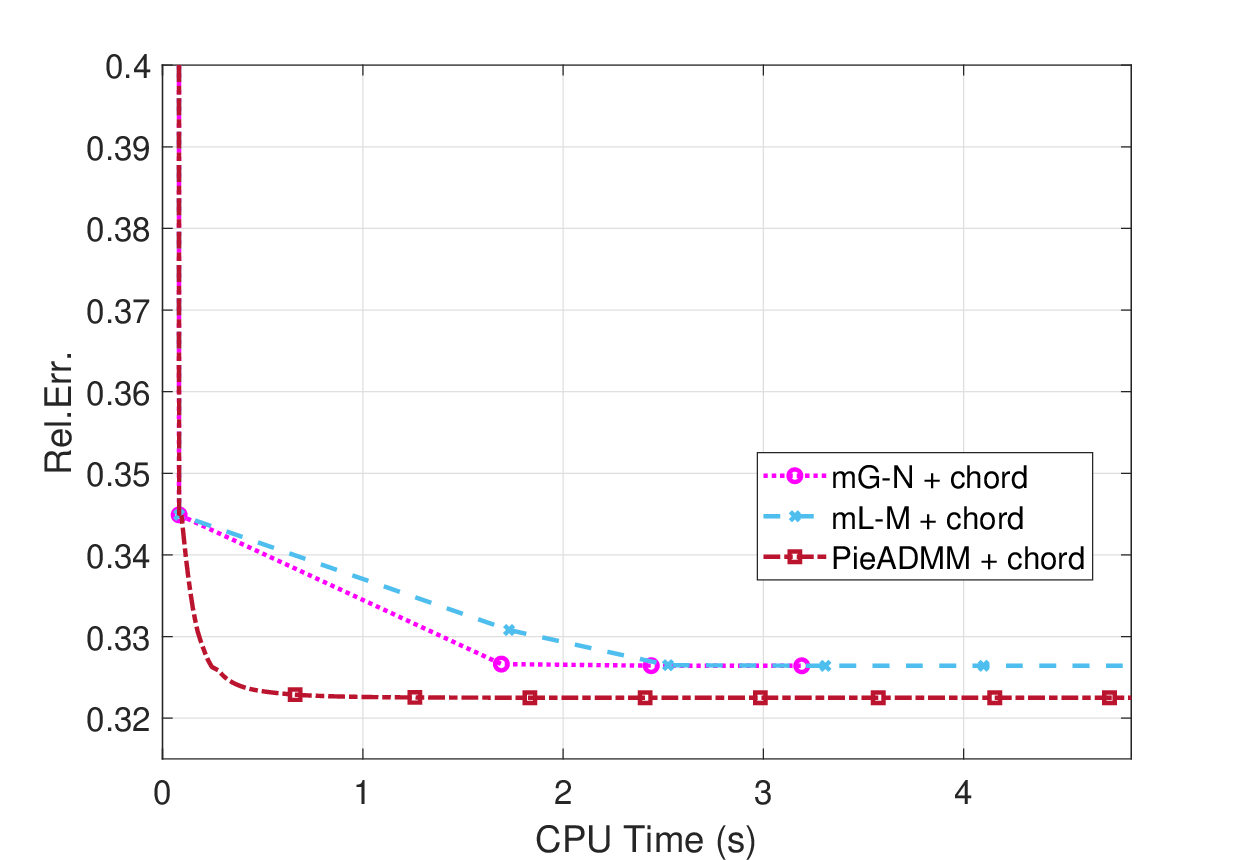}}
	\caption{Performance of the three methods versus CPU time for cube datasets with $\hat{n}=5$ or $8$.}
	\label{cube_RE}
\end{figure}

For cube datasets, let $\sigma_{t}=\sigma_{t}^{rel}/\hat{n}$ where $\sigma_{t}^{rel}$ represent the relative noise level of translation. We first consider two examples with $\hat{n}=5$ or $8$, $\sigma_{r} = 0.1,~\sigma_{t}^{rel}=0.1$ and $p_{cube}=0.3$. Fig.~\ref{cube4_true_trajectory} and \ref{cube7_true_trajectory} shows the real trajectory, in which the blue lines are produced by motions and the red dotted lines are generated by observations. Fig.~\ref{cube4_noisy_trajectory}, \ref{cube4_recovered_trajectory} and \ref{cube7_noisy_trajectory}, \ref{cube7_recovered_trajectory} are the noisy and recovered trajectory corresponding to different $\hat{n}$, respectively. The downward trends of Rel.Err along with CPU time are shown in Fig.~\ref{cube_RE}, in which we omitted the top half of the image to highlight details. Since the PGO model is non-convex, and PieADMM is a non-monotonic algorithm, the curves may oscillate. However, it always converges to the solution with higher precision in less time.

\begin{table}[!tbp]
	\centering
	\footnotesize
	\renewcommand{\arraystretch}{1.1}
	\caption{Numerical results of different $\hat{n}$ of cube datasets  with $\sigma_{r} = 0.1,~\sigma_{t}^{rel}=0.1$ and $p_{cube}=0.3$ or $0.7$.}
        \resizebox{\linewidth}{!}{
	\begin{tabular}{ccccccccccc}
		\toprule[1pt]
		&     & \multicolumn{3}{c}{mG-N}                                                               & \multicolumn{3}{c}{mL-M}                                                               & \multicolumn{3}{c}{PieADMM}                                                            \\
		\cmidrule(r){3-5} \cmidrule(r){6-8} \cmidrule(r){9-11}
		$\hat{n}$ &  $m$    & \multicolumn{1}{c}{Rel.Err.} & \multicolumn{1}{c}{NRMSE} & \multicolumn{1}{c}{Time\,(s)} & \multicolumn{1}{c}{Rel.Err.} & \multicolumn{1}{c}{NRMSE} & \multicolumn{1}{c}{Time\,(s)} & \multicolumn{1}{c}{Rel.Err.} & \multicolumn{1}{c}{NRMSE} & \multicolumn{1}{c}{Time\,(s)} \\
		\midrule[1pt]
		\multicolumn{11}{c}{$p_{cube}=0.3$} \\
		\midrule[0.5pt]
		2  & 12   & 0.0803                       & 0.1062                    & 0.11                        & 0.0803                       & 0.1062                    & 0.14                        & 0.0791                       & 0.1047                    & \underline{0.04}                        \\
		3  & 46   & 0.0854                       & 0.1046                    & 0.22                        & 0.0854                       & 0.1046                    & 0.32                        & 0.0848                       & 0.1039                    & \underline{0.09}                        \\
		4  & 109  & 0.1596                       & 0.1900                    & 0.42                        & 0.1596                       & 0.1900                    & 0.60                        & 0.1576                       & 0.1875                    & \underline{0.21}                        \\
		5  & 233  & 0.1969                       & 0.2309                    & 0.78                        & 0.1969                       & 0.2309                    & 1.19                        & 0.1955                       & 0.2292                    & \underline{0.44}                        \\
		6  & 413  & 0.2343                       & 0.2723                    & 1.03                        & 0.2343                       & 0.2722                    & 1.98                        & 0.2342                       & 0.2715                    & \underline{0.52}                        \\
		7  & 645  & 0.3083                       & 0.3560                    & 1.66                        & 0.3083                       & 0.3560                    & 3.23                        & 0.3074                       & 0.3546                    & \underline{1.05}                        \\
		8  & 1024 & 0.3264                       & 0.3752                    & 2.37                        & 0.3264                       & 0.3752                    & 4.82                        & 0.3225                       & 0.3708                    & \underline{1.66}                        \\
		9  & 1462 & 0.4312                       & 0.4940                    & 3.39                        & 0.4312                       & 0.4940                    & 6.84                        & 0.4309                       & 0.4890                    & \underline{1.96}                        \\
		10 & 2025 & 0.6305                       & 0.7205                    & 4.52                        & 0.6306                       & 0.7205                    & 9.48                        & 0.6237                       & 0.7127                    & \underline{2.22}          
		\\
		\midrule[1pt]
		\multicolumn{11}{c}{$p_{cube}=0.7$} \\
		\midrule[0.5pt]
		2  & 14   & 0.0540                       & 0.0714                    & 0.12                        & 0.0540                       & 0.0714                    & 0.16                        & 0.0539                       & 0.0713                    & \underline{0.05}                        \\
		3  & 69   & 0.0614                       & 0.0752                    & 0.31                        & 0.0614                       & 0.0752                    & 0.42                        & 0.0614                       & 0.0752                    & \underline{0.11}                        \\
		4  & 180  & 0.1338                       & 0.1593                    & 0.63                        & 0.1338                       & 0.1593                    & 0.92                        & 0.1262                       & 0.1503                    & \underline{0.22}                        \\
		5  & 362  & 0.1244                       & 0.1458                    & 0.89                        & 0.1244                       & 0.1458                    & 1.70                        & 0.1177                       & 0.1380                    & \underline{0.46}                        \\
		6  & 647  & 0.1542                       & 0.1792                    & 1.47                        & 0.1542                       & 0.1792                    & 2.88                        & 0.1521                       & 0.1767                    & \underline{0.72}                        \\
		7  & 1096 & 0.2631                       & 0.3038                    & 2.43                        & 0.2631                       & 0.3038                    & 4.87                        & 0.2499                       & 0.2885                    & \underline{1.41}                        \\
		8  & 1650 & 0.3212                       & 0.3692                    & 3.61                        & 0.3212                       & 0.3692                    & 7.34                        & 0.3068                       & 0.3527                    & \underline{2.01}                        \\
		9  & 2418 & 0.2340                       & 0.2681                    & 5.12                        & 0.2340                       & 0.2681                    & 10.68                       & 0.2236                       & 0.2562                    & \underline{2.58}                        \\
		10 & 3347 & 0.4295                       & 0.4908                    & 6.96                        & 0.4295                       & 0.4908                    & 14.49                       & 0.4268                       & 0.4876                    & \underline{3.11}   \\                    
		\bottomrule[1pt]             
	\end{tabular}}
	\label{cube_table}
\end{table}

\begin{figure}[!t]
	\centering
	\subfloat[number of edges]{\label{cube_mn}	\includegraphics[width=0.33\linewidth]{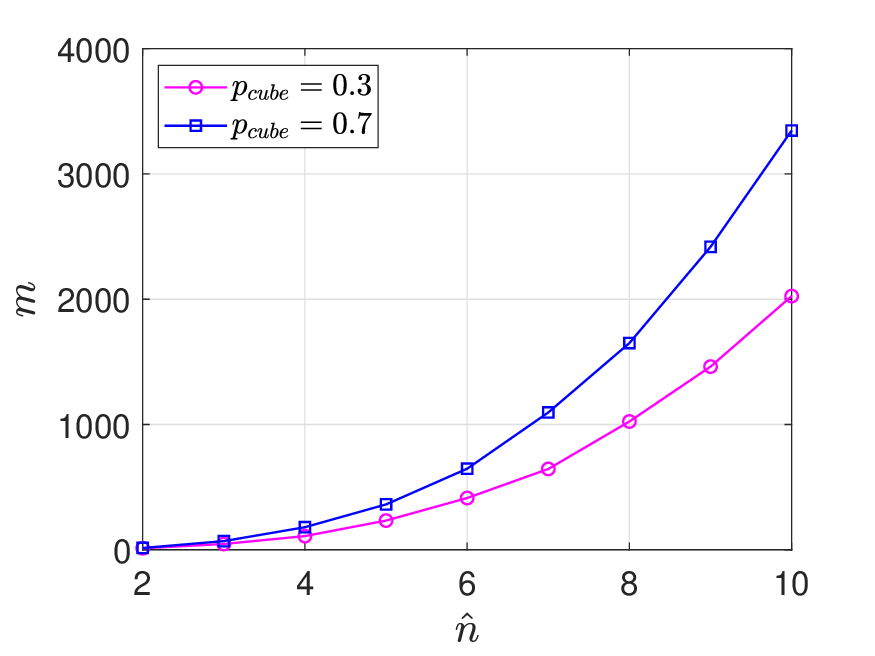}}
	\subfloat[$p_{cube}=0.3$]{\label{cube_n1}	\includegraphics[width=0.33\linewidth]{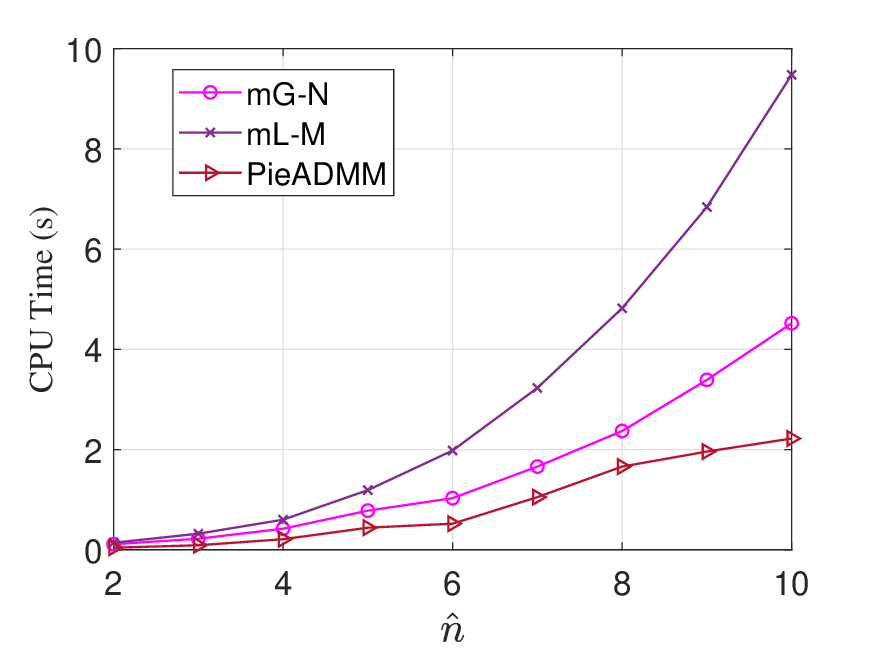}}
	\subfloat[$p_{cube}=0.7$]{\label{cube_n2}	\includegraphics[width=0.33\linewidth]{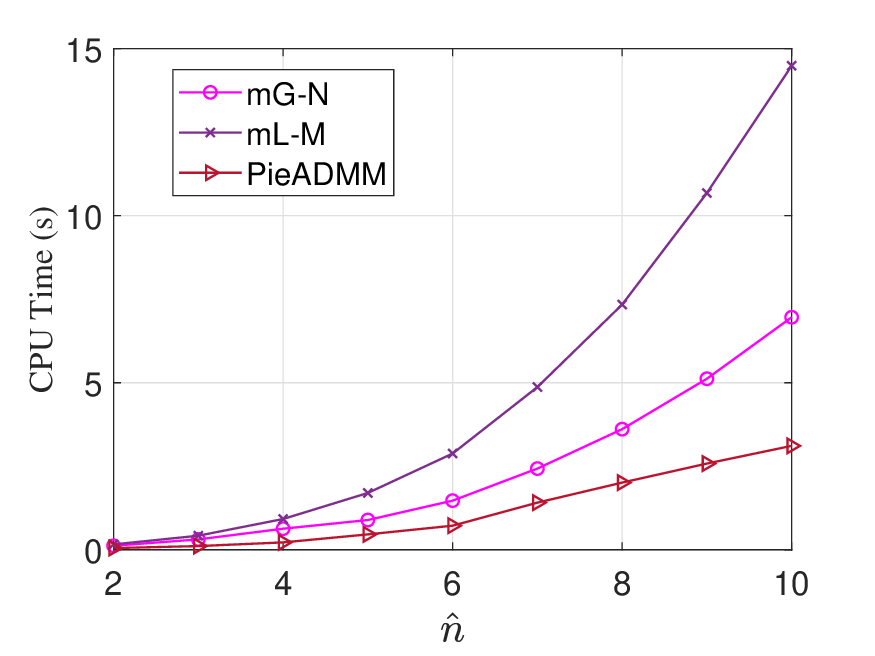}}
	\caption{The trend of the numbers of edges and CPU time along with $\hat{n}$ with $\sigma_{r} = 0.1,~\sigma_{t}^{rel}=0.1$ and $p_{cube}=0.3$ or $0.7$.}
	\label{cube_n}
\end{figure}

We also choose $\hat{n}$ from $2$ to $10$ and show the numerical results   in Table \ref{cube_table}. Fig.~\ref{cube_mn} indicates the relationship between the number of edges and vertices of the cube datasets, and Fig.~\ref{cube_n1} and \ref{cube_n2} illustrate the the upward trend of speed along with $\hat{n}$. The growth of cost of mG-N and mL-M are both cubic, and the growth of PieADMM is slower.

\subsection{SLAM benchmark datasets}

\begin{figure*}[!ht]
	\centering
	\subfloat{\label{garage_trajectory_niose}	\includegraphics[width=0.25\linewidth]{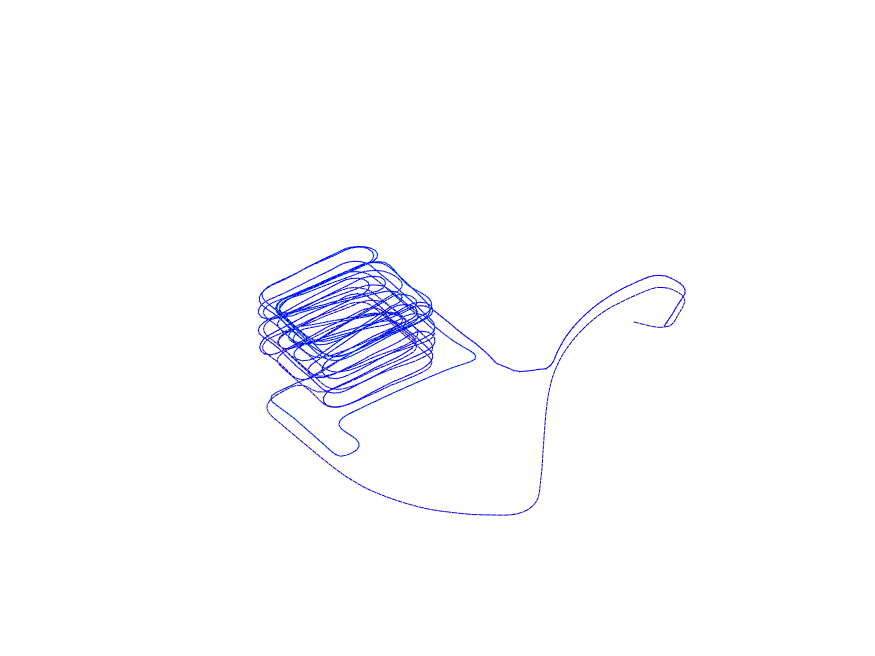}}
	\subfloat{\label{garage_trajectory_gn}	\includegraphics[width=0.25\linewidth]{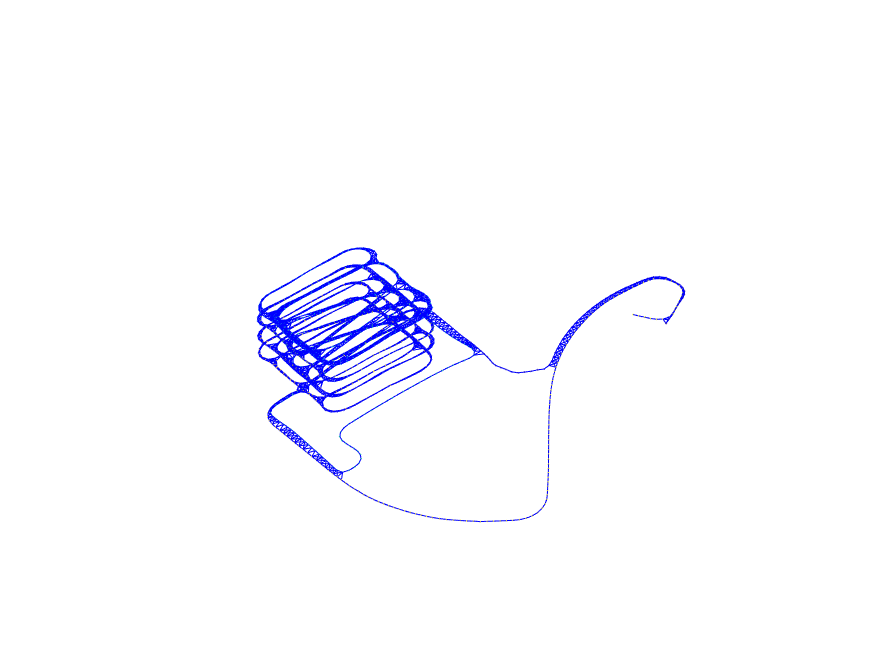}}
	\subfloat{\label{garage_trajectory_lm}	\includegraphics[width=0.25\linewidth]{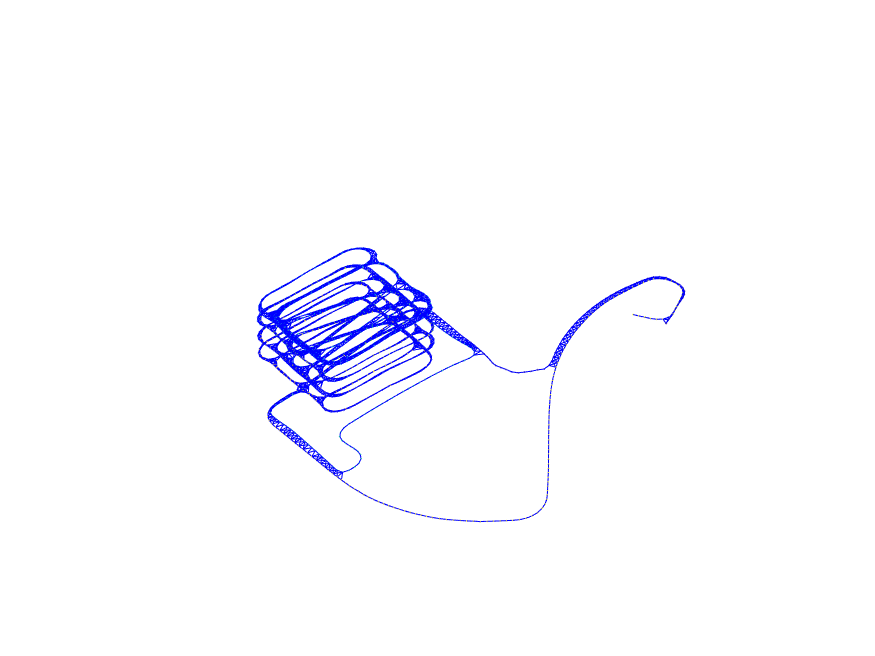}}
	\subfloat{\label{garage_trajectory_admm}	\includegraphics[width=0.25\linewidth]{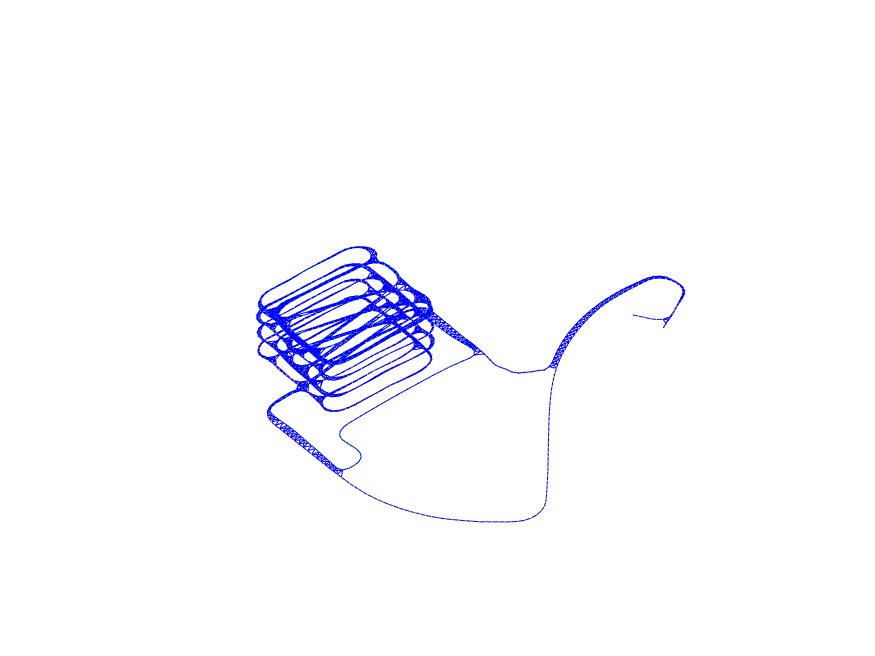}}
	
	\subfloat{\label{torus_trajectory_niose}	\includegraphics[width=0.25\linewidth]{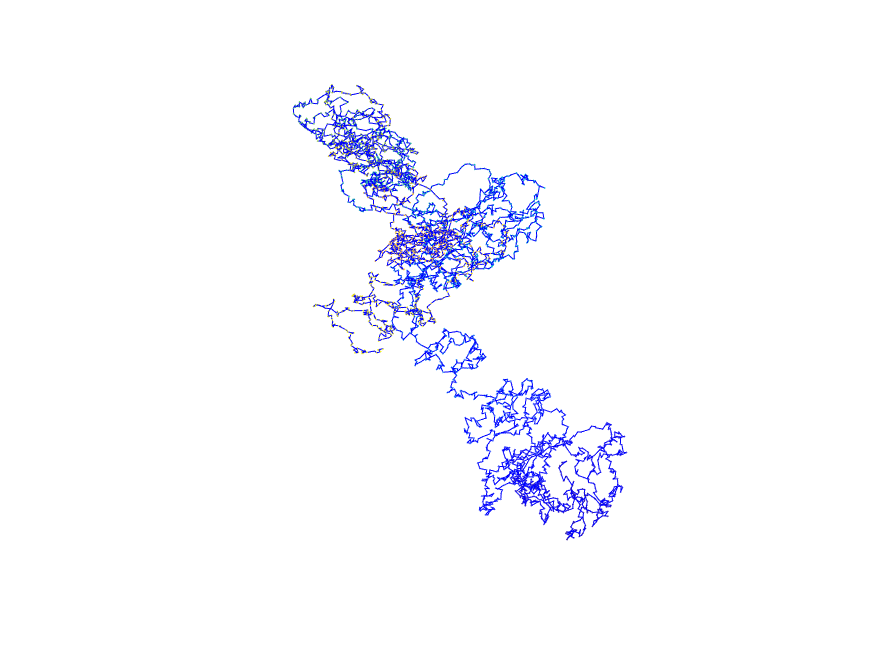}}
	\subfloat{\label{torus_trajectory_gn}	\includegraphics[width=0.25\linewidth]{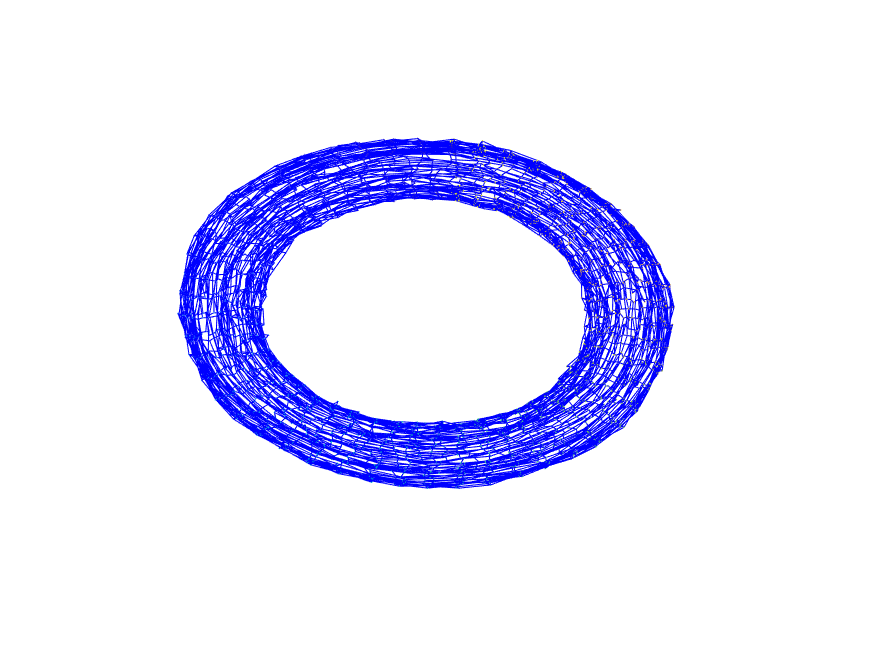}}
	\subfloat{\label{torus_trajectory_lm}	\includegraphics[width=0.25\linewidth]{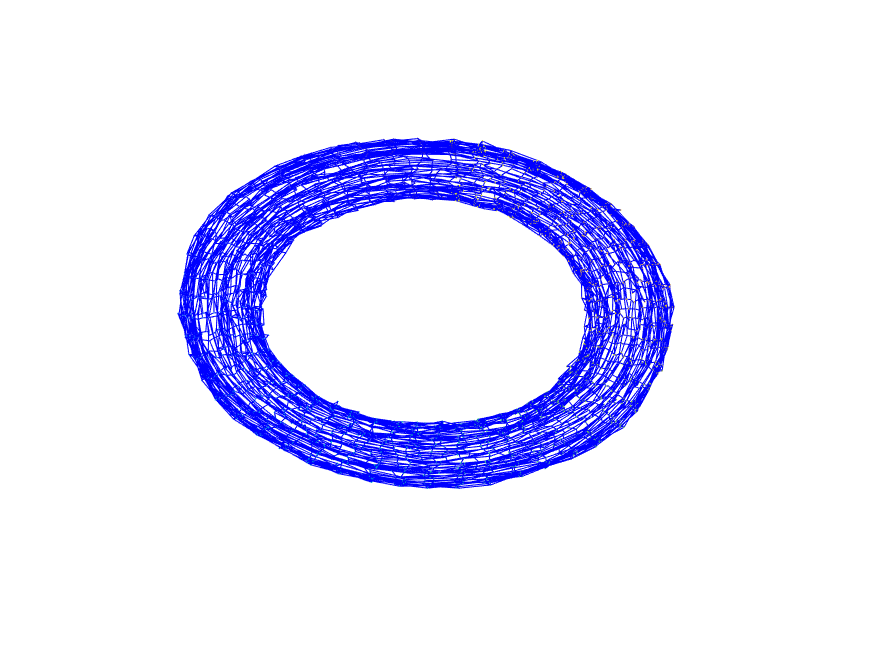}}
	\subfloat{\label{torus_trajectory_admm}	\includegraphics[width=0.25\linewidth]{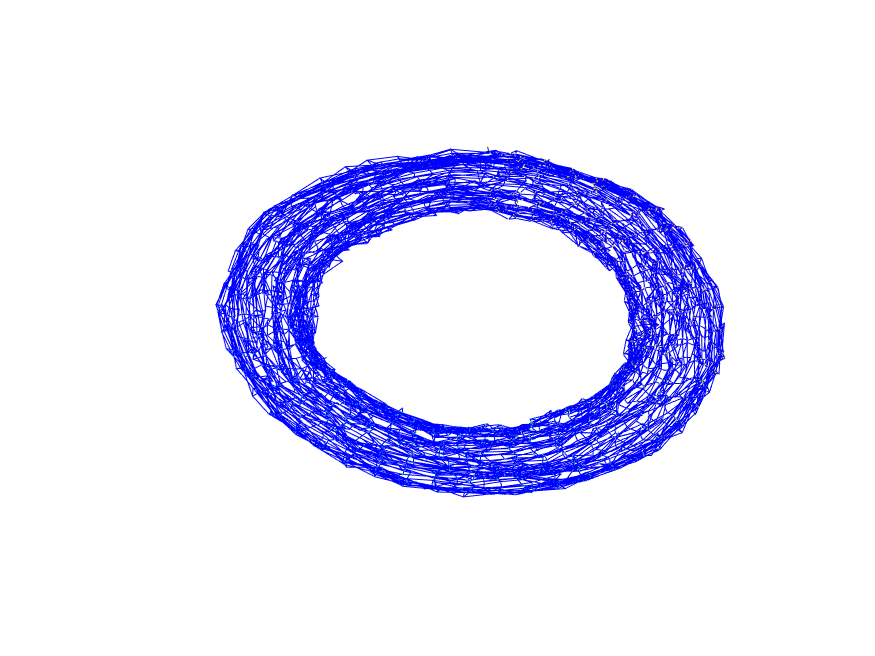}}
	
	\subfloat{\label{sphere2500_trajectory_niose}		
		\includegraphics[width=0.25\linewidth]{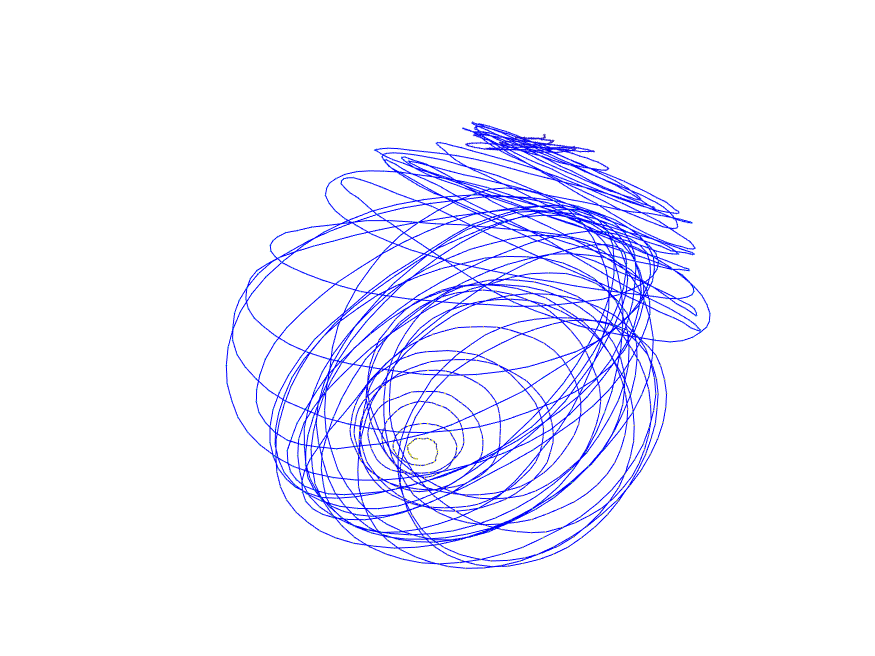}}
	\subfloat{\label{sphere2500_trajectory_gn}		
		\includegraphics[width=0.25\linewidth]{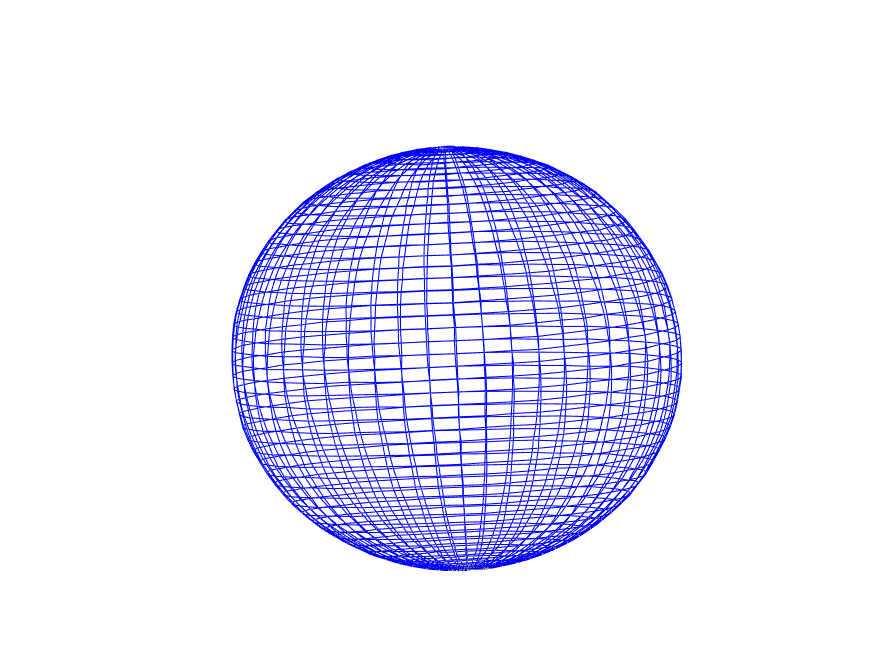}}
	\subfloat{\label{sphere2500_trajectory_lm}		
		\includegraphics[width=0.25\linewidth]{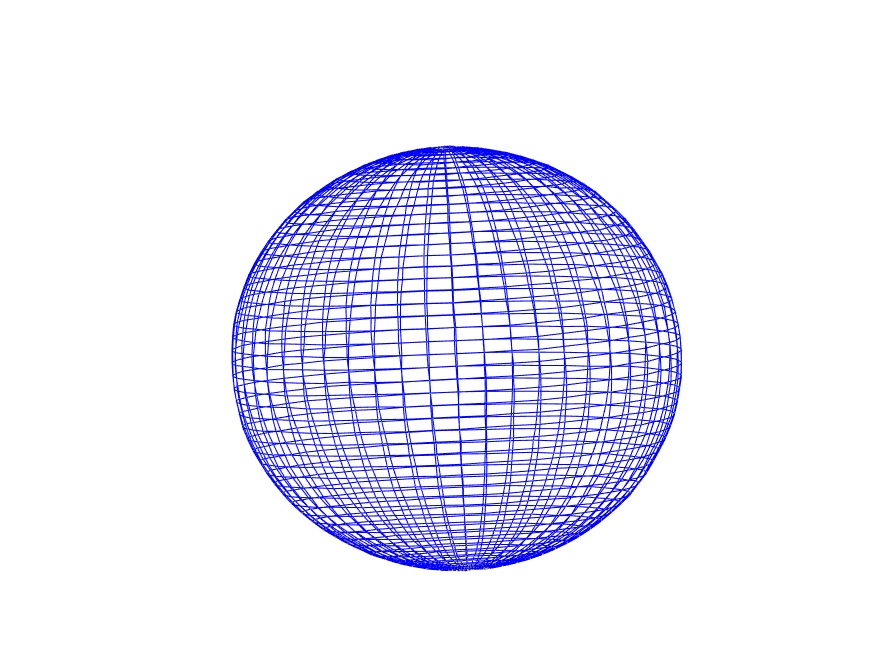}}
	\subfloat{\label{sphere2500_trajectory_admm}		
		\includegraphics[width=0.25\linewidth]{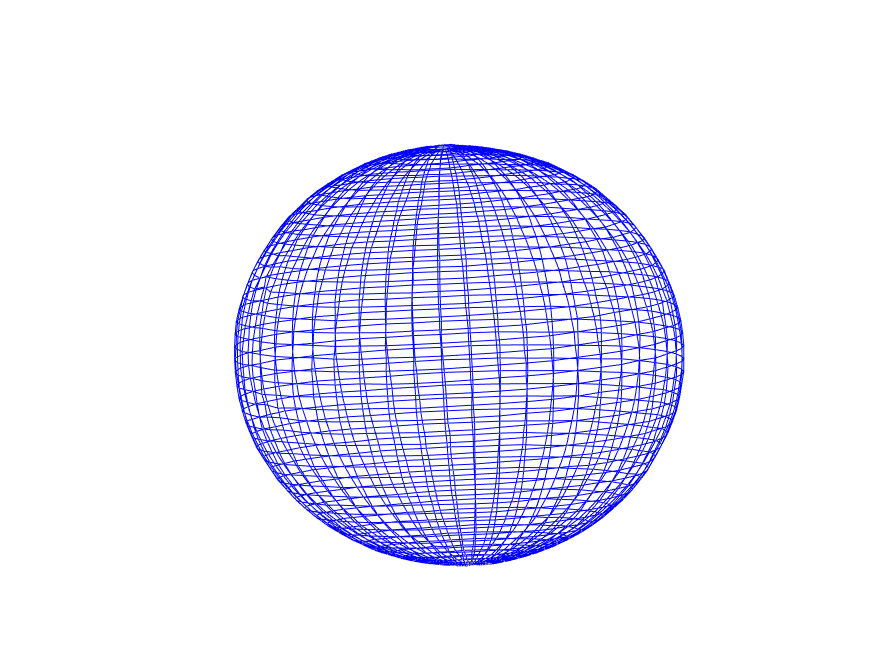}}
	
	\subfloat{\label{spherebig_trajectory_noise}		
		\includegraphics[width=0.25\linewidth]{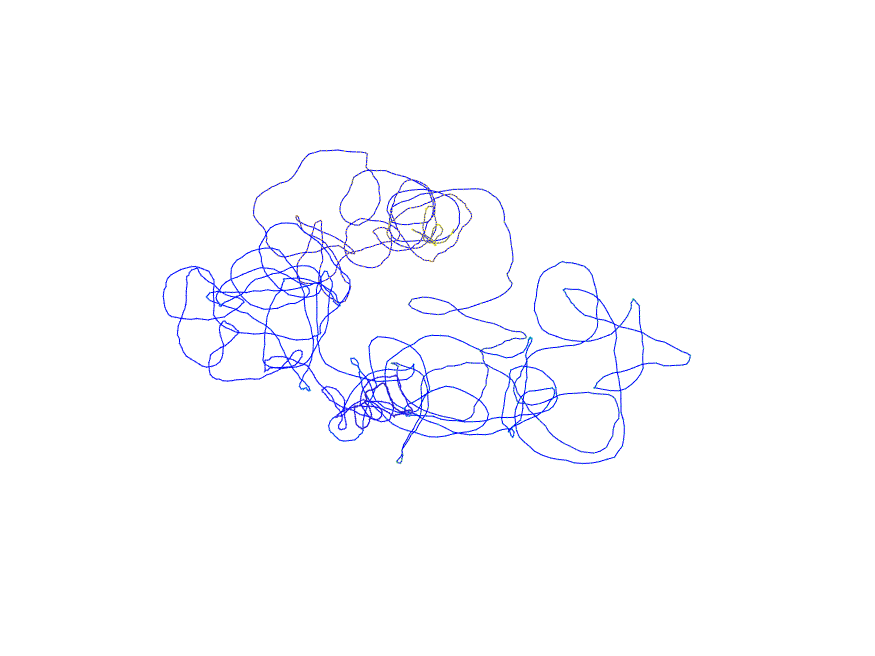}}
	\subfloat{\label{spherebig_trajectory_gn}		
		\includegraphics[width=0.25\linewidth]{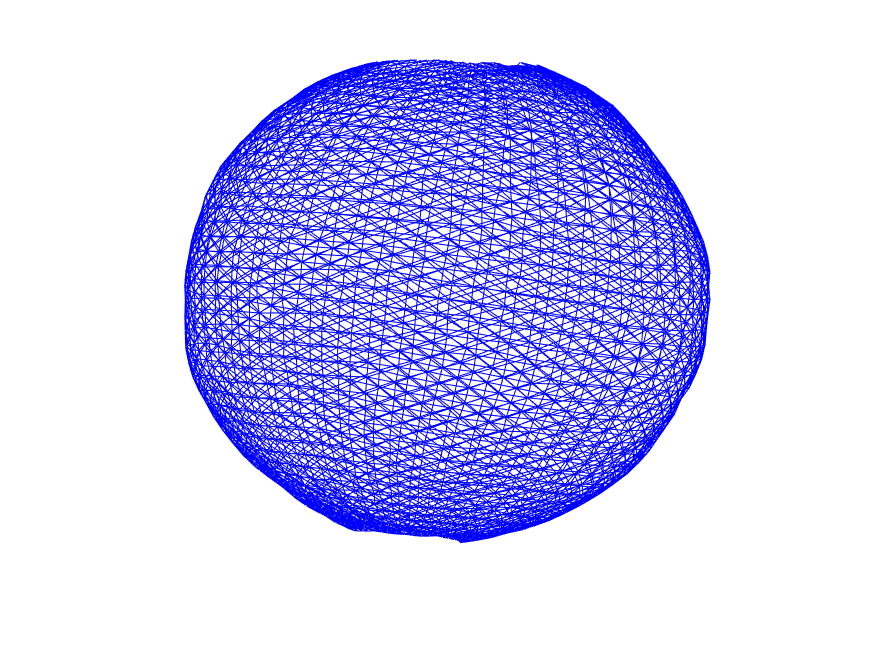}}
	\subfloat{\label{spherebig_trajectory_lm}	
		\includegraphics[width=0.25\linewidth]{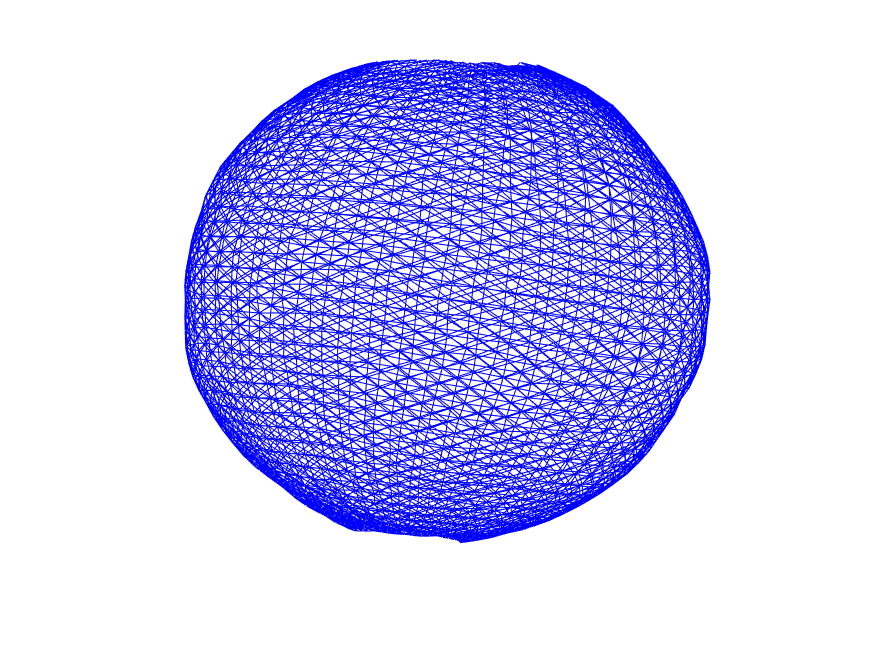}}
	\subfloat{\label{spherebig_trajectory_admm}		
		\includegraphics[width=0.25\linewidth]{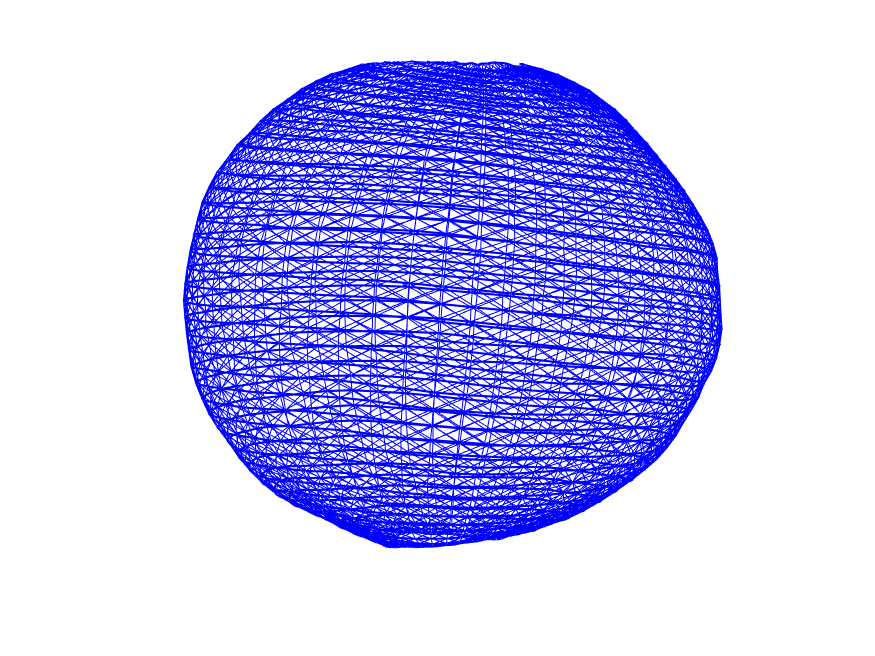}}
	
	\qquad Noised data  \qquad\qquad\qquad\quad mG-N \qquad\qquad\qquad\qquad mL-M \qquad\qquad\qquad\quad PieADMM
	\caption{The results of SLAM benchmark datasets in visual. The rows are the different data sets of garage, torus, sphere $1$ and sphere $2$, respectively. 
		From  the left column to  the right column are the corrupted data and the recovered results by mG-N, mL-M, and PieADMM, respectively.}
	\label{benchmark trajectory}
\end{figure*}

In this section, we test some popular 3D SLAM datasets. The garage dataset is a large-scale real-world example, and the other three (sphere $1$, sphere $2$ and torus) are common datasets used to compare performance. Different from sphere1 dataset, the larger noise is added to sphere2 dataset. We also use chordal initialization technique to compute an initial point for all methods. Fig.~\ref{benchmark trajectory} shows the results of trajectory in visual, and the corresponding numerical results are listed in Table \ref{benchmark_table}.
It is worth noting that our model of rotation is based on the vMF distribution rather than the traditional Gaussian distribution, so the recovered solution is not the same, and it does not make sense to compare objective function values or gradients. We show the CPU time in the table, which indicates that  PieADMM converges faster than mG-N and mL-M.

\begin{table}[t]
	\centering
	\renewcommand{\arraystretch}{1.3}
	\caption{The numerical results of SLAM benchmark datasets.}
	\begin{tabular}{cccccc}
		\toprule[1pt]
		&      &      & \multicolumn{3}{c}{CPU Time\,(s)}\\
		\cmidrule(r){4-6}
		datasets                  & $n$    & $m$    & mG-N        & mL-M        & PieADMM  \\
		\midrule[1pt]
		garage                    & 1661 & 6275 & 17.31       & 32.24       & \underline{14.37}       \\
		sphere $1$                   & 2500 & 4949 & 14.14       & 24.07       & \underline{11.05}       \\
		sphere $2$                   & 2200 & 8647 & 24.08       & 38.02       & \underline{17.95}       \\
		torus                     & 5000 & 9048 & 26.17       & 49.21       & \underline{21.81}      \\
		\bottomrule[1pt]
	\end{tabular}
	\label{benchmark_table}
\end{table}

\section{Conclusions} \label{conclusions}
Pose graph optimization in SLAM is a special non-convex optimization. 
In this paper, we propose a new non-convex pose graph optimization model based on augmented unit quaternion and von Mises-Fisher distribution, which is a large-scale quartic polynomial optimization on unit spheres. By introducing auxiliary variables, we reformulated it into a multi-quadratic polynomial optimization, multi-linear least square problem. Then we introduced a proximal linearized Riemannian ADMM for PGO model, in which the subproblems are simple projection problems, and can be solved in parallel corresponding to the structure of the directed graph, which greatly improve the efficiency. Then, based on the Lipschitz gradient continuity assumption which our PGO model satisfies and the first-order optimality conditions on manifolds, we establish the iteration complexity of $O(1/\epsilon^{2})$ for  finding an $\epsilon$-stationary solution. The numerical experiments on two synthetic datasets with different data scales and noise level and four 3D SLAM benchmark datasets verify the effectiveness of our method.

\bigskip 
{\bf Acknowledgment}  This work was supported by   the R\&D project of Pazhou Lab (Huangpu) (Grant no. 2023K0603),  the National Natural Science Foundation of China (No. 12131004), and the Fundamental Research Funds for the Central Universities (Grant No. YWF-22-T-204). 

{{\bf Data availability} Data will be made available on reasonable request.

	{\bf Conflict of interest} The authors declare no Conflict of interest.}

\bibliographystyle{IEEEtran}
\bibliography{ref}

\end{document}